\documentclass{amsart}

\usepackage{latexsym}

\usepackage{amsfonts,amssymb,euscript,epsfig}

\usepackage[all]{xy}

\newtheorem{theo}{Theorem}[section]
\newtheorem{defi}[theo]{Definition}
\newtheorem{assum}[theo]{Assumption}
\newtheorem{lem}[theo]{Lemma} 
\newtheorem{prop}[theo]{Proposition}
\newtheorem{defiprop}[theo]{Definition-Proposition}
\newtheorem{rem}[theo]{Remark}
\newtheorem{coro}[theo]{Corollary}
\newtheorem{conj}[theo]{Conjecture}

\newcommand{\agot}{\ensuremath{\mathfrak{a}}}
\newcommand{\ggot}{\ensuremath{\mathfrak{g}}}

\newcommand{\kgot}{\ensuremath{\mathfrak{k}}}

\newcommand{\rgot}{\ensuremath{\mathfrak{r}}}
\newcommand{\ngot}{\ensuremath{\mathfrak{n}}}
\newcommand{\pgot}{\ensuremath{\mathfrak{p}}}
\newcommand{\tgot}{\ensuremath{\mathfrak{t}}}
\newcommand{\Rgot}{\ensuremath{\mathfrak{R}}}

\newcommand{\Acal}{\ensuremath{\mathcal{A}}}
\newcommand{\Bcal}{\ensuremath{\mathcal{B}}}

\newcommand{\Dcal}{\ensuremath{\mathcal{D}}}

\newcommand{\Ical}{\ensuremath{\mathcal{I}}}

\newcommand{\Hcal}{\ensuremath{\mathcal{H}}}
\newcommand{\Kcal}{\ensuremath{\mathcal{K}}}
\newcommand{\Lcal}{\ensuremath{\mathcal{L}}}
\newcommand{\Qcal}{\ensuremath{\mathcal{Q}}}

\newcommand{\Ncal}{\ensuremath{\mathcal{N}}}
\newcommand{\Ocal}{\ensuremath{\mathcal{O}}}
\newcommand{\Pcal}{\ensuremath{\mathcal{P}}}
\newcommand{\Rcal}{\ensuremath{\mathcal{R}}}
\newcommand{\Scal}{\ensuremath{\mathcal{S}}}
\newcommand{\Xcal}{\ensuremath{\mathcal{X}}}
\newcommand{\Ycal}{\ensuremath{\mathcal{Y}}}
\newcommand{\Zcal}{\ensuremath{\mathcal{Z}}}
\newcommand{\Ucal}{\ensuremath{\mathcal{U}}}
\newcommand{\Vcal}{\ensuremath{\mathcal{V}}}

\newcommand{\Z}{\ensuremath{\mathbb{Z}}}
\newcommand{\C}{\ensuremath{\mathbb{C}}}
\newcommand{\Lfibre}{\ensuremath{\mathbb{L}}}

\newcommand{\R}{\ensuremath{\mathbb{R}}}
\newcommand{\N}{\ensuremath{\mathbb{N}}}
\newcommand{\tore}{\ensuremath{\mathbb{T}}}

\newcommand{\esp}{\ensuremath{\varepsilon}}
\newcommand{\f}{\ensuremath{\mathcal{C}^{\infty}}}
\newcommand{\fgene}{\ensuremath{\mathcal{C}^{-\infty}}}
\newcommand{\croc}{\ensuremath{\hookrightarrow}}
\newcommand{\indice}{\ensuremath{\hbox{\rm Index}}}
\newcommand{\ch}{\ensuremath{\hbox{\rm H}}}

\newcommand{\mm}{\ensuremath{\hbox{\rm m}}}
\newcommand{\Cr}{\ensuremath{\hbox{\rm Cr}}}
\newcommand{\Vol}{\ensuremath{\hbox{\rm vol}}}
\newcommand{\T}{\ensuremath{\hbox{\bf T}}}
\newcommand{\tr}{\ensuremath{\hbox{\bf Tr}}}
\newcommand{\K}{\ensuremath{\hbox{\bf K}}}

\newcommand{\Char}{\ensuremath{\hbox{\rm Char}}}
\newcommand{\End}{\ensuremath{\hbox{\rm End}}}
\newcommand{\Thom}{\ensuremath{\hbox{\rm Thom}}}
\newcommand{\sthom}{\ensuremath{\hbox{\rm S-Thom}}}

\def \spin {{\rm Spin}}
\def \spinc {{\rm Spin}^{c}}
\def \so {{\rm SO}}
\def \u {{\rm U}}
\def \Det {{\rm Det}}
\def \Pspin {{\rm P}_{{\rm Spin}^{\rm c}}}
\def \Pso {{\rm P}_{{\rm SO}}}
\def \Pu {{\rm P}_{{\rm U}}}
\def \clif {{\bf c}}
\def \Clif {{\rm Cl}}

\def \indT {{\rm Ind}^{^K}_{_T}}

\def \HolT {{\rm Hol}^{^K}_{_T}}
\def \HolB {{\rm Hol}^{^K}_{_{K_{\beta}}}}
\def \HolTB {{\rm Hol}^{^{K_{\beta}}}_{_T}}

\begin{document}
    
\title{Spin$^{\rm c}$-quantization and the $K$-multiplicities of the discrete series}

\author{Paul-Emile  PARADAN}

\maketitle
{\center
UMR 5582, Institut Fourier, B.P. 74, 38402, Saint-Martin-d'H\`eres 
cedex, France\\
e-mail: Paul-Emile.Paradan@ujf-grenoble.fr\\
}

{\center November 2001\\
}

\vspace{7mm}

\begin{abstract}
  We express the $K$-multiplicities of a representation of  
  the discrete series associated to a coadjoint orbit $\Ocal$ in terms 
  of $\spinc$-index on symplectic reductions of $\Ocal$.
\end{abstract}

{\def\thefootnote{\relax}
\addtocounter{footnote}{-1}}

{\small
\tableofcontents
}


\section{Introduction and statement of the results}

The purpose of this paper is to show that the `{\em 
quantization commutes with reduction}' principle of  Guillemin-Sternberg  
\cite{Guillemin-Sternberg82} holds for the  coadjoint orbits that 
parametrize the discrete series of a real connected 
semi-simple Lie group.

\subsection{Discrete series and $K$-multiplicities}\label {intro-1}
Let $G$ be a connected, real, semisimple Lie group with finite center.  
By definition, the {\em discrete series} of $G$ is the set of isomorphism 
classes of irreducible, square integrable, unitary 
representations of $G$. Let $K$ be a maximal compact subgroup of $G$, and $T$ be a 
maximal torus in $K$. Harish-Chandra has shown that $G$ has a discrete
series if and only if $T$ is a Cartan subgroup of $G$ 
\cite{Harish-Chandra65-66}. For the remainder of this paper, we 
may therefore assume that $T$ is a Cartan subgroup of $G$. 

Let us fix some notation. We denote by $\ggot,\kgot,\tgot$ the Lie algebras of 
$G,K,T$, and by  $\ggot^*,\kgot^*,\tgot^*$ their duals. Let
$\Lambda^*\subset\tgot^*$ be the set of real weights: 
$\alpha\in\Lambda^*$ if $i\alpha$ is the differential of a character of $T$.
Let $\Rgot_{c}\subset\Rgot\subset\Lambda^{*}$ be respectively the set of 
roots for the action of $T$ on $\kgot\otimes\C$ and $\ggot\otimes\C$.
We choose a system of positive roots $\Rgot_{c}^{+}$ for 
$\Rgot_{c}$. We denote by $\tgot^{*}_{+}$ the corresponding Weyl chamber, 
and we let $\rho_{c}$ be half the sum of the elements of $\Rgot_{c}^{+}$. 
The set $\Lambda^*_+:=\Lambda^*\cap\tgot^*_+$ parametrizes the unitary dual 
of $K$. For $\mu\in\Lambda^*_+$, let $\chi_{_{\mu}}^{_K}$ be the  
character of the irreducible $K$-representation with highest weight $\mu$.

Harish-Chandra parametrizes the discrete series by a 
discrete subset $\widehat{G}_{d}$ of regular elements of 
the Weyl chamber $\tgot^*_+$ \cite{Harish-Chandra65-66}.  
He associates to any $\lambda\in \widehat{G}_{d}$ an invariant 
eigendistribution on $G$, denoted by $\Theta_{\lambda}$, which is shown 
to be the global trace of an irreducible, square integrable, unitary 
representation $\Hcal_{\lambda}$ of $G$. It is a generalized function 
on $G$, invariant by conjugation, which  admits a restriction to 
$K$ denoted by $\Theta_{\lambda}\vert_{K}$. The distribution 
$\Theta_{\lambda}\vert_{K}$ corresponds to the global trace of
the induced representation of $K$ on $\Hcal_{\lambda}$. It 
admits a decomposition
$$
\Theta_{\lambda}\vert_{K}=\sum_{\mu\in\Lambda^*_+}\mm_{\mu}(\lambda)\, 
\chi_{_{\mu}}^{_K}\ ,
$$
where the integers $\mm_{\mu}(\lambda)$ satisfy certain
 combinatorial identities called the  Blattner 
formulas \cite{Hecht-Schmid}.

\medskip

{\em The  main goal of the paper is to relate the multiplicities $\mm_{\mu}(\lambda)$ 
to the geometry of the coadjoint orbit $G\cdot\lambda\subset\ggot^*$ 
as predicted by the Guillemin-Sternberg principle evoked above.} 

\medskip

Before stating our result we recall how a representation belonging
to the discrete series can be realized as the quantization of
a coadjoint orbit.
 
\subsection{Realisation of the discrete series}\label{intro-2}

In the 60's, Kostant and Langlands conjectured realisations of the
discrete series in terms of ${\rm L}^2$ cohomology 
that fit into the general framework of quantization. 
The proof of this conjecture was given by Schmid somes years later
\cite{Schmid71,Schmid76}. Let us recall the procedure for a 
fixed $\lambda\in\widehat{G}_{d}$.

The manifold  $G\cdot\lambda$ carries several $G$-invariant complex 
structures. For convenience we work with the complex structure $J$ defined by the 
following condition: each weight $\alpha$ for the  $T$-action on the 
tangent space $(\T_{\lambda}(G\cdot\lambda),\, J)$ satisfies 
$(\alpha,\lambda)>0$.

Let  $\Rgot^+\subset\Rgot$ be the set of positive roots defined by $\lambda$: 
$\alpha\in\Rgot^+\Longleftrightarrow (\alpha,\lambda)>0$. Let $\rho$ be 
half the sum of the elements of $\Rgot^{+}$. The condition 
$\lambda\in\widehat{G}_d$ imposes that $\lambda-\rho$ is a weight for
$T$, so we can consider the line bundle 
$$
\tilde{L}:=G\times_{T}\C_{\lambda-\rho}
$$ 
over $G\cdot\lambda\simeq G/T$: this
line bundle carries a canonical holomorphic structure. Let
$\Omega^k(\tilde{L})$ be the space of $\tilde{L}$-valued $(0,k)$ forms on 
$G\cdot\lambda$, and $\overline{\partial}_{\tilde{L}}:\Omega^k(\tilde{L})\to 
\Omega^{k+1}(\tilde{L})$ be the Dolbeault operator. The choice of
$G$-invariant hermitian metrics on $G\cdot\lambda$ and on $\tilde{L}$ 
give meaning to the formal adjoint $\overline{\partial}^*_{\tilde{L}}$ of the 
$\overline{\partial}_{\tilde{L}}$ operator, and to the Dolbeault-Dirac operator
$ \overline{\partial}_{\tilde{L}}+ \overline{\partial}^*_{\tilde{L}}$. 

The ${\rm L}^2$ cohomology of $\tilde{L}$, which we denote by 
$\ch^*_{(2)}(G\cdot\lambda,\tilde{L})$,  is equal to the kernel
of the differential operator $\overline{\partial}_{\tilde{L}}+ 
\overline{\partial}^*_{\tilde{L}}$ acting on the subspace of 
$\Omega^*(\tilde{L})$ formed by the square integrable elements.

\medskip

\begin{theo}(Schmid). \label{theo.schmid}

\ Let $\lambda\in\widehat{G}_{d}$.

(i) $\ch^k_{(2)}(G\cdot\lambda,\tilde{L})=0$ if $k\neq
\frac{\dim(G/K)}{2}$.

(ii) If $k=\frac{\dim(G/K)}{2}$, then
$\ch^k_{(2)}(G\cdot\lambda,\tilde{L})$ is the irreducible
representation $\Hcal_{\lambda}$.
\end{theo}

\medskip

So, the representation $\Hcal_{\lambda}$ is  the 
quantization of the coadjoint orbit $G\cdot\lambda$ beeing the index  
of the Dolbeault-Dirac operator $\overline{\partial}_{\tilde{L}}+ 
\overline{\partial}^*_{\tilde{L}}$ (in the ${\rm L}^2$ sense and
modulo $(-1)^{\frac{\dim(G/K)}{2}}$). In the next subsection, we briefly recall
the `quantization commutes with reduction'  principle of 
Guillemin-Sternberg, and in subsection \ref{intro-4} we state our main 
result.

\subsection{Quantization commutes with reduction}\label{intro-3}

Let $M$ be a Hamiltonian $K$-manifold with symplectic form $\omega$ 
and moment map $\Phi: M\to \kgot^{*}$. The coadjoint orbits
$G\cdot\lambda$ introduced earlier are the key
examples here. Each is equipped with its Kirillov-Kostant-Souriau 
symplectic form $\omega$, and the action of $G$ is Hamiltonian with moment map 
$G\cdot\lambda\croc \ggot^*$ equal to the inclusion. Let $K$ be the
maximal compact Lie subgroup of $G$ introduced in subsection
\ref{intro-2}. The induced action of $K$ on $G\cdot\lambda$ is
Hamiltonian, and the corresponding moment map
$\Phi:G\cdot\lambda\to\kgot^*$ is equal to the composition of the inclusion 
$G\cdot\lambda\croc \ggot^*$ with the projection $\ggot^*\to\kgot^*$.

In the process of quantization one tries to associate a unitary
representation of $K$ to  the data $(M,\omega,\Phi)$. In this general
framework, when $M$ is {\em compact} and under certain
integrability conditions, we associate to these data a virtual representation of $K$ 
defined as the equivariant index of a $\spinc$ Dirac operator: it's
the $\spinc$ quantization. We need two auxilliary data :

 (i)  A  {\em prequantum line bundle} $L\to M$: it is a $K$-equivariant
 Hermitian line bundle equipped with $K$-invariant connection whose
 curvature form is $-i\, \omega$.

 (ii)  A $K$-invariant almost complex structure $J$ on $M$,
{\em compatible} with the symplectic structure: $(v,w)\mapsto \omega(v,Jw)$
defines a metric.

One considers then the $K$-equivariant $\spinc$ Dirac operator $D_L$
corresponding to the $\spinc$ structure on $M$ defined by $J$, and
twisted by the line bundle $L$ \cite{Lawson-Michel,Duistermaat96}. The $\spinc$-quantization of 
$(M,\omega,\Phi)$ is the equivariant index of the differential
operator $\Dcal_L$
$$
RR^{^K}(M,L) := \indice^K_M(\Dcal_L)\quad \in R(K)\ , 
$$
where $R(K)$ is the representation ring of $K$. When $K$ is reduced to 
$\{e\}$, the $\spinc$-quantization of $(M,\omega)$ is just an integer: 
$RR(M,L)\in \Z$.

A fundamental result of Marsden-Weinstein asserts that if
$\xi\in\kgot^{*}$ is a  regular value of the moment map $\Phi$, the 
{\em reduced space} 
$$
M_{\xi}:=\Phi^{-1}(\xi)/K_{\xi}\cong
\Phi^{-1}(K\cdot\xi)/K
$$ 
is an orbifold equipped with a symplectic structure $\omega_{\xi}$
(which one calls also symplectic quotient). For any dominant weight 
$\mu\in \Lambda^{*}_{+}$ which is a regular value of $\Phi$, 
$$
L_{\mu}:=(L\vert_{\Phi^{-1}(\mu)}\otimes \C_{-\mu})/K_{\mu}
$$
is a prequantum orbifold-line bundle over $(M_{\mu},\omega_{\mu})$. The 
definition of $\spinc$-index carries over to the orbifold case, 
hence $RR(M_{\mu},L_{\mu})\in \Z$ is defined. In 
\cite{Meinrenken-Sjamaar}, this is extended further to the case of 
singular symplectic quotients, using partial (or shift) 
desingularization. So the integer $RR(M_{\mu},L_{\mu})\in \Z$ is well
defined for every $\mu\in \Lambda^{*}_{+}$.

The following Theorem was conjectured by 
Guillemin-Sternberg \cite{Guillemin-Sternberg82} and is known as 
``quantization commutes with reduction'' \cite{Meinrenken98,Meinrenken-Sjamaar}. 

\medskip 

\begin{theo} (Meinrenken, Meinrenken-Sjamaar). \label{th.Q-R}
Let $(M,\omega,\Phi)$ be a compact Hamiltonian $K$-manifold prequantized 
by $L$. Let $RR^{^K}(M,-)$ be the equivariant Riemann-Roch 
character defined by means of a compatible almost complex structure on $M$. 
We have the following equality in $R(K)$
$$
RR^{^K}(M,L)=\sum_{\mu\in 
\Lambda^{*}_{+}}RR(M_{\mu},L_{\mu})\, 
\chi_{_{\mu}}^{_K}\ .
$$
\end{theo}

\begin{rem} For a  compact Hamiltonian $K$-manifold $(M,\omega,\Phi)$, 
the Convexity Theorem \cite{Kirwan.84.bis} asserts that 
$\Delta:=\Phi(M)\cap\tgot^*_+$ is a convex rational polytope. 
In Theorem \ref{th.Q-R},  we have  $RR(M_{\mu},L_{\mu})=0$ if 
$\mu\notin\Delta$.
\end{rem}

Other proofs can be found in \cite{pep4,Tian-Zhang98}. For an 
introduction and further references see \cite{Sjamaar96,Vergne01}.

A natural question is to extend Theorem \ref{th.Q-R} to the 
{\em non-compact} Hamiltonian $K$-manifolds which admit a {\em proper}
moment map. In this situation, the reduced space
$M_{\xi}:=\Phi^{-1}(\xi)/K_{\xi}$ is compact for every
$\xi\in\kgot^*$, so the integer $RR(M_{\mu},L_{\mu})\in \Z,\,\mu\in 
\Lambda^{*}_{+},$ can be defined like before. 
\begin{conj}\label{conjecture}
Let $(M,\omega,\Phi)$ be a Hamiltonian $K$-manifold with {\em proper}
moment map, and prequantized by $L$. Let $\overline{\partial}_{L}+ 
\overline{\partial}^*_{L}$ be the Dolbeault-Dirac operator defined by
means of a $K$-invariant compatible almost complex structure, and 
$K$-invariant metric on $M$ and $L$. Then
$$
{\rm L}^2-{\rm Index}^{_K}\Big(\overline{\partial}_{L}+ 
\overline{\partial}^*_{L}\Big)=\sum_{\mu\in 
\Lambda^{*}_{+}}RR(M_{\mu},L_{\mu})\, 
\chi_{_{\mu}}^{_K}\ .
$$
\end{conj}

We present in the next subsection the central result of this paper
that shows that Conjecture \ref{conjecture}
is true, apart from a $\rho_c$-shift,
 for the coadjoint orbits that parametrize the discrete series.

\subsection{The results}\label{intro-4}

Consider the Hamiltonian action of $K$ on the coadjoint orbit  $G\cdot\lambda$.  
Since $G\cdot\lambda$ is closed in $\ggot^*$, the moment map
$\Phi:G\cdot\lambda\to\kgot^*$ is {\em proper} \cite{pep3}. 
Our main Theorem can be stated roughly as follows.
\begin{theo} \label{theo-principal}
Let $\mm_{\mu}(\lambda), \mu\in\Lambda^*_+$, be the 
$K$-multiplicities of the representation $\Hcal_{\lambda}\vert_K$.
For  $\mu\in\Lambda^*_+$ we have:

(i) If $\mu +\rho_c$ is a regular value of
$\Phi$, the orbifold $(G\cdot\lambda)_{\mu +\rho_c}:=\Phi^{-1}(\mu
+\rho_c)/T$, oriented by its symplectic form $\omega_{\mu +\rho_c}$, 
carries a $\spinc$ structure such that
$$
 \mm_{\mu}(\lambda)= \Qcal\left((G\cdot\lambda)_{\mu+\rho_c}\right)\ ,
$$
where the RHS is the index of the corresponding $\spinc$ 
Dirac operator on the reduced space $(G\cdot\lambda)_{\mu +\rho_c}$.

(ii) In general, one can define an integer
    $\Qcal\left((G\cdot\lambda)_{\mu+\rho_c}\right)\in \Z$, as the index of a $\spinc$ 
    Dirac operator on a reduced space $(G\cdot\lambda)_{\xi}$ where  
    $\xi$ is a regular value of $\Phi$, close enough to $\mu+\rho_c$. We still have
$\mm_{\mu}(\lambda)= \Qcal\left(M_{\mu+\rho_c}\right)$.
\end{theo}

Our Theorem states that the decomposition of
$\Theta_{\lambda}\vert_{K}$ into $K$-irreducible
components follows the philosophy of Guillemin-Sternberg:
$$
(1)\quad \quad 
\Theta_{\lambda}\vert_{K}=\sum_{\mu\in\Lambda^*_+}
\Qcal\left((G\cdot\lambda)_{\mu+\rho_c}\right)\,\chi_{_{\mu}}^{_K}\ .
$$
We know from Theorem \ref{theo.schmid} that $\Theta_{\lambda}\vert_{K}=
(-1)^{\frac{\dim(G/K)}{2}} {\rm L{\sp 2}-Index}^{_K}
(\overline{\partial}_{\tilde{L}}+ \overline{\partial}^*_{\tilde{L}})$,
hence  Theorem \ref{theo-principal} states also  that
$$
(2)\quad \quad 
{\rm L{\sp 2}-Index}^{_K}(\overline{\partial}_{\tilde{L}}+ 
\overline{\partial}^*_{\tilde{L}})=(-1)^{\frac{\dim(G/K)}{2}} 
\sum_{\mu\in\Lambda^*_+}
\Qcal\left((G\cdot\lambda)_{\mu+\rho_c}\right)\,
\chi_{_{\mu}}^{_K}\ .
$$

The main difference between Conjecture \ref{conjecture} and Theorem 
\ref{theo-principal} is the $\rho_c$-shift and the choices of 
$\spinc$ structure on the symplectic quotients $M_{\mu}$ 
and $(G\cdot\lambda)_{\mu+\rho_c}$. 

The $\rho_c$-shift is due to the fact that the line bundle $\tilde{L}$ is not a
prequantum line bundle over $(G\cdot\lambda,\omega)$. The difference
on the choice of $\spinc$ structure comes from the fact that the complex
structure $J$ on $G\cdot\lambda$ is not compatible with the symplectic
structure (unless $G=K$ is compact). Hence $J$ does not descend to 
the symplectic reductions $(G\cdot\lambda)_{\mu+\rho_c}$ in general:
the choice of the $\spinc$ structure on them need some care (see Propositions 
\ref{prop.spinc.induit.1} and \ref{prop.spinc.induit.2}).

\begin{rem} For a  Hamiltonian $K$-manifold $M$ with {\em proper} moment 
map $\Phi$, the Convexity Theorem \cite{Kirwan.84.bis,L-M-T-W,Sjamaar98} 
asserts that $\Delta:=\Phi(M)\cap\tgot^*_+$ is a convex rational polyhedron. 
In Theorem \ref{theo-principal}, we have $\Qcal((G\cdot\lambda)_{\mu+\rho_c})=0$ 
if $\mu+\rho_c$ does not belong to the {\em relative interior} of $\Delta$ 
(see Prop. \ref{prop.Q.mu.ro}).
\end{rem}

\subsection{Outline of the Proof}\label{intro-5}
We have  to face the following difficulties:

$[1]$ The symplectic manifold $G\cdot\lambda$ is not compact.

$[2]$ The complex structure on $G\cdot\lambda$ is not compatible with
the symplectic form $\omega$. In other words, the Kirillov-Kostant-Souriau 
symplectic form does not define a K\"ahler structure on $G\cdot\lambda$ 
unless $G=K$ is compact.

$[3]$ The line bundle $\tilde{L}$ is not a
    prequantum line bundle over $(G\cdot\lambda,\omega)$. It's what we call
    in the rest of this paper a $\kappa$-{\em prequantum}\footnote{Formally, 
    $\tilde{L}$ is the tensor product of a prequantum line bundle over 
    $(G\cdot\lambda,\omega)$  with a square root of $\kappa$.} line bundle over 
    $(G\cdot\lambda,\omega, J)$: if $\kappa$ denotes the canonical line bundle 
    of $(G\cdot\lambda,J)$, the tensor product $\tilde{L}^2\otimes \kappa^{-1}$ 
    is a prequantum line bundle over $(G\cdot\lambda,2\omega)$. 

\medskip

The first step of the proof is to solve the difficulties $[2]$ and $[3]$
in the {\em compact} situation. In Section \ref{sec.quantization.compact}, we give  
a modified version of Theorem \ref{th.Q-R} when $(M,\omega,\Phi)$ is
a {\em compact} Hamiltonian $K$-manifold which is equipped with an 
almost complex structure $J$ - not necessarily compatible with
$\omega$ - and a $\kappa$-prequantum line bundle $\tilde{L}$.

\begin{theo}\label{theo.egalite1}
Let $RR^{^K}(M,-)$ be the Riemann-Roch character 
defined by $J$. If the infinitesimal stabilizers for the 
action of $K$ on $M$ are {\em Abelian}, we have
$$
(3)\quad \quad 
RR^{^K}(M,\tilde{L})=\esp
\sum_{\mu\in \Lambda^{*}_{+}}Q(M_{\mu+\rho_{c}})\, 
\chi_{_{\mu}}^{_K}\ ,
$$
where  $\esp=\pm 1$ is the `quotient' of the orientations induced 
by the almost complex structure, and the  symplectic form. 
\end{theo}

In $(3)$, the integer $Q(M_{\mu+\rho_{c}})$ are computed like in Theorem
\ref{theo-principal} (see Def. \ref{prop.Q.mu.ro} for a more precise definition).

In the second step of the proof,  we extend $(3)$ to a 
non-compact setting. Instead of working with the $\rm L{\sp 2}$-Index, we define in 
Section \ref{sec.quant.non.compact} a {\em generalized Riemann-Roch character} 
$RR^{^K}_{\Phi}(M,-)$ when  $(M,\omega,\Phi)$ is a Hamiltonian 
manifold such that the function $\parallel\Phi\parallel^2:M\to\R$ has a 
{\em compact} set of critical points. For every $K$-vector bundle
$E\to M$, the distribution $RR^{^K}_{\Phi}(M,E)$ is defined as the 
index of a transversally elliptic operator on $M$ \cite{Atiyah74}. When the manifold 
is compact, the maps $RR^{^K}_{\Phi}(M,-)$ and $RR^{^K}(M,-)$ coincide.

We prove in Section \ref{section.Q.R} that Theorem \ref{theo.egalite1} 
generalizes to 
\begin{theo}\label{theo.egalite2}
Let $(M,\omega,\Phi)$ be a Hamiltonian $K$-manifold with {\em proper}
moment map and such that the function $\parallel\Phi\parallel^2:M\to\R$ has a 
{\em compact} set of critical points. If the infinitesimal stabilizers
are Abelian, and under Assumption \ref{hypothese.phi.t.carre}, we have

$$
(4)\quad\quad RR^{^K}_{\Phi}(M,\tilde{L})=\esp
\sum_{\mu\in \Lambda^{*}_{+}}
\Qcal(M_{\mu+\rho_{c}})\, \chi_{_{\mu}}^{_K}\ ,
$$
for every $\kappa$-prequantum line bundle. 
\end{theo}

In contrast to $(3)$, the RHS of $(4)$ is in general an infinite sum. 
Assumption \ref{hypothese.phi.t.carre}
is needed to control the data on the non-compact manifold $M$: it asserts
in particular that for any coadjoint orbit $\Ocal$ of $K$, the square of the moment map
$\Phi_{\Ocal}:M\times\Ocal\to\kgot^*,\,(m,\xi)\mapsto\Phi(m)-\xi,$ has a {\em compact}
set of critical points.

\medskip

In the final section we consider, for $\lambda\in\widehat{G}_d$, 
the case of the coadjoint orbit $G\cdot\lambda$ with the Hamiltonian 
$K$-action. The moment map $\Phi$ is {\em proper} and the critical set 
of $\parallel\Phi\parallel^2$ coincides with $K\cdot\lambda$, hence is 
compact. Thus the generalized Riemann-Roch character 
$RR^{^K}_{\Phi}(G\cdot\lambda,-)$ is well defined, and we want to
investigate the index $RR^{^K}_{\Phi}(M,\tilde{L})$ for the 
$\kappa$-prequantum line bundle $\tilde{L}:=G\times_T\C_{\lambda-\rho}$.

On one hand we are able to compute $RR^{^K}_{\Phi}(G\cdot\lambda,\tilde{L})$ 
explicitly in term of the holomorphic induction map $\HolT$.  
Let $\pgot$ be the orthogonal complement of $\kgot$ in $\ggot$.  
It inherits a complex structure and an action of the torus $T$. 
The element $\wedge^{\bullet}_{\C}\pgot\,\in R(T)$ admits a polarized 
inverse $[\wedge^{\bullet}_{\C}\pgot]^{-1}_{\lambda}\in
R^{-\infty}(T)$ (see \cite{pep4}[Section 5]). In Subsection 
\ref{subsec.preuve.theoreme} we prove that 
$$
(5)\quad\quad RR^{^K}_{\Phi}(G\cdot\lambda,\tilde{L})=
(-1)^{\frac{\dim(G/K)}{2}} \HolT\left(t^{\lambda-\rho_c +\rho_n}
\left[\wedge^{\bullet}_{\C}\pgot\right]^{-1}_{\lambda}\right),
$$
where $\rho_n=\rho-\rho_c$ is half the sum of the non-compact roots. On the other 
hand, we show (Lemma \ref{lem.theta.RR}) that the Blattner formulas can be 
reinterpreted through $\HolT$ as follows:
$$
(6)\quad\quad\Theta_{\lambda}\vert_{K}=\HolT\left(t^{\lambda-\rho_c +\rho_n}
 \left[\wedge^{\bullet}_{\C}\pgot\right]^{-1}_{\lambda}\right)  .
$$
From $(5)$ and $(6)$ we obtain
$$
(7)\quad\quad 
RR^{^K}_{\Phi}(G\cdot\lambda,\tilde{L})=(-1)^{\frac{\dim(G/K)}{2}}\, 
\Theta_{\lambda}\vert_{K} 
={\rm L{\sp 2}-Index}^{_K}(\overline{\partial}_{\tilde{L}}+ 
\overline{\partial}^*_{\tilde{L}}) .
$$

Since in this context $\esp=(-1)^{\frac{\dim(G/K)}{2}}$,  
the Theorem follows from $(4)$ and $(7)$, provided one verifies that  Assumption 
\ref{hypothese.phi.t.carre} holds for $G\cdot\lambda$. This is done in the 
final subsection of this paper.

\bigskip

{\bf Acknowledgments.} I would like to thank  Mich\`ele Vergne for 
suggesting this problem, and helpful discussions.

\bigskip

\bigskip

\bigskip

\begin{center}
    {\bf Notation}
\end{center}

Throughout the paper, $K$ will denote a compact, connected Lie group, 
and $\kgot$ its Lie algebra. In Sections 2, 3, and 4, we consider a
$K$-Hamiltonian action on a manifold $M$. And we use there the following 
notation.

\begin{itemize}
  
   \item[]$T$ : maximal torus of $K$ with Lie algebra $\tgot$
   
   \item[]$W$ : Weyl group of $(K,T)$
  
   \item[]$\Lambda=\ker(\exp:\tgot\to T)$ : integral lattice of $\tgot$
   
   \item[]$\Lambda^*=\hom(\Lambda,2\pi\Z)$ : real weight lattice
   
   \item[]$\tgot^*_+,\rho_c$ : Weyl chamber and corresponding half sum 
   of the positive roots
   
   \item[]$\Lambda^*_+=\Lambda^*\cap\tgot^*_+$ : set of positive weights
   
   \item[]$\chi_{_{\mu}}^{_K}$ : character of the irreducible $K$-representation with highest 
   weight $\mu\in \Lambda^*_+$
  
   \item[]$\tore_{\beta}$ : subtorus of $T$ generated by $\beta\in \tgot$
    
   \item[]$M^{\gamma}$ : submanifold of points fixed by $\gamma\in\kgot$
   
   \item[]$\T M$ : tangent bundle of $M$
   
   \item[]$\T_{K} M$ : set of tangent vectors orthogonal to the 
   $K$-orbits in $M$

   \item[]$\Phi$ :  moment map
   
   \item[]$\tilde{L}$ :  $\kappa$-prequantum line bundle
   
   \item[]$\tilde{\C}_{[\mu]}=K\times_T\C_\mu$ :  $\kappa$-prequantum line 
   bundle over the coadjoint orbit $K\cdot(\mu+\rho_c)$
   
   \item[]$\Cr(\parallel\Phi\parallel^2)$ : critical set of the 
   function $\parallel\Phi\parallel^2$
   
   \item[]$\Delta=\Phi(M)\cap\tgot^*_+$ :  moment polytope
    
   \item[]$\Hcal$ : vector field generated by $\Phi$
   
   \item[]$\mm_{\mu}(E)$ : multiplicity of $RR^{^K}_{\Phi}(M,E)$ 
   relatively to $\mu\in \Lambda^*_+$.
\end{itemize}

In the final section, we consider the particular case of the 
$K$-action on $M:=G\cdot\lambda$. Here $G$ is a connected real semi-simple Lie group 
with finite center admitting $K$ as a maximal compact subgroup, 
and $T$ as a compact Cartan subgroup.

Let us recall the definition of the holomorphic induction map 
$\HolT$. Every $\mu\in\Lambda^*$ defines a 
1-dimensional $T$-representation, denoted $\C_{\mu}$, where 
$t=\exp X$ acts by $t^{\mu}:=e^{i\langle\mu,X\rangle}$. 
We denote by $R(K)$ (resp. $R(T)$) the ring of characters of finite-dimensional 
$K$-representations (resp. $T$-representations).
We denote  $R^{-\infty}(K)$ (resp. $R^{-\infty}(T)$)  the 
set of generalized characters of $K$ (resp. $T$). An element 
$\chi\in R^{-\infty}(K)$ is of the form 
$\chi=\sum_{\mu\in\Lambda^{*}_{+}}\mm_{\mu}\, \chi_{_{\mu}}^{_K}\,$, 
where $\mu\mapsto \mm_{\mu}, \Lambda^{*}_{+}\to\Z$ has at most 
polynomial growth. Likewise, an element $\chi\in R^{-\infty}(T)$
is of the form $\chi=\sum_{\mu\in\Lambda^{*}}\mm_{\mu}\, t^{\mu}$, 
where $\mu\mapsto \mm_{\mu}, \Lambda^{*}\to\Z$ has at most 
polynomial growth. We denote  $w\circ\mu=w(\mu+\rho_{c})-\rho_{c}$ 
the affine action of the Weyl group on $\Lambda^*$. The holomorphic induction map 
$$
\HolT: R^{-\infty}(T)\longrightarrow R^{-\infty}(K)
$$
is characterized by the following properties: 

i) $\HolT(t^{\mu})=\chi_{_{\mu}}^{_K}$ for every $\mu\in\Lambda^{*}_{+}$, 

ii) $\HolT(t^{w\circ\mu})=(-1)^w\HolT(t^{\mu})$ for every $w\in W$ and 
$\mu\in\Lambda^{*}$, 

iii) $\HolT(t^{\mu})=0$ if $W\circ\mu\cap \Lambda^{*}_{+}=\emptyset$.

\bigskip

\section{$\spinc$-quantization of compact Hamiltonian 
$K$-manifolds}\label{sec.quantization.compact}
In this Section we give a modified  version of the `quantization commutes 
with reduction' principle.

Let $M$ be a compact Hamiltonian $K$-manifold with symplectic form $\omega$ 
and moment map $\Phi: M\to \kgot^{*}$ characterized by the relation 
$d\langle\Phi,X\rangle=-\omega(X_M,-)$, where $X_M$ is the vector field on $M$ 
generated by $X\in \kgot$ :  $X_{M}(m):= \frac{d}{dt}\exp(-tX).m |_{t=0}$, 
for $m\in M$. 

Let $J$ be a $K$-invariant almost complex structure on $M$ which is not assumed 
to be compatible with the symplectic form. We denote $RR^{^K}(M,-)$ 
the Riemann-Roch character defined by $J$. Let us recall the definition of this map. 

Let $E\to M$ be a complex $K$-vector bundle. The almost complex structure on $M$ 
gives the decomposition $\wedge \T^{*} M \otimes \C =\oplus_{i,j}\wedge^{i,j}\T^* M$
of the bundle of differential forms. Using Hermitian structure in the tangent 
bundle $\T M$ of $M$, and in the fibers of $E$, we define a 
Dolbeault-Dirac operator $\overline{\partial}_E+ 
\overline{\partial}^*_E
:\Acal^{0,even}(M,E)\to\Acal^{0,odd}(M,E)$,
where $\Acal^{i,j}(M,E):=\Gamma(M,\wedge^{i,j}\T^{*}M\otimes_{\C}E)$ 
is the space of $E$-valued forms of type $(i,j)$. The Riemann-Roch
character $RR^{^K}(M,E)$ is defined as the index of the elliptic
operator $\overline{\partial}_E+ \overline{\partial}^*_E$:
$$
RR^{^K}(M,E)= \indice^K_M(\overline{\partial}_E + \overline{\partial}^*_E)
\quad \in R(K)
$$
viewed as an element of $R(K)$, the character ring of $K$. An 
alternative definition goes as follows. The almost complex structure
defines a canonical invariant $\spinc$ structure\footnote{See 
subsection \ref{Spin-c.structure} for a short review on the notion of 
$\spinc$ structure.}. The $\spinc$ Dirac 
operator of $M$ with coefficient in $E$ has the same principal
symbol as $\sqrt{2}(\overline{\partial}_E+ \overline{\partial}^*_E)$
(see e.g. \cite{Duistermaat96}), and therefore has the same 
equivariant index.

In the Kostant-Souriau framework, $M$ is prequantized if there is a 
$K$-equivariant Hermitian line bundle $L$ with a $K$-invariant 
Hermitian connection $\nabla^L$ of curvature $-i\,\omega$.  The line 
bundle $L$ is called a prequantum line bundle for the Hamiltonian 
$K$-manifold $(M,\omega,\Phi)$. Recall that the data $(\nabla^L,\Phi)$ 
are related by the Kostant formula
\begin{equation}\label{eq.kostant}
   \Lcal^L(X)-\nabla^L_{X_M}=i \langle \Phi,X\rangle, \ 
   X\in\kgot\ .
\end{equation}   
Here $\Lcal^L(X)$ is the infinitesimal action of $X$ on the section 
of $L\to M$.

\medskip

The tangent bundle $\T M$ endowed with $J$ is a complex vector bundle 
over $M$, and we consider its complex dual $\T^{*}_{\C}M:=\hom_{\C}(\T M,\C)$. 
We suppose first that the canonical line bundle 
$\kappa:=\det\T^{*}_{\C}M$ admits a $K$-equivariant square 
root $\kappa^{1/2}$.  If $M$ is prequantized by $L$, a standard 
procedure in the geometric quantization literature is to tensor 
$L$ by the bundle of half-forms $\kappa^{1/2}$ \cite{Woodhouse}. 
We consider the index $RR^{^K}(M,L\otimes\kappa^{1/2})$ instead of 
$RR^{^K}(M,L)$. In many contexts, the tensor product 
$\tilde{L}=L\otimes\kappa^{1/2}$ has a meaning even if 
$L$ nor $\kappa^{1/2}$ exist.

\begin{defi}\label{kappa-quantized}
    An Hamiltonian K-manifold $(M,\omega,\Phi)$, equipped with 
    an almost complex structure, is 
    $\kappa$-prequantized by an equivariant line bundle $\tilde{L}$ 
    if $L_{2\omega}:=\tilde{L}^{2}\otimes\kappa^{-1}$ is a prequantum
    line bundle for $(M,2\omega,2\Phi)$. 
\end{defi}    

The basic examples are the 
regular coadjoint orbits of $K$.  For any $\mu\in\Lambda^{*}_{+}$, 
consider the regular coadjoint orbit $\Ocal^{\mu+\rho_{c}}:=
K\cdot(\mu+\rho_{c})$ with the compatible complex structure. 
The line bundle $\tilde{\C}_{[\mu]}= K\times_{T}\C_{\mu}$ is a 
$\kappa$-prequantum line bundle over $\Ocal^{\mu+\rho_{c}}$, and we have 
\begin{equation}\label{V-mu-rho}
RR^{^K}(\Ocal^{\mu+\rho_{c}},\tilde{\C}_{[\mu]})
= \chi_{_\mu}^{_K}
\end{equation}
for any $\mu\in \Lambda^{*}_{+}$. 

Definition \ref{kappa-quantized} can be rewritten in the $\spinc$ 
setting (see subsection \ref{Spin-c.structure} for a brief review on 
$\spinc$-structures). The almost complex structure induces a $\spinc$ structure $P$
with canonical line bundle $\det_{\C}\T M=\kappa^{-1}$. If $(M,\omega,J)$ is 
$\kappa$-prequantized by $\tilde{L}$ one can twist $P$ by 
$\tilde{L}$, and then define a new $\spinc$ structure 
with canonical line bundle $\kappa^{-1}\otimes\tilde{L}^2=L_{2\omega}$ 
(see Lemma \ref{subsec.spinc.red}).

\begin{defi}\label{spinc-quantized}
    A symplectic manifold $(M,\omega)$ is 
    $\spinc$-prequantized if there exists a $\spinc$ structure 
    with canonical line bundle $L_{2\omega}$ which is a prequantum
    line bundle on $(M,2\omega)$. If a compact Lie group acts on $M$, 
    the $\spinc$-structure is required to be equivariant. Here we take 
    the symplectic orientation on $M$.
\end{defi}

When $(M,\omega,J)$ is $\kappa$-prequantized by $\tilde{L}$, one wants to
compute the $K$-multiplicities of $RR^{^K}(M,\tilde{L})$ 
in geometrical terms, like in Theorem \ref{th.Q-R}.

\begin{defi}\label{def.regular}
An element $\xi\in\kgot^{*}$ is a {\em quasi-regular} value of 
$\Phi$ if all the $K_{\xi}$-orbits in $\Phi^{-1}(\xi)$ have the same dimension. 
A quasi-regular value is  {\em generic} if the submanifold 
$\Phi^{-1}(\xi)$ is of maximal dimension.
\end{defi}

For any quasi-regular value $\xi\in\kgot^{*}$, the reduced space 
$M_{\xi}:=\Phi^{-1}(\xi)/K_{\xi}$ is an orbifold equipped with a 
symplectic structure $\omega_{\xi}$. Let $\tilde{L}$ be a $\kappa$-prequantum line 
bundle over $M$, and let $L_{2\omega}:=\tilde{L}^{2}\otimes\kappa^{-1}$ 
be the corresponding prequantum line bundle for $(M,2\omega)$. For any dominant 
weight $\mu\in \Lambda^{*}_{+}$ such that $\mu+\rho_{c}$ is a quasi-regular value of $\Phi$, 
$$
(L_{2\omega}\vert_{\Phi^{-1}(\mu +\rho_{c})}\otimes \C_{-2(\mu+\rho_{c})})/T
$$
is a prequantum orbifold-line bundle over $(M_{\mu+\rho_{c}},2\omega_{\mu+\rho_{c}})$.

The following Proposition is the main point for computing 
the $K$-multiplicities of $RR^{^K}(M,\tilde{L})$ in terms of  the
reduced spaces $M_{\mu+\rho_{c}}:=\Phi^{-1}(\mu+\rho_{c})/T$, 
$\mu\in\Lambda^{*}_{+}$. It deals with the coherence of the 
definition of an integer valued map 
$\mu\in\Lambda^{*}_{+}\mapsto \Qcal(M_{\mu+\rho_{c}})$. In the next 
proposition we suppose that $(M,\omega,\Phi)$ is a Hamiltonian 
$K$-manifold with {\em proper} moment map. The set 
$\Phi(M)\cap\tgot^*_+$ is denoted by $\Delta$. By the Convexity 
Theorem \cite{Kirwan.84.bis,L-M-T-W,Sjamaar98} it is a convex rational 
polyhedron, referred to as the {\em moment polyhedron}. 

\medskip

\begin{defiprop}\label{prop.Q.mu.ro}
Let $(M,\omega,\Phi)$ be a Hamiltonian $K$-manifold, with proper 
moment map. We denote $\Delta^o$  the relative interior of the moment polyhedron 
$\Delta:=\Phi(M)\cap\tgot^*_+$. Let $\tilde{L}$ be a $\kappa$-prequantum line 
bundle relative to an almost complex structure $J$. Let $\mu\in\Lambda^{*}_{+}$.
\begin{itemize}
  \item If $\mu +\rho_c\notin\Delta^o$, we set 
  $\Qcal(M_{\mu+\rho_{c}})=0$.
  
  \item If $\mu +\rho_{c}$ is a generic quasi-regular value of $\Phi$, 
  then the $\spinc$ prequantization defined by the data $(J,\,\tilde{L})$ 
induces a  $\spinc$ prequantization on the symplectic quotient 
$(M_{\mu +\rho_c}\, , \,\omega_{\mu +\rho_c})$. We denote 
$\Qcal(M_{\mu+\rho_{c}})\in\Z$ the index of the corresponding $\spinc$ 
Dirac operator.

  \item If $\mu +\rho_c\in\Delta^o$, we take $\xi$ generic and quasi-regular 
  sufficiently close to $\mu +\rho_{c}$. The reduced space 
  $M_{\xi}:=\Phi^{-1}(\xi)/T$ 
  inherits a $\spinc$-structure with canonical line bundle 
  $(L_{2\omega}\vert_{\Phi^{-1}(\xi)}\otimes\C_{-2(\mu +\rho_{c})})/T$. 
  The index $\Qcal(M_{\xi})$ of the corresponding $\spinc$ Dirac operator 
  on $M_{\xi}$ does not depend of $\xi$, when $\xi$ is sufficiently close to 
  $\mu +\rho_{c}$ : it is denoted $\Qcal(M_{\mu+\rho_{c}})$. 
\end{itemize}
\end{defiprop}

When $\xi=\mu +\rho_{c}$ is a generic quasi-regular of 
$\Phi$, the line bundle 
$(L_{2\omega}\vert_{\Phi^{-1}(\xi)}\otimes\C_{-2(\mu +\rho_{c})})/T$
is a prequantum line bundle over $(M_{\mu +\rho_c}\, , \,2\omega_{\mu +\rho_c})$:
so the second point of this `definition' is in fact a particular case of 
the third point. But we prefer to keep it since it outlines the main 
point:  
$\spinc$ prequantization is preserved under symplectic reductions.

The existence of $\spinc$-structures on symplectic quotient is proved 
in Subsection \ref{subsec.spinc.red}. The hard part is to show that 
the index $\Qcal(M_{\xi})$ does not depend of $\xi$, for $\xi$ sufficiently 
close to $\mu +\rho_{c}$: it is done in Subsection \ref{def.Q.M.mu}.

Note that Definition \ref{prop.Q.mu.ro} becomes trivial when $\Delta^o$ is not 
included in the interior of the Weyl chamber: 
$\Qcal(M_{\mu+\rho_{c}})=0$ for all $\mu\in\Lambda^{*}_{+}$. However, 
in this paper we work under the assumption that the infinitesimal 
stabilizers for the $K$-action are {\em Abelian}. And that imposes 
$\Delta^o\subset {\rm Interior}\{ {\rm Weyl\ chamber}\}$ (see 
Lemma \ref{lem.Y.non.nul}). 

\medskip

The following 
`{\em quantization commutes with reduction}' Theorem 
holds for the $\kappa$-prequantum line bundles.

\begin{theo}\label{Th.Q-R.pep}
Let $(M,\omega,\Phi)$ be a compact Hamiltonian $K$-manifold 
equipped with an almost complex structure $J$.
Let $\tilde{L}$ be a $\kappa$-prequantum line bundle
over $M$, and let $RR^{^K}(M,-)$ be the Riemann-Roch character  
defined by $J$. If the infinitesimal stabilizers for the 
action of $K$ on $M$ are {\em Abelian}, we have the following 
equality in $R(K)$
\begin{equation}\label{eq.Q-R-tilde}
RR^{^K}(M,\tilde{L})=\esp
\sum_{\mu\in \Lambda^{*}_{+}}Q(M_{\mu+\rho_{c}})\, 
\chi_{_{\mu}}^{_K}\ ,
\end{equation}
where  $\esp=\pm 1$ is the `quotient' of the orientations defined by the almost complex 
structure, and by the symplectic form. 
\end{theo}

Theorem \ref{Th.Q-R.pep} will be proved in a stronger form in Section \ref{section.Q.R}. 

\medskip

Let us now give an example where the 
stabilizers for the action of $K$ on $M$ are {\em not Abelian}, and 
where (\ref{eq.Q-R-tilde}) does not hold. Suppose that the group
$K$ is not Abelian, so we can consider a face $\sigma\neq \{0\}$ 
of the Weyl chamber. Let $\rho_{c,\sigma}$ be half the sum of the 
positive roots which vanish on $\sigma$, and consider the 
coadjoint orbit $M:=K\cdot(\rho_{c}-\rho_{c,\sigma})$ equipped
with its compatible complex structure.  
Since $\rho_{c}-\rho_{c,\sigma}$ belongs to $\sigma$, the
trivial line bundle $M\times\C\to M$ is  
$\kappa$-prequantum, and the 
image of the moment map $\Phi :M\to \kgot^{*}$ does not 
intersect the interior of the Weyl chamber. So 
$M_{\mu+\rho_{c}}=\emptyset$ for every $\mu$, thus 
the RHS of (\ref{eq.Q-R-tilde}) is equal to zero. But the
LHS of (\ref{eq.Q-R-tilde}) is $RR^{^K}(M,\C)$ which is 
equal to $1$, the character of the trivial representation. $\Box$

\medskip

Theorem \ref{Th.Q-R.pep} can be extended in two directions. 
First one can bypass the condition on the stabilizers by the 
following trick. Starting from a $\kappa$-prequantum line bundle
$\tilde{L}\to M$, one can form the product 
$M\times(K\cdot\rho_{c})$ with the coadjoint orbit through 
$\rho_{c}$. The Kunneth formula gives 
$$
RR^{^K}(M\times(K\cdot\rho_{c}),\tilde{L}\boxtimes\C)=
RR^{^K}(M,\tilde{L})\otimes RR^{^K}(K\cdot\rho_{c},\C)=
RR^{^K}(M,\tilde{L})
$$ 
since $RR^{^K}(K\cdot\rho_{c},\C)=1$. Now we can apply Theorem 
\ref{Th.Q-R.pep} to compute the multiplicities of 
$RR^{^K}(M\times(K\cdot\rho_{c}),\tilde{L}\boxtimes\C)$ 
since $\tilde{L}\boxtimes\C$ is a $\kappa$-prequantum line bundle 
over $M\times(K\cdot\rho_{c})$, and the stabilizers for the 
$K$-action on $M\times(K\cdot\rho_{c})$ are Abelian. Finally we see that 
the multiplicity of the irreducible representation with highest 
weight $\mu$ in $RR^{^K}(M,\tilde{L})$ is equal to 
$\esp\Qcal((M\times(K\cdot\rho_{c}))_{\mu +\rho_{c}})$.

On the other hand, we can extend Theorem \ref{Th.Q-R.pep} to the $\spinc$
setting. It will be treated in a forthcoming paper. 

\medskip


\section{Quantization of non-compact Hamiltonian $K$-manifolds}
\label{sec.quant.non.compact}

In this section $(M,\omega,\Phi)$ denotes a Hamiltonian 
$K$-manifold, not necessarily compact, but with {\em proper}  moment map 
$\Phi$. Let $J$ be an almost complex structure on $M$, and 
let $\tilde{L}$ be a $\kappa$-prequantum line bundle 
over $(M,\omega,J)$ (see Def. \ref{kappa-quantized}). 
From Proposition \ref{prop.Q.mu.ro} the infinite sum
\begin{equation}\label{somme-2}
\sum_{\mu\in 
\Lambda^{*}_{+}}\Qcal(M_{\mu+\rho_{c}})\, 
\chi_{_{\mu}}^{_K}   
\end{equation}
is a  well defined element of $\widehat{R}(K):=\hom_{\Z}(R(K),\Z)$.

The aim of this section is to realize this sum as the index of a  
{\em transversally elliptic} symbol naturally associated to the data
$(M,\Phi,J,\tilde{L})$.


\subsection{Transversally elliptic symbols}

Here we give the basic definitions of the theory of
transversally elliptic symbols (or operators) 
defined by Atiyah in \cite{Atiyah74}. For an axiomatic treatment
of the index morphism see Berline-Vergne 
\cite{B-V.inventiones.96.1,B-V.inventiones.96.2} and for a short
introduction see \cite{pep4}.

Let $M$ be a {\em compact} $K$-manifold. Let $p:\T M\to M$ 
be the projection, and let $(-,-)_M$ be a $K$-invariant Riemannian metric.
If $E^{0},E^{1}$ are $K$-equivariant vector bundles over $M$, a 
$K$-equivariant morphism $\sigma \in \Gamma(\T M,\hom(p^{*}E^{0},p^{*}E^{1}))$  
is called a {\em symbol}. The subset of all $(m,v)\in \T M$ where 
$\sigma(m,v): E^{0}_{m}\to E^{1}_{m}$ is not invertible is called 
the {\em characteristic set} of $\sigma$, and is denoted by $\Char(\sigma)$. 

Let $\T_{K}M$ be the following subset of $\T M$ :
$$
   \T_{K}M\ = \left\{(m,v)\in \T M,\ (v,X_{M}(m))_{_{M}}=0 \quad {\rm for\ all}\ 
   X\in\kgot \right\} .
$$

A symbol $\sigma$ is {\em elliptic} if $\sigma$ is 
invertible outside a compact subset of $\T M$ ($\Char(\sigma)$ is
compact), and is {\em transversally elliptic} if the restriction of $\sigma$ 
to $\T_{K}M$ is invertible outside a compact subset  of $\T_{K}M$ 
($\Char(\sigma)\cap \T_{K}M$ is compact). An elliptic symbol $\sigma$ defines 
an element in the equivariant $K$-theory of $\T M$ with compact 
support, which is denoted by $\K_{K}(\T M)$, and the 
index of $\sigma$ is a virtual finite dimensional representation of $K$
\cite{Atiyah-Segal68,Atiyah-Singer-1,Atiyah-Singer-2,Atiyah-Singer-3}.

A {\em transversally elliptic} symbol $\sigma$ defines an element of 
$\K_{K}(\T_{K}M)$, and the index of $\sigma$ is defined as a trace class virtual 
representation of $K$ (see \cite{Atiyah74} for the analytic index and 
\cite{B-V.inventiones.96.1,B-V.inventiones.96.2} for the cohomological one). 
Remark that any elliptic symbol of $\T M$ is transversally elliptic, hence 
we have a restriction map $\K_{K}(\T M)\to \K_{K}(\T_{K}M)$, and 
a commutative diagram
\begin{equation}\label{indice.generalise}
\xymatrix{
\K_{K}(\T M)\ar[r]\ar[d]_{\indice_{M}^K} & 
\K_{K}(\T_{K}M)\ar[d]^{\indice_{M}^K}\\
R(K)\ar[r] &  R^{-\infty}(K)\ .
   }
\end{equation} 

\medskip

Using the {\em excision property}, one can easily show that the index map 
$\indice_{\Ucal}^{K}:\K_{K}(\T_{K}\Ucal)\to R^{-\infty}(K)$
is still defined when $\Ucal$ is a $K$-invariant relatively compact 
open subset of a $K$-manifold (see \cite{pep4}[section 3.1]).


\subsection{Thom symbol deformed by the moment map}\label{ssection.thom}

To a $K$-invariant almost complex 
structure $J$ one associates the Thom symbol $\Thom_{_{K}}(M,J)$, 
and the corresponding Riemann-Roch character 
$RR^{^K}$ when $M$ is compact \cite{pep4}. Let us recall the 
definitions.

Consider a $K$-invariant Riemannian metric $q$ on $M$ such that 
$J$ is orthogonal relatively to $q$, and let 
$h$ be the  Hermitian structure on  $\T M$ defined by : 
$h(v,w)=q(v,w) -i q(Jv,w)$ for 
$v,w\in \T M$. The symbol 
$$
\Thom_{_{K}}(M,J)\in 
\Gamma\left(M,\hom(p^{*}(\wedge_{\C}^{even} \T M),\,p^{*}
(\wedge_{\C}^{odd} \T M))\right)
$$  
at $(m,v)\in \T M$ is equal to the Clifford map
\begin{equation}\label{eq.thom.complex}
 \Clif_{m}(v)\ :\ \wedge_{\C}^{even} \T_m M
\longrightarrow \wedge_{\C}^{odd} \T_m M,
\end{equation}
where $\Clif_{m}(v).w= v\wedge w - c_{h}(v).w$ for $w\in 
\wedge_{\C}^{\bullet} \T_{x}M$. Here $c_{h}(v):\wedge_{\C}^{\bullet} 
\T_{m}M\to\wedge^{\bullet -1} \T_{m}M$ denotes the 
contraction map relative to $h$. Since the map $\Clif_{m}(v)$ is invertible for all
$v\neq 0$, the symbol $\Thom_{_{K}}(M,J)$ is 
elliptic when $M$ is {\em compact}. 

The important point is that for any $K$-vector bundle $E$, 
$\Thom_{_{K}}(M,J)\otimes p^{*}E$ corresponds to the {\em principal symbol} of the  
twisted $\spinc$ Dirac operator $\Dcal_{E}$ \cite{Duistermaat96}. 
So, when $M$ is a {\em compact} manifold, 
the Riemann-Roch character $RR^{^K}(M,-):\K_{K}(M)\to R(K)$ is 
defined by the following relation
\begin{equation}\label{def.RR}
RR^{^K}(M,E)=\indice^{K}_{M}\left(\Thom_{_{K}}(M,J)\otimes p^{*}E\right)\ .
\end{equation}
Since the class of $\Thom_{_{K}}(M,J)$ in $\K_{K}(\T M)$ is 
independent of the choice of the Riemannian structure, the Riemann-Roch character 
$RR^{^K}(M,-)$ also does not depend on this choice.

\medskip

Consider now the case of a {\em non-compact} Hamiltonian $K$-manifold 
$(M,\omega,\Phi)$. We choose a 
$K$-invariant scalar product on $\kgot^{*}$, and  we consider the 
function $\parallel\Phi\parallel^2: M\to \R$. Let $\Hcal$ be the 
Hamiltonian vector field for $\frac{-1}{2}\parallel\Phi\parallel^2$, i.e. 
the contraction of the symplectic form by $\Hcal$ is equal to the $1$-form 
$\frac{-1}{2}d\parallel\Phi\parallel^2$.  In fact the vector field $\Hcal$ 
only depends on $\Phi$. The scalar product on 
$\kgot^{*}$ gives an identification $\kgot^{*}\simeq\kgot$,  
hence $\Phi$ can be consider as a map from $M$ to $\kgot$. 
We have then
\begin{equation}\label{def.H}
\Hcal_{m}=(\Phi(m))_{M}\vert_{m},\quad m\in M\ ,
\end{equation}
where $(\Phi(m))_{M}$ is the vector field on $M$ generated by 
$\Phi(m)\in\kgot$.

\begin{defi}\label{def.thom.loc}
    The Thom symbol deformed by the moment map, which is denoted by 
    $\Thom_{_{K}}^{\Phi}(M,J)$, is defined by the relation
    $$
   \Thom_{_{K}}^{\Phi}(M,J)(m,v):=\Thom_{_{K}}(M,J)(m,v-\Hcal_{m})
   $$
   for any $(m,v)\in\T M$. Likewise, any equivariant map $S:M\to\kgot$ defines a 
   Thom symbol $\Thom_{_{K}}^{S}(M,J)$ deformed by the vector field $S_M : 
   m \to S(m)_M\vert_{m}$ : $\Thom_{_{K}}^{S}(M,J)(m,v):=\Thom_{_{K}}(M,J)(m,v-S_M(m))$.
\end{defi}

Atiyah first proposed to `deform' the symbol of an elliptic operator by 
the vector field induced by an $S^1$-action in order to 
localize its index on the fixed point submanifold, giving then 
another proof of the Lefschetz fixed-point theorem 
\cite{Atiyah74}[Lecture 6]. Afterwards the idea was exploited by Vergne 
to give a proof of the `quantization commutes with reduction' theorem 
in the case of an $S^1$-action \cite{Vergne96}. In \cite{pep4}, we 
extended this procedure for an action of a compact Lie group. 
Here, we use this idea to produce a transversally 
elliptic symbol on a non-compact manifold.

\medskip

The characteristic set of $\Thom_{_{K}}^{\Phi}(M,J)$ corresponds to
$\{(m,v)\in \T M,\ v=\Hcal_m\}$, the graph of the vector field
$\Hcal$. Since $\Hcal$ belongs to the set of tangent vectors to the 
$K$-orbits, we have
\begin{eqnarray*}
\Char\left(\Thom_{_{K}}^{\Phi}(M,J)\right)\cap \T_{K}M&=&\{(m,0)\in \T M,\ 
\Hcal_m=0 \}\\
&\cong& \{m\in M,\ d\parallel\Phi\parallel^2_m=0 \} \ .
\end{eqnarray*}
Therefore the symbol $\Thom_{_{K}}^{\Phi}(M,J)$ is transversally elliptic 
if and only if the set \break 
$\Cr(\parallel\Phi\parallel^2)$ of critical points 
of the function $\parallel\Phi\parallel^2$ is {\em compact}.

\begin{defi}\label{def.RR.phi}
    Let $(M,\omega,\Phi)$ be a Hamiltonian $K$-manifold with 
    $\Cr(\parallel\Phi\parallel^2)$ compact. 
    For any invariant almost complex structure $J$, the symbol 
    $\Thom_{_{K}}^{\Phi}(M,J)$ is {\em transversally elliptic}.
    For any $K$-vector bundle $E\to M$, the tensor product 
    $\Thom_{_{K}}^{\Phi}(M,J)\otimes p^{*}E$ 
    is transversally elliptic and we denote by  
    $$
    RR^{^K}_{\Phi}(M,E)\in R^{-\infty}(K)
    $$ 
    its index\footnote{Here we take a $K$-invariant relatively compact 
    open subset $\Ucal$ of $M$ such that $\Cr(\parallel\Phi\parallel^2)
    \subset \Ucal$. Then the restriction of $\Thom_{_{K}}^{\Phi}(M,J)$ 
    to $\Ucal$ defines a class $\Thom_{_{K}}^{\Phi}(M,J)\vert_{\Ucal}
    \in \K_{K}(\T_{K}\Ucal)$. Since the index map is well defined on 
    $\Ucal$, one sets $RR^{^K}_{\Phi}(M,E):=
    \indice_{\Ucal}^K(\Thom_{_{K}}^{\Phi}(M,J)\vert_{\Ucal}\otimes 
    p^{*}E\vert_{\Ucal})$. A simple application of the excision property
     shows us that the definition does not depend on the choice of 
     $\Ucal$. In order to simplify our notation (when the almost 
     complex structure is understood), we write
     $RR^{^K}_{\Phi}(M,E):=\indice_{M}^K(\Thom_{_{K}}^{\Phi}(M)\otimes p^{*}E)$.}.
     In the same way, an equivariant map $S:M\to\kgot$ defines a 
   transversally elliptic symbol $\Thom_{_{K}}^{S}(M,J)$ if and only 
   if $\{m\in M,\, S_M(m)=0\}$ is compact. If this holds one  
   defines the localized Riemann-Roch character $RR^{^K}_{S}(M,E):=
   \indice_{M}^K(\Thom_{_{K}}^{S}(M)\otimes p^{*}E)$.
\end{defi}

\begin{rem}
  If $M$ is compact the symbols $\Thom_{_{K}}(M,J)$ and 
  $\Thom_{_{K}}^{\Phi}(M,J)$ are homotopic as elliptic symbols, 
  thus the maps $RR^{^K}(M,-)$ and $RR^{^K}_{\Phi}(M,-)$ coincide
  (see section 4 of \cite{pep4}). 
\end{rem}


We end up this subsection with some technical remarks about the 
symbols $\Thom_{_{K}}^{S}(M,J)$ associated to an equivariant map 
$S:M\to\kgot$, and an almost complex structure. 

Let $\Ucal$ be a $K$-invariant open subspace of $M$. The restriction 
$\Thom_{_{K}}^{S}(M,J)\vert_{\Ucal}=\Thom_{_{K}}^{S}(\Ucal,J)$ is transversally 
elliptic if and 
only if $\{m\in M,\, S_M(m)=0\}\cap\Ucal$ is compact. Let 
$j^{^{\Ucal,\Vcal}}:\Ucal\croc\Vcal$ be two $K$-invariant open 
subspaces of $M$, where $j^{^{\Ucal,\Vcal}}$ denotes the inclusion. If
$\{m\in M,\, S_M(m)=0\}\cap\Ucal=\{m\in M,\, S_M(m)=0\}\cap\Vcal$ is compact, 
the excision property tells us that 
$$
j^{^{\Ucal,\Vcal}}_*\left(\Thom_{_{K}}^{S}(\Ucal,J)\right)=
\Thom_{_{K}}^{S}(\Vcal,J)\ ,
$$
where $j^{^{\Ucal,\Vcal}}_*:\K_K(\T_K\Ucal)\to \K_K(\T_K\Vcal)$ is the 
pushforward map (see \cite{pep4}[Section 3]).

\medskip 

\begin{lem}\label{lem.deformation}
  \begin{enumerate}
    \item If $\{m\in M,\, S_M(m)=0\}\cap\Ucal$ is compact, then the class
    defined by $\Thom_{_{K}}^{S}(\Ucal,J)$ in $\K_K(\T_K\Ucal)$ does 
    not depend on the choice of a Riemannian metric.
    
    \item  Let $S^0,S^1:M\to\kgot$ be two equivariant maps. Suppose there exist 
  an open subset $\Ucal\subset M$, and a vector field $\theta$ on 
  $\Ucal$ such that $(S^0_M,\theta)_{_M}$ and $(S^1_M,\theta)_{_M}$ 
  are $>0$ outside a compact subset $\Kcal$ of $\Ucal$. Then, the 
  equivariant symbols $\Thom_{_{K}}^{S^1}(\Ucal,J)$ and 
  $\Thom_{_{K}}^{S^0}(\Ucal,J)$ are transversally elliptic and define 
  the same class in $\K_K(\T_K\Ucal)$.
  
    \item  Let $J^0,J^1$ be two almost complex structures on $\Ucal$, 
    and suppose that \break $\{m\in M,\, S_M(m)=0\}\cap\Ucal$ is compact. 
    The transversally elliptic symbols \break $\Thom_{_{K}}^{S}(\Ucal,J^0)$ 
    and $\Thom_{_{K}}^{S}(\Ucal,J^1)$ define the same class if 
    there exists a homotopy $J^t,\,t\in[0,1]$ of $K$-equivariant 
    almost complex structures between $J^0$ and $J^1$.
  \end{enumerate} 
\end{lem}

{\em Proof.} Two $J$-invariant Riemannian metrics $q_0,q_1$ 
are connected by $q_t:=$ \break $(1-t)q_0+ tq_1$.  Hence the transversally  
elliptic symbols $\Thom_{_{K}}^{S}(\Ucal,J,q_0)$ and \break  
$\Thom_{_{K}}^{S}(\Ucal,J,q_1)$ are tied by the homotopy $t\mapsto 
\Thom_{_{K}}^{S}(\Ucal,J,q_t)$. The point 1. is then proved. 
The proof of 2. is similar to our 
deformation process in \cite{pep2}. Here we consider the maps
$S^t:=t S^1+(1-t)S^0,\, t\in[0,1]$, and the corresponding symbols 
$\Thom_{_{K}}^{S^t}(\Ucal,J)$. The vector field $\theta$, ensures 
that $\Char(\Thom_{_{K}}^{S^t}(\Ucal,J))\cap\T_K\Ucal\subset\Kcal$ is compact.
Hence $t\to \Thom_{_{K}}^{S^t}(\Ucal,J)$ defines a homotopy of 
transversally elliptic symbols. The proof of 3. is identical to the proof 
of Lemma 2.2 in \cite{pep4}. $\Box$

\begin{coro}
When $\{m\in M,\, S_M(m)=0\}$ is compact, the generalized  Riemann-Roch character 
$RR^{^K}_{S}(M,-)$ does not depend on the choice of a Riemannian metric. 
$RR^{^K}_{S}(M,-)$ does not change either if the almost complex structure 
is deformed smoothly and equivariantly in a neighborhood of $\{m\in M,\, S_M(m)=0\}$.
\end{coro}


\medskip

In Subsections \ref{subsec.K.mult} and 
\ref{subsec.localisation},  we set up the 
technical preliminaries that are needed to compute 
the $K$-multiplicity of $RR^{^K}_{\Phi}(M,\tilde{L})$. 

In Section \ref{section.Q.R}, we compute the
$K$-multiplicity of $RR^{^K}_{\Phi}(M,\tilde{L})$, when the moment map 
is {\em is proper}, in terms of the symplectic quotients 
$M_{\mu+\rho_c},\,\mu\in\Lambda_+^*$.


\subsection{Counting the $K$-multiplicities}\label{subsec.K.mult}

Let $E$ be a $K$-vector bundle over 
a Hamiltonian manifold $(M,\omega,\Phi)$ and suppose that 
$\Cr(\parallel\Phi\parallel^2)$ is compact. One wants to compute the 
$K$-multiplicities of $RR^{^K}_{\Phi}(M,E)\in R^{-\infty}(K)$, i.e.
the integers $\mm_{\mu}(E)\in\Z,\, \mu\in \Lambda^{*}_{+}$ such that
\begin{equation}\label{ed.mm.mu.E}
RR^{^K}_{\Phi}(M,E)=\sum_{\mu\in \Lambda^{*}_{+}}
\mm_{\mu}(E)\, \chi_{_{\mu}}^{_K}\ .
\end{equation}
For this purpose one  use the classical `shifting trick'. By 
definition, one has
$\mm_{\mu}(E)=[RR^{^K}_{\Phi}(M,E)\otimes V_{\mu}^*]^K$,
where $V_{\mu}$ is the irreducible $K$-representation with highest 
weight $\mu$, and $V_{\mu}^*$ is its dual. We know from 
(\ref{V-mu-rho}) that the $K$-trace of $V_{\mu}$ is $ \chi_{_\mu}^{_K}=
RR^{^K}(\Ocal^{\tilde{\mu}},\tilde{\C}_{[\mu]})$, where 
\begin{equation}\label{mu-tilde}
\tilde{\mu}=\mu +\rho_c\ .
\end{equation}
Hence the $K$-trace of the dual $V_{\mu}^*$ is equal to 
$RR^{^K}(\overline{\Ocal^{\tilde{\mu}}},\tilde{\C}_{[-\mu]})$, 
where $\overline{\Ocal^{\tilde{\mu}}}$ is the coadjoint orbit 
$\Ocal^{\tilde{\mu}}$ with opposite symplectic structure and opposite 
complex structure. Let $\Thom_{_{K}}(\overline{\Ocal^{\tilde{\mu}}})$ 
be the equivariant Thom symbol on $\overline{\Ocal^{\tilde{\mu}}}$. 
Then the trace of $V_{\mu}^*$ is equal to $\indice^K_{\Ocal^{\tilde{\mu}}}(
\Thom_{_{K}}(\overline{\Ocal^{\tilde{\mu}}})\otimes \tilde{\C}_{[-\mu]})$, and 
finally the multiplicative property of the index 
\cite{Atiyah74}[Theorem 3.5]  gives 
$$
\mm_{\mu}(E)=\left[
\indice^K_{M\times \Ocal^{\tilde{\mu}}}\left(
(\Thom_{_K}^{\Phi}(M)\otimes p^* E)
\odot
(\Thom_{_K}(\overline{\Ocal^{\tilde{\mu}}})\otimes \tilde{\C}_{[-\mu]})
\right)\right]^K\ .
$$
See \cite{Atiyah74,pep4}, for the definition of the exterior product 
$\odot:\K_K(\T_K M)\times\K_K(\T\Ocal^{\tilde{\mu}})\to 
\K_K(\T_K (M\times\Ocal^{\tilde{\mu}}))$.

\medskip

The moment map relative to the Hamiltonian $K$-action on 
$M\times\overline{\Ocal^{\tilde{\mu}}}$ is 
\begin{eqnarray}\label{Phi-mu}
  \Phi_{\tilde{\mu}}:M\times 
  \overline{\Ocal^{\tilde{\mu}}}&\longrightarrow&\kgot^*\nonumber\\
  (m,\xi)&\longmapsto &\Phi(m)-\xi 
\end{eqnarray}

For any $t\in\R$, we consider the map $\Phi_{t\tilde{\mu}}:
M\times \overline{\Ocal^{\tilde{\mu}}}\to \kgot^{*},\ 
  \Phi_{t\tilde{\mu}}(m,\xi):=\Phi(m)-t\,\xi$.
\begin{assum}\label{hypothese.phi.t.carre}
There exists a {\em compact} subset $\Kcal\subset M$, such that, for every 
$t\in[0,1]$, the critical set of the function 
$\parallel\Phi_{t\tilde{\mu}}\parallel^2:M\times\Ocal^{\tilde{\mu}}\to\R$ is 
contained in $\Kcal\times \Ocal^{\tilde{\mu}}$.
\end{assum}

If $M$ satisfies Assumption \ref{hypothese.phi.t.carre} at 
$\tilde{\mu}$, one has a generalized Riemann-Roch character 
$RR_{\Phi_{\tilde{\mu}}}^{^K}(M\times \overline{\Ocal^{\tilde{\mu}}},-)$ 
since $\Cr(\parallel\Phi_{\tilde{\mu}}\parallel^2)$ is compact.

\begin{prop}\label{m-mu-e}
  Let $\mm_{\mu}(E)$ be the multiplicity of $RR^{^K}_{\Phi}(M,E)$ relatively 
  to the highest weight $\mu\in \Lambda^{*}_{+}$. If $M$ satisfies 
  Assumption \ref{hypothese.phi.t.carre} at $\tilde{\mu}$, then 
$$
\mm_{\mu}(E)=\left[ RR_{\Phi_{\tilde{\mu}}}^{^K}(M\times 
\overline{\Ocal^{\tilde{\mu}}},E\boxtimes\tilde{\C}_{[-\mu]})\right]^K\ .
$$
\end{prop}

{\em Proof.}  One has to show that the transversally elliptic symbols 
$\Thom_{_K}^{\Phi}(M)\odot\Thom_{_K}(\overline{\Ocal^{\tilde{\mu}}})$ 
and $\Thom_{_K}^{\Phi_{\tilde{\mu}}}(M\times\overline{\Ocal^{\tilde{\mu}}})$ 
define the same class in $\K_K(\T_K(M\times \Ocal^{\tilde{\mu}}))$ 
when $M$ satisfies Assumption \ref{hypothese.phi.t.carre} at $\tilde{\mu}$. 
Let $\sigma_{1},\sigma_{2}$ be respectively the Thom 
symbols $\Thom_{_K}(M)$ and $\Thom_{_K}(\overline{\Ocal^{\tilde{\mu}}})$. The symbol 
$\sigma_{I}=\Thom_{_K}^{\Phi}(M)\odot\Thom_{_K}(\overline{\Ocal^{\tilde{\mu}}})$ 
is defined by
$$
\sigma_{I}(m,\xi,v,w)=\sigma_{1}(m,v-\Hcal_{m})\odot
\sigma_{2}(\xi,w)\ ,
$$
where $(m,v)\in\T M$, $(\xi,w)\in\T \Ocal^{\tilde{\mu}}$, and
$\Hcal$ is defined in (\ref{def.H}). 
Let $\Hcal^t$ be the vector field on 
$M\times\Ocal^{\tilde{\mu}}$ generated by the map
$\Phi_{t\tilde{\mu}}:M\times\Ocal^{\tilde{\mu}}\to\kgot$.
For $(m,\xi)\in M\times\Ocal^{\tilde{\mu}}$, we have
$\Hcal^t_{(m,\xi)}=(\Hcal^{a,t}_{(m,\xi)},
\Hcal^{b,t}_{(m,\xi)})$ where $\Hcal^{a,t}_{(m,\xi)}\in
\T_{m}M$ and 
$\Hcal^{b,t}_{(m,\xi)}\in\T_{\xi}\Ocal^{\tilde{\mu}}$. 
The symbol $\sigma_{II}=\Thom_{_K}^{\Phi_{\tilde{\mu}}}(M\times 
\overline{\Ocal^{\tilde{\mu}}})$ is defined by 
$$
\sigma_{II}(m,\xi,v,w)=\sigma_{1}(m,v-\Hcal^{a,1}_{m,\xi})
\odot \sigma_{2}(\xi,w-\Hcal^{b,1}_{(m,\xi)} )\ .
$$
We connect $\sigma_{I}$ and $\sigma_{II}$ through 
two homotopies. First we consider the symbol $A$
on $[0,1]\times\T(M\times\Ocal^{\tilde{\mu}})$ defined 
by 
$$
A(t;m,\xi,v,w)=\sigma_{1}(m,v-\Hcal^{a,t}_{m,\xi})
\odot \sigma_{2}(\xi,w-\Hcal^{b,t}_{(m,\xi)})\  ,
$$
for $t\in [0,1]$, and 
$(m,\xi,v,w)\in\T(M\times\Ocal^{\tilde{\mu}})$. We have 
$\Char(A)=\{(t;m,\xi,v,w)\ \vert\ v=\Hcal^{a,t}_{m,\xi}, 
\ {\rm and}\ w=\Hcal^{b,t}_{(m,\xi)}\}
$ 
and 
\begin{eqnarray*}
\Char(A)\bigcap 
[0,1]\times\T_{K}(M\times\Ocal^{\tilde{\mu}})
&=&\{(t;m,\xi,0,0)\ \vert\ (m,\xi)\in
\Cr(\parallel\Phi_{t\tilde{\mu}}\parallel^2)\ \}\\
&\subset& [0,1]\times\Kcal\times\Ocal^{\tilde{\mu}}\ ,
\end{eqnarray*}
where $\Kcal\subset M$ is the compact subset of Assumption
\ref{hypothese.phi.t.carre}. Thus $A$ defines a homotopy of 
transversally elliptic symbols. The restriction of $A$ to $t=1$ 
is equal to $\sigma_{II}$. The restriction of $A$ to $t=0$ 
defines the following transversally elliptic symbol
$$
\sigma_{III}(m,\xi,v,w)=\sigma_{1}(m,v-\Hcal_{m})
\odot \sigma_{2}(\xi,w-\Hcal^{b,0}_{(m,\xi)})
$$
since $\Hcal^{a,0}_{m,\xi}=\Hcal_{m}$ for every 
$(m,\xi)\in M\times\Ocal^{\tilde{\mu}}$. Next, we consider 
the symbol $B$ on $[0,1]\times\T(M\times\Ocal^{\tilde{\mu}})$ defined 
by 
$$
B(t;m,\xi,v,w)=\sigma_{1}(m,v-\Hcal_{m})
\odot \sigma_{2}(\xi,w-t\,\Hcal^{b,0}_{(m,\xi)})\  .
$$
We have 
$\Char(B)=\{(t;m,\xi,v,w)\ \vert\ v=\Hcal_{m}, 
\ {\rm and}\ w=t\,\Hcal^{b,0}_{(m,\xi)}\}
$ 
and 
\begin{eqnarray*}
\lefteqn{\Char(B)\bigcap 
[0,1]\times\T_{K}(M\times\Ocal^{\tilde{\mu}})
\subset}\\
&  &\left\{
(t;m,\xi,v=\Hcal_{m},w=t\,\Hcal^{b,0}_{(m,\xi)})
\quad , \quad \parallel\Hcal_{m}\parallel^2 +
\, t \parallel\Hcal^{b,0}_{(m,\xi)}\parallel^2=0\ 
\right\} \ .
\end{eqnarray*}
In particular $\Char(B)\cap 
[0,1]\times\T_{K}(M\times\Ocal^{\tilde{\mu}})$ is contained 
in \break
$\{(t;m,\xi,0,w=t\,\Hcal^{b,0}_{(m,\xi)})
\, , \, m\in\Cr(\parallel\Phi\parallel^2)\}$
which is compact since  $\Cr(\parallel\Phi\parallel^2)$ 
is compact. So, $B$ defines a homotopy of transversally elliptic 
symbols between $\sigma_{I}=B\vert_{t=0}$ and  
$\sigma_{III}=B\vert_{t=1}$. We have finally proved that 
$\sigma_{I},\sigma_{II},\sigma_{III}$ define the same class
in $\K_{K}(\T_{K}(M\times\Ocal^{\tilde{\mu}}))$. $\Box$

\bigskip

When $E=\tilde{L}$ is a $\kappa$-prequantum line bundle over $M$, the line bundle 
$\tilde{L}\boxtimes\tilde{\C}_{[-\mu]}$ is a $\kappa$-prequantum line bundle 
over $M\times \overline{\Ocal^{\tilde{\mu}}}$. Therefore Proposition \ref{m-mu-e} 
shows  that under Assumption  \ref{hypothese.phi.t.carre} 
the $K$-multiplicities of $RR_{\Phi}^{^K}(M,\tilde{L}) $ have the form
\begin{equation} \label{forme-generale}
\left[RR_{\Phi}^{^K}(\Xcal,\tilde{L}_{\Xcal}) \right]^K\ , 
\end{equation}
where $(\Xcal,\omega_{\Xcal},\Phi)$ is a Hamiltonian K-manifold with 
$\Cr(\parallel\Phi\parallel^2)$ compact, and $\tilde{L}_{\Xcal}$ is a 
$\kappa$-prequantum line bundle over $\Xcal$ relative to a $K$-invariant 
almost complex structure. In order to compute the quantity (\ref{forme-generale}), 
we exploit in the next subsection the localization techniques developed in  
\cite{pep4}.


\subsection{Localization of the map $RR^{^K}_{\Phi}$} 
\label{subsec.localisation}

For a detailed account on the procedure of localization that we use here, 
see Sections 4 and 6 of \cite{pep4}. In this section $(\Xcal,\omega_{\Xcal},\Phi)$ 
is a Hamiltonian $K$-manifold which is equipped with a   $K$-invariant almost 
complex structure, and a $\kappa$-prequantum line bundle $\tilde{L}$. We suppose 
furthermore that \break
$\Cr(\parallel\Phi\parallel^2)$ is compact. We give here a condition 
under which $[RR_{\Phi}^{^K}(\Xcal,\tilde{L})]^K$  depends only on the data 
in the neighborhood of $\Phi^{-1}(0)$. 

\medskip 

For any $\beta\in\kgot$, let $\Xcal^{\beta}$ be 
the symplectic submanifold of points of $\Xcal$ fixed by the torus 
$\tore_{\beta}$ generated by $\beta$.
Following Kirwan \cite{Kirwan84}, the critical set 
$\Cr(\parallel\Phi\parallel^2)$ decomposes as 
\begin{equation}\label{decomp.kirwan}
\Cr(\parallel\Phi\parallel^2)=
\bigcup_{\beta\in\Bcal} C^{^{K}}_{\beta},\quad {\rm with}\quad 
C^{^{K}}_{\beta}=K.(\Xcal^{\beta}\cap \Phi^{-1}(\beta)),
\end{equation}  
where $\Bcal$ is the subset of $\tgot^*_{+}$ defined by 
$\Bcal:=\{\beta\in\tgot^*_{+},\ \Xcal^{\beta}\cap \Phi^{-1}(\beta)\neq\emptyset\}$. 
Since $\Cr(\parallel\Phi\parallel^2)$ is supposed to be compact,  
$\Bcal$ is finite.

For each $\beta\in\Bcal$, let $\Ucal^{^\beta}\croc \Xcal$ be a 
$K$-invariant {\em relatively compact} open neighborhood of 
$C^{^{K}}_{\beta}$ such that $\overline{\Ucal^{^\beta}}\cap 
\Cr(\parallel\Phi\parallel^2)= C^{^{K}}_{\beta}$.  The restriction of 
the transversally elliptic symbol $\Thom_{_{K}}^{\Phi}(\Xcal)$ to the 
subset $\Ucal^{^\beta}$ defines 
$\Thom_{_{K}}^{\Phi}(\Ucal^{^\beta})\in \K_{K}(\T_{K}
\Ucal^{^\beta})$.

\begin{defi}\label{def.RR.beta}
For every $\beta\in \Bcal$, we denote by $RR_{\beta}^{^K}(\Xcal,-)$
the Riemann-Roch character localized near $C^{^{K}}_{\beta}$, which is 
defined  by 
$$
RR_{\beta}^{^K}(\Xcal,E)=\indice^{K}_{\Ucal^{^\beta}}
\left(\Thom_{_{K}}^{\Phi}(\Ucal^{^\beta})\otimes 
p^* E\vert_{\Ucal^{^\beta}}\right)\ ,
$$ 
for every $K$-vector bundle $E\to \Xcal$. 
\end{defi}

The {\em excision property} tells us that 
\begin{equation}\label{localisation-pep}
  RR^{^K}_{\Phi}(\Xcal,E)=\sum_{\beta\in\Bcal}RR_{\beta}^{^K}(\Xcal,E)
\end{equation}  
for every $K$-vector bundle $E\to \Xcal$ (see \cite{pep4}[Section 
4]). In particular, \break 
$[RR^{^K}_{\Phi}(\Xcal,\tilde{L})]^K=
\sum_{\beta}[RR_{\beta}^{^K}(\Xcal,\tilde{L})]^K$, and our 
main point here is to find suitable conditions under which
$[RR_{\beta}^{^K}(\Xcal,\tilde{L})]^K=0$ for $\beta\neq 0$.

\medskip

Let $\beta$ be a non-zero element in $\kgot$.  For every connected component $\Zcal$ of 
$\Xcal^{\beta}$, let $\Ncal_{\Zcal}$ be the normal bundle of ${\Zcal}$ 
in $\Xcal$.  Let $\alpha^{\Zcal}_1,\cdots, \alpha^{\Zcal}_l$ 
be the real infinitesimal 
weights for the action of $\tore_{\beta}$ on the fibers of 
$\Ncal_{\Zcal}\otimes\C$.  The infinitesimal action of $\beta$ on 
$\Ncal_{\Zcal}\otimes\C$ is a linear map with trace equal to 
$\sqrt{-1}\sum_i\langle \alpha^{\Zcal}_i,\beta\rangle$.

\begin{defi}\label{trace-beta} 
    Let us denote by $\tr_{\beta}|\Ncal_{\Zcal}|$ the following 
positive number
$$
\tr_{\beta}|\Ncal_{\Zcal}|:=\sum_{i=1}^l |\langle
\alpha^{\Zcal}_i,\beta\rangle|\ ,
$$
where $\alpha^{\Zcal}_1,\cdots,\alpha^{\Zcal}_l$ are the the real 
infinitesimal weights for the action of $\tore_{\beta}$ on the fibers of
$\Ncal_{\Zcal}\otimes\C$. For any $\tore_{\beta}$-equivariant real vector 
bundle $\Vcal\to\Zcal$ (resp. real $\tore_{\beta}$-equivariant real 
vector space $E$), we define in the same way $\tr_{\beta}|\Vcal|\geq 
0$ (resp. $\tr_{\beta}|E|\geq 0$).
\end{defi} 

\begin{rem}\label{rem.tr.beta}
  If $\Vcal=\Vcal^1\oplus\Vcal^2$, we have $\tr_{\beta}|\Vcal|=
  \tr_{\beta}|\Vcal^1|+\tr_{\beta}|\Vcal^2|$, and if $\Vcal'$ is 
  an equivariant real subbundle of $\Vcal$, we get $\tr_{\beta}|\Vcal|
  \geq \tr_{\beta}|\Vcal'|$. In particular one see
  that $\tr_{\beta}|\Ncal_{\Zcal}|=\tr_{\beta}|\T \Xcal\vert_{\Zcal}|$, 
  and then, if $E_m\subset\T_m \Xcal$ is a $\tore_{\beta}$-invariant real vector 
  subspace for some $m\in \Zcal$, we have $\tr_{\beta}|\Ncal_{\Zcal}|
  \geq \tr_{\beta}|E_m|$.
\end{rem}

The following Proposition and Corollary give us an essential condition under which  
the number  $[RR_{\Phi}^{^K}(\Xcal,\tilde{L})]^K$ only 
depends on data localized in a neighborhhood of $\Phi^{-1}(0)$.

\begin{prop} \label{RR-beta=0}
 Let $\tilde{L}$ be a $\kappa$-prequantum line bundle over 
 $\Xcal$. The multiplicity of the trivial representation in
$RR_{\beta}^{^K}(\Xcal,\tilde{L})$ is equal to zero if
\begin{equation}\label{eq.RR-beta=0}
\parallel\beta\parallel^2 +
\frac{1}{2}\tr_{\beta}|\Ncal_{\Zcal}|-2(\rho_c,\beta)>0
\end{equation}
for every connected component ${\Zcal}$ of $\Xcal^{\beta}$ which 
intersects $\Phi^{-1}(\beta)$. Condition (\ref{eq.RR-beta=0}) 
always holds if $\beta\in\kgot-\{0\}$ is 
$K$-invariant or if $\parallel\beta\parallel>\parallel\rho_c\parallel$. 
\end{prop}

Since every $\beta\in\Bcal$ belongs to the Weyl chamber, we have  
$2(\rho_c,\beta)=\tr_{\beta}|\kgot/\tgot|$, and then (\ref{eq.RR-beta=0}) 
can be rewritten as $\parallel\beta\parallel^2 +\frac{1}{2} 
\tr_{\beta}|\Ncal_{\Zcal}|-\tr_{\beta}|\kgot/\tgot|>0$. 
From  (\ref{localisation-pep}), we get 

\medskip

\begin{coro} \label{coro.RR.loc.0}
  If condition (\ref{eq.RR-beta=0}) holds for all non-zero
$\beta\in\Bcal$, we have 
$$
\left[RR_{\Phi}^{^K}(\Xcal,\tilde{L}) \right]^K=
\left[RR_{0}^{^K}(\Xcal,\tilde{L}) \right]^K
$$
where $RR_{0}^{^K}(\Xcal,-)$ is the Riemann-Roch character localized near
$\Phi^{-1}(0)$ (see Definition \ref{def.RR.beta}). In particular,
$[RR_{\Phi}^{^K}(\Xcal,\tilde{L})]^K=0$ if  (\ref{eq.RR-beta=0})  
holds for all non-zero $\beta\in\Bcal$, and $0\notin{\rm Image}(\Phi)$.
\end{coro}


\subsection{Proof of Proposition 
\ref{RR-beta=0}}\label{subsec.preuve.RR.beta}

When $\beta\in\kgot$ is $K$-invariant, the scalar product 
$(\rho_c,\beta)$ vanishes and then (\ref{eq.RR-beta=0}) trivially holds. 
Let us show that (\ref{eq.RR-beta=0}) holds when 
$\parallel\beta\parallel>\parallel\rho_c\parallel$. Let $\Zcal$ a 
connected component of $\Xcal^{\beta}$ which 
intersects $\Phi^{-1}(\beta)$. Let  $m\in \Phi^{-1}(\beta)\cap \Zcal$, 
 and let $E_m\subset\T_m \Xcal$ be the subspace spanned by $X_\Xcal(m),\,X\in\kgot$. 
We have $E_m\simeq \kgot/\kgot_m$, where 
$\kgot_m:=\{X\in\kgot,\, X_\Xcal(m)=0\}$. Since $\Phi(m)=\beta$, and 
$\Phi$ is equivariant $\kgot_m\subset \kgot_{\beta}:=
\{X\in\kgot,\, [X,\beta]=0\}$, so $\T_m \Xcal$ contains a 
$\tore_{\beta}$-equivariant subspace isomorphic to 
$\kgot/\kgot_{\beta}$. So we have 
$\tr_{\beta}|\Ncal_{\Zcal}|\geq\tr_{\beta}|\kgot/\kgot_{\beta}|=2(\rho_c,\beta)$, 
and then
\begin{eqnarray*}
\parallel\beta\parallel^2 +
\frac{1}{2}\tr_{\beta}|\Ncal_{\Zcal}|-2(\rho_c,\beta)
&\geq&
\parallel\beta\parallel^2 -(\rho_c,\beta)\\
&>&
0
\end{eqnarray*}
since $\parallel\beta\parallel>\parallel\rho_c\parallel$. $\Box$

\medskip

We prove now that condition (\ref{eq.RR-beta=0}) forces 
$[RR_{\beta}^{^K}(\Xcal,\tilde{L})]^K$ to be equal to $0$. Let 
$\mm_{\beta,\mu}(E)\in\Z$ be the $K$-multiplicities of the localized 
Riemann-Roch character $RR_{\beta}^{^K}(\Xcal,E)$  introduced in 
Definition \ref{def.RR.beta} : $RR_{\beta}^{^K}(\Xcal,E)=\sum_{\mu\in \Lambda^{*}_{+}}
\mm_{\beta,\mu}(E)\, \chi_{_{\mu}}^{_K}$. 
We show now that $\mm_{\beta,0}(\tilde{L})=0$, by using the 
formulas of localization that we proved in \cite{pep4} for the maps 
$RR_{\beta}^{^K}(\Xcal,-)$.

\medskip

{\em  First case :} $\beta\in\Bcal$ {\em is a non-zero} 
$K${\em -invariant element of} $\kgot^{*}$.

\medskip

We show here the following relation for the 
multiplicities $\mm_{\beta,\mu}(\tilde{L})$ :
\begin{equation}\label{eq.RR-beta-mu}
    \mm_{\beta,\mu}(\tilde{L})\neq 0
    \ \Longrightarrow \ 
(\mu,\beta) \geq\, \parallel\beta\parallel^2
+\frac{1}{2} \tr_{\beta}|\Ncal_{\Zcal}|\quad {\rm 
for\ some}\ \Zcal\subset\Xcal^{\beta},
\end{equation}
in particular $\mm_{\beta,0}(\tilde{L})=0$.

\medskip 

Since $\tore_{\beta}$  belongs to the center of $K$, $\Xcal^{\beta}$ 
is a symplectic $K$-invariant submanifold of $\Xcal$.  
Let $\Ncal$ be the normal bundle of $\Xcal^{\beta}$ in $\Xcal$. The $K$-invariant
almost complex structure of $\Xcal$ induces a $K$-invariant 
almost complex structure on $\Xcal^{\beta}$,  and a complex structure on the 
fibers of $\Ncal\to \Xcal^{\beta}$. Then we have  a Riemann-Roch character 
$RR_{\beta}^{^K}(\Xcal^{\beta},-)$ localized along $\Xcal^{\beta}\cap 
\Phi^{-1}(\beta)$ with the decomposition 
$RR_{\beta}^{^K}(\Xcal^{\beta},F)=\sum_{\Zcal}
RR_{\beta}^{^K}(\Zcal,F\vert_\Zcal)$, 
where the sum is taken over the connected components  
$\Zcal\subset\Xcal^{\beta}$ which intersect $\Phi^{-1}(\beta)$.
The torus $\tore_{\beta}$ acts 
linearly on the fibers of the complex vector bundle $\Ncal$, thus we 
can associate the polarized complex $K$-vector bundle 
$\Ncal^{+,\beta}$ and $(\Ncal\otimes\C)^{+,\beta}$ (see Definition 
5.5 in \cite{pep4}): for any real $\tore_{\beta}$-weight $\alpha$ 
on $\Ncal^{+,\beta}$, or on $(\Ncal\otimes\C)^{+,\beta}$, we have 
\begin{equation}\label{eq.polarise}
(\alpha,\beta)>0\ .
\end{equation}
We proved the following localization formula in Section 6.2 of 
\cite{pep4} which holds in $\widehat{R}(K)$ for any $K$-vector bundle $E$
over $\Xcal$ :
\begin{equation}\label{RR-beta-1}
RR_{\beta}^{^{K}}(\Xcal,E)=(-1)^{r_{\Ncal}}\sum_{k\in\N}
RR_{\beta}^{^{K}}\left(\Xcal^{\beta},E\vert_{\Xcal^{\beta}}\otimes\det
    \Ncal^{+,\beta}\otimes S^k((\Ncal\otimes\C)^{+,\beta})\right) \ .
\end{equation}
Here $r_{\Ncal}$ is the locally constant function on $\Xcal^{\beta}$ 
equal to the complex rank of $\Ncal^{+,\beta}$, and $S^k(-)$ is the 
$k$-th symmetric product over $\C$.

Let $i:\tgot_{\beta}\croc\tgot$ be the inclusion of the Lie algebra of
$\tore_{\beta}$, and let $i^*: \tgot^*\to \tgot_{\beta}^*$ be the 
canonical dual map. Let us recall the basic relationship between the 
$\tore_{\beta}$-weight on the fibers of a $K$-vector bundle 
$F\to \Xcal^{\beta}$ and the $K$-multiplicities of 
$RR_{\beta}^{^{K}}(\Xcal^{\beta},F)\in \widehat{R}(K)$: if the irreducible 
representation $V_{\mu}$ occurs in $RR_{\beta}^{^{K}}(\Xcal^{\beta},F)$, 
then $i^*(\mu)$ is a $\tore_{\beta}$-weight on the fibers of $F$ 
(see Appendix B in \cite{pep4}). 

If one now uses (\ref{RR-beta-1}), one sees that 
$\mm_{\mu,\beta}(\tilde{L})\neq 0$ 
only if $i^*(\mu)$ is a $\tore_{\beta}$-weight on the fibers of 
some $\tilde{L}\vert_{\Zcal}\otimes\det\Ncal_{\Zcal}^{+,\beta}\otimes 
S^k((\Ncal_{\Zcal}\otimes\C)^{+,\beta})$. Since 
$(i^*(\mu),\beta)=(\mu,\beta)$, (\ref{eq.RR-beta-mu}) will be proved 
if one shows that each $\tore_{\beta}$-weight $\gamma_{_\Zcal}$ on 
$\tilde{L}\vert_{\Zcal}\otimes\det\Ncal_{\Zcal}^{+,\beta}\otimes 
S^k((\Ncal_{\Zcal}\otimes\C)^{+,\beta})$ satisfies
\begin{equation}\label{eq.gamma.Z}
  (\gamma_{_\Zcal},\beta)\geq\parallel\beta\parallel^2
  +\frac{1}{2} \tr_{\beta}|\Ncal_{\Zcal}|\ .
\end{equation}  

Let $\alpha_{_\Zcal}$ be the $\tore_{\beta}$-weight on the fiber of the 
line bundle $\tilde{L}\vert_{\Zcal}\otimes\det\Ncal_{\Zcal}^{+,\beta}$. Since 
any $\tore_{\beta}$-weight on $S^k((\Ncal\otimes\C)^{+,\beta})$ 
satisfies (\ref{eq.polarise}), (\ref{eq.gamma.Z}) holds if
\begin{equation}\label{alpha.beta.Z}
(\alpha_{_\Zcal},\beta)\geq \parallel\beta\parallel^2 
+\frac{1}{2} \tr_{\beta}|\Ncal_{\Zcal}| 
\end{equation}
for every $\Zcal\subset \Xcal^{\beta}$ which intersects $\Phi^{-1}(\beta)$. 
Let $L_{2\omega}$ be the prequantum line bundle on 
$(M,2\omega,2\Phi)$ such that $\tilde{L}^2=L_{2\omega}\otimes \kappa$ 
(where $\kappa$ is by definition equal to 
$\det(\T_{\C}^* \Xcal)\cong \det(\T \Xcal)^{-1}$). We have
$$
(\tilde{L}\vert_{\Zcal}\otimes\det(\Ncal^{+,\beta}))^2=
L_{_{2\omega}}\vert_{\Zcal}\otimes \det(\T \Xcal)^{-1}\vert_{\Zcal}\otimes
\det(\Ncal^{+,\beta})^2\ .
$$
So $2\alpha_{_\Zcal}=\alpha_1 +\alpha_2$ where 
$\alpha_1$, $\alpha_2$  are respectively $\tore_{\beta}$-weights on 
$L_{_{2\omega}}\vert_{\Zcal}$ and \break 
$\det(\T \Xcal)^{-1}\vert_{\Zcal}\otimes\det(\Ncal^{+,\beta})^2$. 
The Kostant formula (\ref{eq.kostant}) on $L_{_{2\omega}}\vert_{\Zcal}$ gives 
$(\alpha_1,X)=2(\beta,X)$ for every $X\in\tgot_{\beta}$, in particular
\begin{equation}\label{alpha-1.beta}
(\alpha_1,\beta)= 2\parallel\beta\parallel^2\ . 
\end{equation}
On $\Zcal$, the complex vector bundle $\T \Xcal$ has the following decomposition, 
$\T \Xcal\vert_{\Zcal}=\T\Zcal \oplus \Ncal^{-,\beta} \oplus 
\Ncal^{+,\beta}$, 
where $\Ncal^{-,\beta}$ is the orthogonal complement of 
$\Ncal^{+,\beta}$ in $\Ncal$: every $\tore_{\beta}$-weight $\delta$ 
on $\Ncal^{-,\beta}$ verifies $(\delta,\beta)<0$. So we get the 
decomposition 
$\det(\T \Xcal)^{-1}\vert_{\Zcal}\otimes\det(\Ncal^{+,\beta})^2=\det(\T\Zcal)\otimes
\det(\Ncal^{-,\beta})^{-1}\otimes\det(\Ncal^{+,\beta})\ ,$
which gives
\begin{equation}\label{alpha-2.beta}
(\alpha_{2},\beta)=\tr_{\beta}|\Ncal_{\Zcal}| 
\end{equation}
since $\tore_{\beta}$ acts trivially on $\T\Zcal$. Finally 
(\ref{alpha.beta.Z}) follows trivially from (\ref{alpha-1.beta}) and 
(\ref{alpha-2.beta}).

\medskip

{\em  Second case :} $\beta\in\Bcal$ {\em such that} $K_{\beta}\neq K$.

\medskip

Consider the induced Hamiltonian action of $K_{\beta}$ on $\Xcal$, with moment 
map $\Phi_{K_{\beta}}:\Xcal\to\kgot_{\beta}^*$. Let $\Bcal'$ be the indexing 
set for the critical point of $\parallel\Phi_{K_{\beta}}\parallel^2$ 
(see (\ref{decomp.kirwan})). Following Definition \ref{def.RR.beta}, for each 
$\beta'\in\Bcal'$ we consider the $K_{\beta}$-Riemann-Roch character 
$RR_{\beta'}^{^{K_{\beta}}}(\Xcal,-)$ localised along        
$C^{^{K_{\beta}}}_{\beta'}=K_{\beta}.(\Xcal^{\beta'}\cap 
\Phi_{K_{\beta}}^{-1}(\beta'))$. Here $\beta$ is a 
$K_{\beta}$-invariant element of $\Bcal'$ with  
$C^{^{K_{\beta}}}_{\beta}=\Xcal^{\beta}\cap \Phi^{-1}(\beta)$.

Let $\HolT:R^{-\infty}(T)\to R^{-\infty}(K)$, 
$\HolTB:R^{-\infty}(T)\to R^{-\infty}(K_{\beta})$, and 
$\HolB:R^{-\infty}(K_{\beta})\to R^{-\infty}(K)$ be the holomorphic 
induction maps (see Appendix B in \cite{pep4}). Recall that 
$\HolT=\HolB\circ\HolTB$. The choice of a Weyl chamber determines a 
complex structure on the real vector space $\kgot/\kgot_{\beta}$. We 
denote by $\overline{\kgot/\kgot_{\beta}}$ the vector space endowed 
with the opposite complex structure.

The induction formula that we proved in 
\cite{pep4}[Section 6] states that
\begin{equation}\label{eq.induction.RR}
RR_{\beta}^{^{K}}(\Xcal,E)=\HolB
\left(RR_{\beta}^{^{K_{\beta}}}(\Xcal,E)\wedge_{\C}^{\bullet}
\overline{\kgot/\kgot_{\beta}}\right)  
\end{equation}
for every equivariant vector bundle $E$. Let us first write 
the decomposition \break 
$RR_{\beta}^{^{K_{\beta}}}(\Xcal,\tilde{L}) 
=\sum_{\mu\in\Lambda_{\beta}^{+}}\mm_{\beta,\mu}(\tilde{L})
\chi_{_{\mu}}^{_{K_{\beta}}}$  
into irreducible characters of $K_{\beta}$. Since $\beta$ is 
$K_{\beta}$-invariant we can use the result of the  First case.
In particular (\ref{eq.RR-beta-mu}) tells us that 
\begin{equation}\label{eq.RR-beta-mu.bis}
    \mm_{\beta,\mu}(\tilde{L})\neq 0
    \ \Longrightarrow \ 
(\mu,\beta) \geq\, \parallel\beta\parallel^2
+\frac{1}{2} \tr_{\beta}|\Ncal_{\Zcal}|
\end{equation}
for some connected component $\Zcal\subset\Xcal^{\beta}$ which 
intersects $\Phi^{-1}(\beta)$.
 
Each irreducible character 
$\chi_{_{\mu}}^{_{K_{\beta}}}$ is equal to $\HolTB(t^{\mu})$, 
so from (\ref{eq.induction.RR}) we get 
$RR_{\beta}^{^{K}}(M,\tilde{L})=
\HolT\Big((\sum_{\mu}m_{\beta,\mu}(\tilde{L})\, 
t^{\mu})\Pi_{\alpha\in\Rgot^+(\kgot/\kgot_{\beta})}(1-t^{-\alpha})\Big)$
where $\Rgot^+(\kgot/\kgot_{\beta})$ is the set of positive $T$-weights on
$\kgot/\kgot_{\beta}$ : so $\langle\alpha,\beta\rangle>0$ for all 
$\alpha\in \Rgot^{+}(\kgot/\kgot_{\beta})$. Finally , we see that 
$RR_{\beta}^{^{K}}(M,\tilde{L})$
is a sum of terms of the form 
$m_{\beta,\mu}(\tilde{L})\,
\HolT(t^{\mu-\alpha_{I}})$ where $\alpha_{I}=\sum_{\alpha\in 
I}\alpha$ and $I$ is a subset of $\Rgot^+(\kgot/\kgot_{\beta})$.
We know that $\HolT(t^{\mu'})$ is either $0$ or the character of an 
irreducible representation (times $\pm 1$) ; in particular 
$\HolT(t^{\mu'})$ is equal to $\pm 1$  only if 
$(\mu',X)\leq 0$ for every $X\in\tgot_{+}$ 
(see Appendix B in \cite{pep4}). So 
$[RR_{\beta}^{^{K}}(M,\tilde{L})]^K\neq 0$ only if there 
exists a weight $\mu$ such that 
$m_{\beta,\mu}(\tilde{L})\neq 0$ and that 
$\HolT(t^{\mu-\alpha_{I}})=\pm 1$. The first condition imposes 
$(\mu,\beta) \geq \parallel\beta\parallel^2
+\frac{1}{2} \tr_{\beta}|\Ncal_{\Zcal}|$ for some connected component 
$\Zcal\subset\Xcal^{\beta}$, and the second one gives 
$(\mu,\beta) \leq (\alpha_{I},\beta)$. 
Combining the two we end up with 
$$
\parallel\beta\parallel^2 +\frac{1}{2} 
\tr_{\beta}|\Ncal_{\Zcal}|\,\leq\, (\alpha_{I},\beta)\,\leq\, 
\sum_{\alpha\in \Rgot^+(\kgot/\kgot_{\beta})}(\alpha,\beta)
=2(\rho_c,\beta)\ ,
$$
for some connected component $\Zcal\subset\Xcal^{\beta}$ which 
intersects $\Phi^{-1}(\beta)$. This completes the proof that 
$[RR_{\beta}^{^{K}}(M,\tilde{L})]^K = 0\,$ 
if $\parallel\beta\parallel^2 +\frac{1}{2} 
\tr_{\beta}|\Ncal_{\Zcal}|> 2(\rho_c,\beta)$ for every 
component $\Zcal\subset\Xcal^{\beta}$ which 
intersects $\Phi^{-1}(\beta)$. $\Box$


\section{Quantization commutes with reduction}\label{section.Q.R}

Let $(M,\omega,\Phi)$ be a Hamiltonian $K$-manifold 
 equipped with an almost complex structure $J$. 
In this section, we assume that the moment map $\Phi$  
is {\em proper} and that the set $\Cr(\parallel\Phi\parallel^2)$ of 
critical points of $\parallel\Phi\parallel^2: M\to\R$ is {\em compact}. 
We denote the corresponding Riemann-Roch character by $RR^{^K}_{\Phi}(M,-)$ 
(see Definition \ref{def.RR.phi}). Let 
$\Delta:=\Phi(M)\cap\tgot^*_+$ be the moment polyhedron.

The main result of this section is the following

\begin{theo}\label{quantization-non-compact-1}
Suppose that $M$ satisfies Assumption 
\ref{hypothese.phi.t.carre} at every  $\tilde{\mu}\in \tgot^*$, and that the
infinitesimal stabilizers for the $K$-action on $M$ are {\em Abelian}. 
If $\tilde{L}$ is a $\kappa$-prequantum line bundle 
over $(M,\omega,\Phi,J)$, we have
\begin{equation}\label{Q-R-noncompact-tilde}
RR^{^K}_{\Phi}(M,\tilde{L})=\esp \sum_{\mu\in\Lambda^{*}_{+}}
\Qcal(M_{\mu+\rho_{c}})\, \chi_{_{\mu}}^{_K}\ ,
\end{equation}
where $\esp=\pm 1$ is the 
`quotient' of the orientation $o(J)$ defined by the almost complex 
structure and the orientation $o(\omega)$ defined by the symplectic 
form. Here the integer $\Qcal(M_{\mu+\rho_{c}})$ is computed by 
Proposition \ref{prop.Q.mu.ro}. In particular, 
$\Qcal(M_{\mu+\rho_{c}})=0$ if $\mu+\rho_{c}$ does not belong to 
the relative interior of $\Delta$.
\end{theo}

\medskip

The same result holds in the traditional `prequantum' case.
Suppose that $M$ satisfies Assumption 
\ref{hypothese.phi.t.carre} at every  $\mu\in \tgot^*$, and that the
almost complex structure $J$ is compatible with $\omega$. If $L$ is prequantum 
line bundle over $(M,\omega,\Phi)$, we have
$RR^{^K}_{\Phi}(M,L)=\sum_{\mu\in\Lambda^{*}_{+}}
RR(M_{\mu},L_{\mu})\, \chi_{_{\mu}}^{_K}$.

\medskip

The next Lemma is the first step in computing the $K$-multiplicities 
$\mm_{\mu}(\tilde{L})$ of $RR^{^K}_{\Phi}(M,\tilde{L})$. Since $(M,\Phi)$ 
satisfies Assumption \ref{hypothese.phi.t.carre} at every $\tilde{\mu}$, 
we know from Proposition \ref{m-mu-e}  that 
$\mm_{\mu}(\tilde{L})=[RR^{^K}_{\Phi_{\tilde{\mu}}}
(M\times\overline{\Ocal^{\tilde{\mu}}},\tilde{L}\boxtimes
\tilde{\C}_{[-\mu]})]^K$ for every $\mu\in\Lambda^{*}_{+}$. 

Let $RR^{^K}_{0}(M\times\overline{\Ocal^{\tilde{\mu}}},-)$ be the 
Riemann-Roch character localized near 
$\Phi^{-1}_{\tilde{\mu}}(0)\simeq \Phi^{-1}(\mu +\rho_c)$ (see 
Definition \ref{def.RR.beta}). This map is the {\em zero map} if $\Phi^{-1}(\mu +\rho_c)=\emptyset$. 

\begin{lem} \label{lem.RR.0.mu}
   Let $\tilde{L}$ be a $\kappa$-prequantum line bundle over $M$. 
   Suppose that the infinitesimal stabilizers for the $K$-action are {\em Abelian} and 
   that Assumption \ref{hypothese.phi.t.carre} is satisfied at 
   $\tilde{\mu}$. We have then
\begin{equation}
\mm_{\mu}(\tilde{L})=\left[ 
RR^{^K}_{0}(M\times\overline{\Ocal^{\tilde{\mu}}},\tilde{L}\boxtimes
\tilde{\C}_{[-\mu]})\right]^K \ .
\end{equation}
In particular $\mm_{\mu}(\tilde{L})=0$ if $\mu+\rho_c$ does not belong 
to the moment polyhedron $\Delta$.
\end{lem}

\medskip

{\em Proof.}  The lemma follows from Corollary \ref{coro.RR.loc.0}, 
applied to the Hamiltonian manifold 
$\Xcal:=M\times\overline{\Ocal^{\tilde{\mu}}}$, with moment map 
$\Phi_{\tilde{\mu}}$ and $\kappa$-prequantum line bundle $\tilde{L}\boxtimes
\tilde{\C}_{[-\mu]}$. Let $\beta\neq 0$ such that 
$\Xcal^{\beta}\cap\Phi_{\tilde{\mu}}^{-1}(\beta)\neq \emptyset$. Let 
$\Ncal$ be the normal bundle of $\Xcal^{\beta}$ in $\Xcal$, and let
$x\in \Xcal^{\beta}\cap\Phi_{\tilde{\mu}}^{-1}(\beta)$. From 
the criterion of Proposition \ref{RR-beta=0}, it is sufficient to show
that
\begin{equation}\label{eq.condition.trace}
\parallel\beta\parallel^2 +
\frac{1}{2}\tr_{\beta}|\Ncal_{x}|-2(\rho_c,\beta)>0\ .
\end{equation}
Write $x=(m,\xi)$ with $m\in M^{\beta}$ and 
$\xi\in (\Ocal^{\tilde{\mu}})^{\beta}$. We know that 
$\tr_{\beta}|\Ncal_{x}|=\tr_{\beta}|\T_{x}\Xcal|=\tr_{\beta}|\T_{m}M|+
\tr_{\beta}|\T_{\xi}\Ocal^{\tilde{\mu}}|$. 

Since the stabilizer $\kgot_{\xi}\simeq\tgot$ is Abelian and 
$\beta\in\kgot_{\xi}$ we have $\kgot_{\xi}\subset\kgot_{\beta}$. 
Then the tangent space 
$\T_{\xi}\Ocal^{\tilde{\mu}}\simeq \kgot/\kgot_{\xi}$ 
contains a copy of $\kgot/\kgot_{\beta}$, so 
$\tr_{\beta}|\T_{\xi}\Ocal^{\tilde{\mu}}|\geq \tr_{\beta}|\kgot/\kgot_{\beta}|=
2(\rho_c,\beta)$. On the other hand, $\T_{m}M$ contains the vector space 
$E_m\simeq\kgot/\kgot_{m}$ spanned by $X_M(m),\, X\in\kgot$. We have 
assume that the stabilizer subalgebra $\kgot_{m}$ is Abelian, and 
since $\beta\in\kgot_{m}$, we get $\kgot_{m}\subset\kgot_{\beta}$. 
Thus $\kgot/\kgot_{\beta}\subset E_m\subset\T_{m}M$ and    
$\tr_{\beta}|\T_{m}M|\geq 2(\rho_c,\beta)$. 
Finally (\ref{eq.condition.trace}) 
is proved since $\frac{1}{2}(\tr_{\beta}|\T_{m}M|+
\tr_{\beta}|\T_{\xi}\Ocal^{\tilde{\mu}}|)\geq 2(\rho_c,\beta)$. $\Box$

\bigskip

The remaining part of this section 
is devoted to the proof of Theorem  \ref{quantization-non-compact-1}. 
Following the preceding Lemma we have to show that 
\begin{equation}\label{equation-finale}
  \left[RR^{^K}_{0}(M\times\overline{\Ocal^{\tilde{\mu}}},\tilde{L}\boxtimes
\tilde{\C}_{[-\mu]})\right]^K=\esp \Qcal(M_{\mu+\rho_{c}})\ , 
\end{equation}
where $\Qcal(M_{\mu+\rho_{c}})$ is defined in Proposition \ref{prop.Q.mu.ro}. 

In  Subsection \ref{Spin-c.structure}, we recall the basic notions about 
$\spinc$-structures. The existence of induced $\spinc$-structures on symplectic quotient 
is proved in Subsection \ref{subsec.spinc.red}. The proof of 
(\ref{equation-finale}) is settled in Subsection \ref{def.Q.M.mu}. We  
give in the same time the proof of the `hard part' of Proposition 
\ref{prop.Q.mu.ro}: the fact that the index $\Qcal(M_{\xi})$ does not depend on $\xi$, for 
$\xi$ sufficiently close to $\mu +\rho_{c}$.


\subsection{\protect $\spinc$ structures and symbols} \label{Spin-c.structure}
We refer to Lawson-Michelson \cite{Lawson-Michel} for background on $\spinc$-structures, and to 
Duistermaat \cite{Duistermaat96} for a discussion of the symplectic case.

The group $\spin_n$ is the connected double cover of the group 
$\so_n$. Let $\eta :\spin_n\to\so_n$ be the covering map, and let 
$\esp$ be the element who generates the kernel. The group $\spinc_n$ 
is the quotient $\spin_n\times_{\Z_2}\u_1$, where $\Z_2$ acts by 
$(\esp,-1)$. There are two canonical group homomorphisms
$$
\eta:\spinc_n\to\so_n\quad ,\quad \Det :\spinc_n\to \u_1\ .
$$
Note that  $\eta^{\rm c}=(\eta,\Det):\spinc_n\to\so_n\times \u_1$ is a 
double covering map.

Let $p:E\to M$ be a oriented Euclidean vector bundle of rank $n$, and 
let $\Pso(E)$ be its bundle of oriented orthonormal frames. A 
$\spinc$-structure on $E$ is a $\spinc_n$ principal bundle 
$\Pspin(E)\to M$, together with a $\spinc$-equivariant map 
$\Pspin(E)\to\Pso(E)$. The line bundle 
\begin{equation}\label{eq.canonical.line}
  \Lfibre:=\Pspin(E)\times_{\Det}\C
\end{equation}
is called the {\em canonical line bundle} associated to $\Pspin(E)$.  
We have then a double covering map\footnote{If $P$, $Q$ are  principal 
bundle over $M$ respectively for the groups $G$ and $H$, we denote simply 
by $P\times Q$ their fibering product over $M$ which is a $G\times H$ 
principal bundle over $M$.}
\begin{equation}\label{eq.spin.covering}
\eta^{\rm c}_E \, : \, \Pspin(E)\longrightarrow\Pso(E)\times\Pu(\Lfibre) \ ,
\end{equation}
where $\Pu(\Lfibre):=\Pspin(E)\times_{\Det}\u_1$ is the associated 
$\u_1$-principal bundle over $M$.

A $\spinc$-structure on a oriented Riemannian manifold is a 
$\spinc$-structure on its tangent bundle. If a group $K$ acts on the 
bundle $E$, preserving the orientation and the Euclidean structure, 
we define a $K$-equivariant $\spinc$-structure by requiring $\Pspin(E)$ 
to be a $K$-equivariant principal bundle, and (\ref{eq.spin.covering}) 
to be $(K\times\spinc_n)$-equivariant.

\medskip

Let $\Delta_{2m}$ be the complex Spin representation of 
$\spinc_{2m}$. Recall that 
$\Delta_{2m}=\Delta_{2m}^+\oplus\Delta_{2m}^-$ 
inherits a canonical Clifford action 
$\clif :\R^{2m}\to\End_{\C}(\Delta_{2m})$ 
which is $\spinc_{2m}$-equivariant, and which interchanges the 
grading : $\clif(v):\Delta_{2m}^{\pm}\to\Delta_{2m}^{\mp}$, for 
every $v\in\R^{2m}$. Let 
\begin{equation}\label{eq.spinor.bundle}
  \Scal(E):=\Pspin(E)\times_{\spinc_{2m}}\Delta_{2m}
 \end{equation} 
be the spinor bundle over $M$, with the grading  
$\Scal(E):=\Scal(E)^+\oplus\Scal(E)^-$. Since 
$E=\Pspin(E)\times_{\spinc_{2m}}\R^{2m}$, the bundle $p^*\Scal(E)$
is isomorphic to $\Pspin(E)\times_{\spinc_{2m}}(\R^{2m}\oplus\Delta_{2m})$.

Let $\overline{E}$ be the bundle $E$ with opposite orientation. A 
$\spinc$ structure on $E$ induces a $\spinc$ on $\overline{E}$, with 
the same canonical line bundle, and such that 
$\Scal(\overline{E})^{\pm}=\Scal(E)^{\mp}$.

\begin{defi}
  Let $\sthom(E): p^*\Scal(E)^+\to p^*\Scal(E)^-$ be the  symbol 
  defined by 
  \begin{eqnarray*}
    \Pspin(E)\times_{\spinc_{2m}}(\R^{2m}\oplus\Delta_{2m}^+)
    &\longrightarrow&
   \Pspin(E)\times_{\spinc_{2m}}(\R^{2m}\oplus\Delta_{2m}^-)\\
  {} [p;v,w]&\longmapsto &[p,v,\clif(v)w]\ .
  \end{eqnarray*} 

When $E$ is the tangent bundle of a manifold $M$, the symbol 
$\sthom(E)$ is denoted by $\sthom(M)$. If a group $K$ acts equivariantly on 
the $\spinc$-structure, we denote by $\sthom_{K}(E)$ the equivariant 
symbol.
\end{defi}

The characteristic set of $\sthom(E)$ is $M\simeq\{{\rm zero\ section\ 
of}\ E\}$, hence it defines a class in $\K(E)$ if $M$ is compact 
(this class is a free generator of the $\K(M)$-module  $\K(E)$ 
\cite{Atiyah-Bott-Shapiro}). When
$E=\T M$, the symbol $\sthom(M)$ corresponds to the {\em principal symbol} 
of the $\spinc$ Dirac operator associated to the $\spinc$-structure 
\cite{Duistermaat96}. If moreover $M$ is compact, the number 
$\Qcal(M)\in\Z$ is defined as the index of $\sthom(M)$. If we change the 
orientation, note that $\Qcal(\overline{M})=-\Qcal(M)$.

\begin{rem}\label{rem.q.change}
  It should be noted that the choice of the metric on the fibers 
  of $E$ is not essential in the construction. Let $g_0,g_1$ be two 
  metric on the fibers of $E$, and suppose that $(E,g_0)$ admits a 
  $\spinc$-structure denoted by $\Pspin(E,g_0)$. The trivial homotopy
  $g_t=(1-t).g_0 +t.g_1$ between the metrics, induces a homotopy 
  between the principal bundles $\Pso(E,g_0)$, $\Pso(E,g_1)$ which
  can be lifted to a homotopy between $\Pspin(E,g_0)$ and a 
  $\spinc$-bundle over $(E,g_1)$. When the base $M$ is compact, the
  corresponding symbols $\sthom(E,g_0)$ and $\sthom(E,g_1)$ define 
  the same class in $\K(E)$.
\end{rem}

These notions extend to the orbifold case. Let $M$ be a manifold with 
a locally free action of a compact Lie group $H$. The quotient 
$\Xcal:=M/H$ is an orbifold, a space with finite quotient 
singularities. A $\spinc$ structure on $\Xcal$ is by definition a
$H$-equivariant $\spinc$ structure on the bundle $\T_{H}M\to M$, where 
$\T_{H}M$ is identified with the pullback of $\T\Xcal$ via the 
quotient map $\pi:M\to\Xcal$. We define in the same way 
$\sthom(\Xcal)\in\K_{orb}(\T\Xcal)$, such that $\pi^{*}\sthom(\Xcal)=
\sthom_{H}(\T_{H}M)$. Here $\K_{orb}$ denotes the 
$K$-theory of {\em proper} vector bundles \cite{Kawasaki81}. The 
pullback by $\pi$ induces an isomorphism $\pi^{*}: 
\K_{orb}(\T\Xcal)\simeq\K_{H}(\T_{H}M)$. The number 
$\Qcal(\Xcal)\in\Z$ is defined as the index of $\sthom(\Xcal)$, or 
equivalently as the multiplicity of the trivial representation in 
$\indice_{M}^{H}(\sthom_{H}(\T_{H}M))$.

\medskip

Consider now the case of a {\em Hermitian} vector bundle $E\to M$, of 
complex rank $m$. The orientation on the fibers of $E$ is given by the 
complex structure $J$. Let $\Pu(E)$ be the bundle of unitary frames on 
$E$. We denote by  ${\rm i}:\u_m\croc\so_{2m}$ the 
canonical inclusion map. We have a 
morphism ${\rm j}:\u_m\to\spinc_{2m}$ which makes the   
diagram 
\begin{equation}
\xymatrix@C=2cm{
 \u_m\ar[r]^{\rm j} \ar[dr]_{{\rm i}\times\det} & 
 \spinc_{2m}\ar[d]^{\eta^{\rm c}}\\
     & \so_{2m}\times\u_1\ .}
\end{equation}
commutative \cite{Lawson-Michel}. Then 
\begin{equation}\label{eq.J.spin}
  \Pspin(E):=\spinc_{2m}\times_{\rm j}\Pu(E)
\end{equation}
defines a 
$\spinc$-structure over $E$, with canonical line bundle equal to 
$\det_{\C}E$. 

\begin{lem} Let $M$ be a manifold equipped with an almost complex 
structure $J$. The symbol $\sthom(M)$  defined by the $\spinc$-structure 
(\ref{eq.J.spin}), and the Thom symbol $\Thom(M,J)$  defined in Section 
\ref{ssection.thom}  coincide.
\end{lem}

{\em Proof.}  The Spinor bundle $\Scal$ is of the form 
$\Pspin(\T M)\times_{\spinc_{2m}}\Delta_{2m}=$ \break 
$\Pu(\T M)\times_{\u_m}\Delta_{2m}$. 
The map $\clif :\R^{2m}\to\End_{\C}(\Delta_{2m})$, when restricted to 
the $\u_m$-equivariant action through ${\rm j}$, is equivalent to the 
Clifford map $\Clif : \R^{2m}\to\End_{\C}(\wedge\C^m)$ (with the 
canonical action of $\u_m$ on $\R^{2m}$ and $\wedge\C^m$). Then 
$\Scal=\wedge_{\C}\T M$ endowed with the Clifford action. $\Box$

\begin{lem} \label{spinc-tordu} Let $P$ be a $\spinc$-structure over $M$, with bundle 
  of spinors $\Scal$, and canonical line bundle $\Lfibre$. For every 
  Hermitian line bundle $L\to M$, there exists a unique 
  $\spinc$-structure $P_L$ with bundle 
  of spinors $\Scal\otimes L$, and canonical line bundle 
  $\Lfibre\otimes L^2$ ($P_L$ is called the $\spinc$-structure $P$ twisted 
  by $L$).
\end{lem}

{\em Proof.}  Take $P_L=P\times_{\u_1}\Pu(L)$.

\bigskip

We finish this subsection with the following definitions. Let 
$(M,o)$ be an oriented manifold. Suppose that 
\begin{itemize}
  \item a connected compact Lie $K$ acts on $M$
  
  \item $(M,o,K)$ carries a $K$-equivariant $\spinc$-structure

  \item one has an equivariant map $\Psi :M\to\kgot$.
\end{itemize}

Suppose first that $M$ is compact. The symbol $\sthom_{_K}(M)$ is then 
elliptic and defines a map 
$$
\Qcal^{^K}(M,-) :\K_{K}(M)\to R(K)
$$
by the relation $\Qcal^{^K}(M,V):=\indice_{M}^K(\sthom_{_K}(M)\otimes 
V)$. Thus $\Qcal^{^K}(M,V)$ is the equivariant index of the $\spinc$ Dirac 
operator on $M$ twisted by $V$.

Let $\Psi_M$ be the equivariant vector field on $M$ defined by 
$\Psi_M(m):=\Psi(m)_M\vert_m$.

\begin{defi}\label{def.S-thom.loc}
  The symbol $\sthom_{_K}(M)$ deformed by the map $\Psi$, 
  which is denoted by $\sthom_{_K}^{\Psi}(M)$, is defined by the relation
  $$
  \sthom_{_{K}}^{\Psi}(M)(m,v)
  :=\sthom_{_{K}}(M)(m,v-\Psi_{M}(m))
  $$
  for any $(m,v)\in\T M$. The symbol $\sthom_{_{K}}^{\Psi}(M)$ is 
  transversally elliptic  if and only 
   if $\{m\in M,\, \Psi_M(m)=0\}$ is compact. When this holds one  
   defines the localized map $\Qcal^{^K}_{\Psi}(M,V):=
   \indice_{M}^K(\sthom_{_{K}}^{\Psi}(M)\otimes V)$.
\end{defi}

\medskip

We finish this section with an adaptation of Lemma 9.4 and Corollary
9.5 of \cite{pep4}[Appendix B]  to the localized map 
$\Qcal^{^K}_{\Psi}(M,-)$. Let $\beta\in \tgot^*_{+}$ be a non-zero element in the 
center of the Lie algebra $\kgot\cong\kgot^*$ of $K$. We suppose here that 
the subtorus $i:\tore_{\beta}\croc K$, which is equal to the closure of 
$\{\exp(t.\beta),\ t\in\R\}$, acts trivially on $M$. 
Let $\mm_{\mu}(V),\, \mu\in \Lambda^{*}_+$ be the 
$K$-multiplicities of $\Qcal^{^K}_{\Psi}(M,V)$.

\begin{lem}\label{lem.poids.tore}
If $\mm_{\mu}(V)\neq 0$, $i^{*}(\mu)$ is a weight for the action of 
$\tore_{\beta}$ on $V\otimes\Lfibre^{\frac{1}{2}}$. If each weight $\alpha$ for 
the action of $\tore_{\beta}$ on $V\otimes\Lfibre^{\frac{1}{2}}$ satisfies 
$(\alpha,\beta)>0$, then \break $[\Qcal^{^K}_{\Psi}(M,V)]^K=0$.
\end{lem}

  
\subsection{\protect $\spinc$ structures on symplectic 
reductions}\label{subsec.spinc.red}

Let $(M,\omega,\Phi)$ be a Hamiltonian $K$-manifold, such that $\Phi$ 
is proper. Let $J$ be a $K$-invariant almost complex structure on $M$. 
And let $\tilde{L}$ be a $\kappa$-prequantum line bundle over 
$(M,\omega,J)$. Since we do not impose a compatibility condition 
between $J$ and $\omega$, the almost complex structure does not descend 
to the symplectic quotients in general. 
Nevertheless we prove in this section that the $\spinc$ 
prequantization defined by the data $(\tilde{L},J)$ induces a 
$\spinc$ prequantization on the symplectic quotients $M_{\mu+\rho_c}$.

\medskip

Let $\Ycal$ be the subset $\Phi^{-1}({\rm interior}(\tgot^{*}_{+}))$. 
When $\Ycal\neq\emptyset$, the Principal-cross-section Theorem  tells us that
$\Ycal$ is a Hamiltonian $T$-submanifold of $M$, with moment map the 
restriction of $\Phi$ to $\Ycal$ \cite{L-M-T-W}.

\begin{lem}\label{lem.Y.non.nul}
If the infinitesimal stabilizers 
for the $K$-action on $M$ are {\em Abelian}, the symplectic slice 
$\Ycal$ is not empty.
\end{lem}  
  
{\em Proof.}  There exists a unique relatively open face $\tau$ of the Weyl chamber 
$\tgot^*_+$ such that $\Phi(M)\cap\tau$ is dense in $\Phi(M)\cap\tgot^*_+$.
The face $\tau$ is called the principal face of $(M,\Phi)$ \cite{L-M-T-W}.
All points in the open face $\tau$ have the same connected 
centralizer $K_{\tau}$. The Principal-cross-section Theorem  tells us that
$\Ycal_{\tau}:=\Phi^{-1}(\tau)$ is a Hamiltonian $K_{\tau}$-manifold, 
where $[K_{\tau},K_{\tau}]$ acts trivially \cite{L-M-T-W}. 
Here we have assumed that the subalgebras
$\kgot_m:=\{X\in \kgot,\, X_M(m)=0\}$, $m\in M$,
are  Abelian. Hence $[\kgot_{\tau},\kgot_{\tau}]\subset \kgot_m$ for 
every $m\in\Ycal_{\tau}$, and this imposes 
$[\kgot_{\tau},\kgot_{\tau}]=0$. Therefore the subgroup $K_{\tau}$ is 
Abelian, and this is the case only if $\tau$ is the interior of 
the Weyl chamber. $\Box$

\medskip

For the remaining of this section, we assume that $\Ycal\neq\emptyset$, 
so that the relative interior $\Delta^o$ of the moment polyhedron 
is a dense subset of $\Phi(\Ycal)$. On $M$, we have the orientation $o(J)$ defined 
by the almost complex structure and the orientation $o(\omega)$ defined by the symplectic 
form.  We denote their `quotient' by $\esp=\pm 1$. On the symplectic 
quotients we will have also two orientations, one induces by 
$\omega$, and the other induces by $J$, with the same `quotient' $\esp$.

\begin{prop}\label{prop.spinc.induit.1}
   The almost complex structure $J$ induces 
   
   i) an orientation $o(\Ycal)$ on $\Ycal$, and
   
   ii) a $T$-equivariant $\spinc$ structure on $(\Ycal,o(\Ycal))$ with canonical line 
   bundle \break $\det_{\C}(\T M\vert_{\Ycal})\otimes\C_{-2\rho_{c}}$.
\end{prop}

{\em Proof.}  On $\Ycal$, we have the decomposition $\T 
M\vert_{\Ycal}=\T\Ycal\oplus[\kgot/\tgot]$, where $[\kgot/\tgot]$ 
denotes the trivial bundle $\Ycal\times\kgot/\tgot$ corresponding of 
the subspace of $\T M\vert_{\Ycal}$ formed by the vector fields 
generated by the infinitesimal action of $\kgot/\tgot$. The choice of 
the Weyl chamber induces a complex structure on $\kgot/\tgot$, and 
hence an orientation $o([\kgot/\tgot])$. This orientation can be 
also defined by a symplectic form of the type 
$\omega_{\kgot/\tgot}(X,Y)=\langle\xi,[X,Y]\rangle$, where $\xi$ belongs 
to the interior $\tgot^*_+$. Let $o(\Ycal)$ be 
the orientation on $\Ycal$ defined by $o(J)\vert_{\Ycal}=o(\Ycal)
o([\kgot/\tgot])$. On $\Ycal$, we have also the orientation 
$o(\omega_{\Ycal})$ defined by the symplectic form $\omega_{\Ycal}$. 
Note that if $o(J)=\esp\, o(\omega)$, we have also $o(\Ycal)=\esp\, 
o(\omega_{\Ycal})$.

Let ${\rm P}:=\spinc_{2n}\times_{\u_{n}}\Pu(\T M)$ be the $\spinc$ structure 
on $M$ induced by $J$ (see (\ref{eq.J.spin})). When restricted to 
$\Ycal$, ${\rm P}\vert_{\Ycal}$ defines a $\spinc$ structure on the 
bundle $\T \Ycal\oplus[\kgot/\tgot]$. Let $q$ be a $T$-invariant 
Riemannian structure on $\T \Ycal\oplus[\kgot/\tgot]$ such that 
$\T \Ycal$ is orthogonal with $[\kgot/\tgot]$, and $q$ equals the 
Killing form on $[\kgot/\tgot]$. Following Remark \ref{rem.q.change}, 
${\rm P}\vert_{\Ycal}$ induces a $\spinc$ structure ${\rm P}'$ on 
$(\T \Ycal\oplus[\kgot/\tgot],q)$, with the same canonical line 
bundle $\Lfibre=\det_{\C}(\T M\vert_{\Ycal})$. Since the 
$\so_{2k}\times\u_{l}$-principal bundle $\Pso(\T\Ycal)\times\u(\kgot/\tgot)$ 
is a reduction\footnote{Here $2n=\dim M$, $2k=\dim\Ycal$ and 
$2l=\dim(\kgot/\tgot)$, so $n=k+l$.}
of the $\so_{2n}$ principal bundle $\Pso(\T \Ycal\oplus[\kgot/\tgot])$, 
we have the commutative diagram

\begin{equation}
\xymatrix@C=2cm{
 {\rm Q}\ar[r]\ar[d] & 
 \Pso(\T\Ycal)\times\u(\kgot/\tgot)\times\Pu(\Lfibre)\ar[d]\\
  {\rm P}'\ar[r]& \Pso(\T 
  \Ycal\oplus[\kgot/\tgot])\times\Pu(\Lfibre)\ ,
  }
\end{equation}
where ${\rm Q}$ is a $(\eta^{\rm c})^{-1}(\so_{2k}\times\u_{l})
\simeq\spinc_{2k}\times\u_{l}$-principal bundle. Finally we see that 
${\rm Q}'={\rm Q}/\u_{l}$ is a $\spinc$ structure on $\T\Ycal$. Since 
$(\u(\kgot/\tgot)\times\Pu(\Lfibre))/\u_{l}\simeq\Pu(\Lfibre\otimes\C_{-2\rho_{c}})$, 
the corresponding canonical line bundle is 
$\Lfibre'=\Lfibre\otimes\C_{-2\rho_{c}}$. $\Box$

\medskip

Let ${\rm Aff}(\Delta)$ be the affine subspace generated by moment 
polyhedron $\Delta$, 
and let $\overrightarrow{\Delta}$ be the subspace of $\tgot^*$ generated by 
$\{m-n\mid m,n\in\Delta\}$. Let $T_{\Delta}$ the subtorus of $T$ with Lie 
algebra $\tgot_{\Delta}$ equal to  the 
orthogonal (for the duality) of $\overrightarrow{\Delta}$. It is not difficult 
to see that $T_{\Delta}$ corresponds to the connected component of the 
principal stabilizer for the $T$-action on $\Ycal$.

Here we consider the symplectic quotient 
$M_{\xi}:=\Phi^{-1}(\xi)/T$ for {\em generic quasi-regular values} 
$\xi\in\Delta^o$ (see Definition \ref{def.regular}). For such $\xi$, 
the fiber $\Phi^{-1}(\xi)$ is a smooth submanifold of $M$, with a locally 
free action of $T/T_{\Delta}$, and with a tubular neighborhood equivariantly 
diffeomorphic to $\Phi^{-1}(\xi)\times\overrightarrow{\Delta}$. Recall 
that $M_{\xi}$ inherits a canonical symplectic form $\omega_{\xi}$.

\begin{prop}\label{prop.spinc.induit.2}
   Let $\mu\in\Lambda^*_+$ such that $\tilde{\mu}=\mu +\rho_c$ 
   belongs to $\Delta$. Let 
   $\tilde{L}$ be a $\kappa$-prequantum line bundle. 
   For every generic quasi-regular value 
   $\xi\in\Delta^o$, the $\spinc$ structures on $\Ycal$, when twisted by 
   $\tilde{L}\vert_{\Ycal}\otimes\C_{-\mu}$,  induces a
    $\spinc$ structure on the reduced space 
    $M_{\xi}:=\Phi^{-1}(\xi)/T$ with canonical line bundle 
    $(L_{2\omega}\vert_{\Phi^{-1}(\xi)}\otimes\C_{-2\tilde{\mu}})/T$.
    Here we have two choices for the orientations : $o(M_{\xi})$ 
    induced by $o(\Ycal)$, and $o(\omega_{\xi})$ defined by the 
    symplectic form $\omega_{\xi}$. They are related by   
    $o(M_{\xi})=\esp\,o(\omega_{\xi})$.
\end{prop}

\begin{rem}The preceding Proposition will be used 
 
  i) when $\xi=\mu +\rho_c$ is a generic quasi-regular value of 
  $\Phi$: the symplectic quotient $(M_{\mu +\rho_c},\omega_{\mu +\rho_c})$ is then 
  $\spinc$ prequantized. Or 
  
  ii) for general $\mu +\rho_c\in\Delta$. One takes then $\xi$ generic 
   quasi-regular  close enough to $\mu +\rho_c$.
\end{rem}

{\em Proof of the Proposition.}  Let $\xi\in\Delta^o$ be a generic quasi-regular value 
of $\Phi$, and $\Zcal:=\Phi^{-1}(\xi)$. This is a submanifold of $\Ycal$ 
with a trivial action of $T_{\Delta}$ and a locally free action of 
$T/T_{\Delta}$. We denote the quotient map by $\pi:\Zcal\to M_{\xi}$. We 
identify $\pi^{*}(\T M_{\xi})$ with the orthogonal complement 
(with respect to a Riemannian metric) of the trivial bundle 
$[\tgot/\tgot_{\Delta}]$ formed by the vector fields 
generated by the infinitesimal action of $\tgot/\tgot_{\Delta}$. 
On the other hand the tangent bundle 
$\T\Ycal$, when restricted to $\Zcal$, decomposes as 
$\T\Ycal\vert_{\Zcal}=\T\Zcal\oplus[\overrightarrow{\Delta}]$, so 
we have
\begin{eqnarray}\label{eq.decomposition.Y.Z}
\T\Ycal\vert_{\Zcal}&=&\pi^{*}(\T 
M_{\xi})\oplus[\tgot/\tgot_{\Delta}]\oplus[\overrightarrow{\Delta}]\nonumber\\
&=&\pi^{*}(\T M_{\xi})\oplus[\tgot/\tgot_{\Delta}\otimes\C]\ ,
\end{eqnarray}
with the convention 
$\tgot/\tgot_{\Delta}=\tgot/\tgot_{\Delta}\otimes i\R$ and
$\overrightarrow{\Delta}=\tgot/\tgot_{\Delta}\otimes\R$.
Since $\tgot/\tgot_{\Delta}\otimes\C$ is canonically oriented by 
the complex multiplication by $i$, the orientation $o(\Ycal)$  
determines an orientation $o(M_{\xi})$ on $\T M_{\xi}$ 
through (\ref{eq.decomposition.Y.Z}). 

Now we proceed like the proof of Proposition \ref{prop.spinc.induit.2}. 
Let ${\rm Q}'$ be the $\spinc$ structure on $\Ycal$ 
introduced in Proposition \ref{prop.spinc.induit.2}, and let 
${\rm Q}^{\mu}$ be ${\rm Q}'$ twisted by the line bundle 
$\tilde{L}\vert_{\Ycal}\otimes\C_{-\mu}$ : its canonical line 
bundle is $\det_{\C}(\T M)\vert_{\Ycal}\otimes\C_{-2\rho_c}
\otimes(\tilde{L}\vert_{\Ycal}\otimes\C_{-\mu})^2
=L_{2\omega}\vert_{\Ycal} \otimes\C_{-2\tilde{\mu}}$. 
The $\so_{2k'}\times\u_{l'}$-principal bundle 
$\Pso(\pi^{*}(\T M_{\xi}))\times\u(\tgot/\tgot_{\Delta}\otimes\C)$ 
is a reduction\footnote{Here $2k=\dim \Ycal$, $2k'=\dim M_{\xi}$ and 
$l'=\dim(\tgot/\tgot_{\Delta})$, so $k=k'+l'$.}
of the $\so_{2k}$ principal bundle 
$\Pso(\pi^{*}(\T M_{\xi})\oplus[\tgot/\tgot_{\Delta}\otimes\C])$;  
we have the commutative diagram

\begin{equation}\label{diagram.Q.seconde}
\xymatrix@C=2cm{
 {\rm Q}''\ar[r]\ar[d] & 
 \Pso(\pi^{*}(\T M_{\xi}))\times\u(\tgot/\tgot_{\Delta}\otimes\C)
 \times\Pu(\Lfibre\vert_{\Zcal})\ar[d]\\
  {\rm Q}^{\mu}\vert_{\Zcal}\ar[r]   & 
  \Pso(\pi^{*}(\T M_{\xi})\oplus[\tgot/\tgot_{\Delta}\otimes\C])
  \times\Pu(\Lfibre\vert_{\Zcal})\ ,
  }
\end{equation}
where $\Lfibre=L_{2\omega}\vert_{\Ycal} \otimes\C_{-2\tilde{\mu}}$. Here
${\rm Q}''$ is a $(\eta^{\rm c})^{-1}(\so_{2k'}\times\u_{l'})
\simeq\spinc_{2k'}\times\u_{l'}$-principal bundle. The Kostant 
formula (\ref{eq.kostant}) tells us that the action of 
$T_{\Delta}$ is trivial on $\Lfibre\vert_{\Zcal}$, since $\xi-\tilde{\mu}\in 
\overrightarrow{\Delta}$. Thus the action of $T_{\Delta}$ is trivial on 
${\rm Q}''$. Finally we see that 
${\rm Q}_{\xi}={\rm Q}''/(\u_{l'}\times T)$ is a $\spinc$ structure 
on $M_{\xi}$ with canonical line bundle 
$\Lfibre_{\xi}=(L_{2\omega}\vert_{\Zcal} \otimes\C_{-2\tilde{\mu}})/T$. $\Box$


\subsection{Definition of \protect 
$\Qcal(M_{\mu+\rho_{c}})$}\label{def.Q.M.mu}

First we give three different ways to define the quantity  
$\Qcal(M_{\mu+\rho_{c}})\in\Z$ for any $\mu\in\Lambda^{*}_{+}$.
The compatibility of these different definitions proves 
Theorem \ref{quantization-non-compact-1} and the `hard part' of 
Proposition \ref{prop.Q.mu.ro} simultaneously. First of all 
$\Qcal(M_{\mu+\rho_{c}})=0$ if $\mu+\rho_{c}\notin \Delta$. 

{\em First definition}.

If $\mu+\rho_{c}\in\Delta^o$ is a {\em generic quasi regular} value of $\Phi$, 
$M_{\mu+\rho_{c}}:=\Phi^{-1}(\mu+\rho_{c})/T$ is a symplectic 
orbifold.  We know from Proposition \ref{prop.spinc.induit.2} that 
$M_{\mu+\rho_{c}}$ inherits $\spinc$-structures, with the same 
canonical line bundle $(L_{2\omega}\vert_{\Phi^{-1}(\tilde{\mu})} 
\otimes\C_{-2\tilde{\mu}})/T$, for the two choices of orientation 
$o(M_{\mu+\rho_{c}})$ and $o(\omega_{\mu+\rho_{c}})$.  We denote 
the index of the $\spinc$ Dirac operator associated to the $\spinc$ structure on 
$(M_{\mu+\rho_{c}},o(\omega_{\mu+\rho_{c}}))$ by $\Qcal(M_{\mu+\rho_{c}})\in\Z$ 
and the index of the $\spinc$ Dirac operator associated to the $\spinc$ structure 
on $(M_{\mu+\rho_{c}},o(M_{\mu+\rho_{c}}))$ by 
$\Qcal(M_{\mu+\rho_{c}},o(M_{\mu+\rho_{c}}))$.  Since 
$o(M_{\mu+\rho_{c}})=\esp\,o(\omega_{\mu+\rho_{c}})$, we have 
$\Qcal(M_{\mu+\rho_{c}})=\esp\,\Qcal(M_{\mu+\rho_{c}},o(M_{\mu+\rho_{c}}))$.

{\em Second definition}.

We can also define $Q(M_{\mu+\rho_{c}})$ by shift `desingularization' 
as follows.  If $\mu+\rho_{c}\in\Delta$, one considers generic quasi 
regular values $\xi\in\Delta^o$, close enough to $\mu+\rho_{c}$.  
Following Proposition \ref{prop.spinc.induit.2}, 
$M_{\xi}=\Phi^{-1}(\xi)/T$ inherits a $\spinc$ structure, with 
canonical line bundle $(L_{2\omega}\vert_{\Phi^{-1}(\xi)} 
\otimes\C_{-2\tilde{\mu}})/T$.  Then we set 
$\Qcal(M_{\mu+\rho_{c}}):=\Qcal(M_{\xi})$, where the RHS is the index 
of the $\spinc$ Dirac operator associated to the $\spinc$ structure on 
$(M_{\xi},o(\omega_{\xi}))$.  If we take the orientation 
$o(M_{\xi})$ induced by $o(\Ycal)$ we have another index  
$\Qcal(M_{\xi},o(M_{\xi}))= \esp\ \Qcal(M_{\xi})$. Here one has to
show that these quantities does not depend on the choice of $\xi$ when 
$\xi$ is close enough to $\mu+\rho_{c}$.  We will 
see  that $\Qcal(M_{\mu+\rho_{c}})=0$ when 
$\mu+\rho_{c}\notin\Delta^o$.

{\em Third definition}.

We can use the characterization of the multiplicity 
$\mm_{\mu}(\tilde{L})$ given in Lemma \ref{lem.RR.0.mu}. The number 
$\Qcal(M_{\mu+\rho_{c}})$ is the multiplicity
of the trivial representation in 
$\esp\,
RR_{0}^{^K}(M\times\overline{\Ocal^{\tilde{\mu}}},\tilde{L}\otimes\tilde{\C}_{[-\mu]})$. 

\medskip

We have to show the compatibility of these definitions, that is
\begin{equation}\label{def.compatible.1}
  {\rm if}\ \mu+\rho_c\in\Delta^o:\quad
\left[RR_{0}^{^K}(M\times\overline{\Ocal^{\tilde{\mu}}},\tilde{L}
\otimes\tilde{\C}_{[-\mu]})\right]^K = \Qcal(M_{\xi},o(M_{\xi}))
\end{equation}
for any generic quasi regular value $\xi\in\Delta^o$ close enough 
to $\mu+\rho_{c}$. And 
\begin{equation}\label{def.compatible.2}
{\rm if}\ \mu+\rho_c\notin\Delta^o:\quad
\left[RR_{0}^{^K}(M\times\overline{\Ocal^{\tilde{\mu}}},
\tilde{L}\otimes\tilde{\C}_{[-\mu]})\right]^K=0\ .
\end{equation}

We have proved already (Lemma \ref{lem.RR.0.mu}) that 
$[RR_{0}^{^K}(M\times\overline{\Ocal^{\tilde{\mu}}},
\tilde{L}\otimes\tilde{\C}_{[-\mu]})]^K=0$ if $\mu+\rho_{c}\notin \Delta$.

\bigskip

{\bf We work now with a fixed element $\mu\in\Lambda^*_{+}$ such 
that $\tilde{\mu}=\mu+\rho_c$ belongs to $\Delta$. During the
remaining part of this section, $\Ycal$ will denote a small
$T$-invariant open neighborhood of $\Phi^{-1}(\mu+\rho_{c})$ in the 
symplectic slice $\Phi^{-1}({\rm interior}(\tgot^*_+))$.}

\medskip

We will check in Subsection \ref{preuve.prop.2}, that the functions  
$\parallel\Phi-\tilde{\mu}\parallel^2$ and 
$\parallel\Phi-\xi\parallel^2$ have compact critical set on $\Ycal$ 
when $\xi\in{\rm Aff}(\Delta)$ is close enough to $\tilde{\mu}$.
Since the manifold $(\Ycal,o(\Ycal))$ carries a $T$-invariant 
$\spinc$-structure, we consider the localized  maps 
$\Qcal^{^T}_{\Phi-\tilde{\mu}}(\Ycal,-)$ 
and $\Qcal^{^T}_{\Phi-\xi}(\Ycal,-)$ (see Definition \ref{def.S-thom.loc}).
The proof of (\ref{def.compatible.1}) and (\ref{def.compatible.2}) is
divided into two steps. We first relate the maps 
$RR_{0}^{^K}(M\times\overline{\Ocal^{\tilde{\mu}}},-)$ and 
$\Qcal^{^T}_{\Phi-\tilde{\mu}}(\Ycal,-)$ through the induction map
\begin{equation}\label{eq.induction.T.K}
    \indT : \fgene(T)\longrightarrow \fgene(K)^{K}\ .   
\end{equation}
Here $\fgene(T)$, $\fgene(K)$ denote respectively the set of generalized 
functions on $T$ and $K$, and the $K$ invariants are taken with the conjugation action.
The map $\indT$ is defined as follows : for $\phi\in\fgene(T)$, 
we have $\int_{K}\indT(\phi)(k)f(k)dk
=\frac{\Vol(K,dk)}{\Vol(T,dt)}\int_{T}\phi(t)f|_{T}(t)dt$, 
for every $f\in\f(K)^{K}$.

\begin{prop}\label{prop.1}
  Let $E$ and $F$ be respectively $K$-equivariant complex vector 
  bundles over $M$ and $\Ocal^{\tilde{\mu}}$. We have the following 
  equality 
  $$
  RR_{0}^{^K}(M\times\overline{\Ocal^{\tilde{\mu}}},E \boxtimes F)
  =\indT\left(\Qcal^{^T}_{\Phi-\tilde{\mu}}(\Ycal,E\vert_{\Ycal}\otimes 
  F\vert_{\bar{e}})\right)
  $$  
  in $R^{-\infty}(K)$. It gives in particular that 
$$
\left[RR_{0}^{^K}
(M\times\overline{\Ocal^{\tilde{\mu}}},E\boxtimes F)\right]^K=
\left[\Qcal^{^T}_{\Phi-\tilde{\mu}}(\Ycal,E\vert_{\Ycal}\otimes 
F\vert_{\bar{e}})\right]^T\ .
$$
\end{prop}

After we compute the map $\Qcal^{^T}_{\Phi-\tilde{\mu}}(\Ycal,-)$ by 
making the shift $\tilde{\mu}\to\xi$.

\begin{prop}\label{prop.2}
Suppose $\xi\in{\rm Aff}(\Delta)$ is close enough to $\tilde{\mu}$. Then 

i) the maps $\Qcal^{^T}_{\Phi-\tilde{\mu}}(\Ycal,-)$ and 
$\Qcal^{^T}_{\Phi-\xi}(\Ycal,-)$ are equal,
        
ii) if furthermore $\xi\in\Delta^o$ is a generic quasi-regular value of $\Phi$ we 
get 
$$
\left[\Qcal^{^T}_{\Phi-\xi}
(\Ycal,\tilde{L}\vert_{\Ycal}\otimes \C_{-\mu})\right]^T=
\Qcal(M_{\xi},o(M_{\xi})) \ ,
$$

iii) and if $\xi\notin\Delta$, $[\Qcal^{^T}_{\Phi-\xi}
(\Ycal,\tilde{L}\vert_{\Ycal}\otimes \C_{-\mu})]^T=0$.
\end{prop}

Finally, if $\xi\in{\rm Aff}(\Delta)$ is close enough to $\tilde{\mu}$, 
Propositions \ref{prop.1} and \ref{prop.2} give
\begin{eqnarray}\label{eq.bord.delta}
\left[RR_{0}^{^K}
(M\times\overline{\Ocal^{\tilde{\mu}}},\tilde{L}\otimes\tilde{\C}_{[-\mu]})\right]^K
&=& \left[\Qcal^{^T}_{\Phi-\tilde{\mu}}(\Ycal,\tilde{L}\vert_{\Ycal}
\otimes\C_{-\mu})\right]^T\\
&=& \left[\Qcal^{^T}_{\Phi-\xi}(\Ycal,\tilde{L}\vert_{\Ycal}
\otimes\C_{-\mu})\right]^T \nonumber .
\end{eqnarray}

If $\mu +\rho_{c}\in\Delta^o$, we choose $\xi\in\Delta^o$ close to $\mu+\rho_{c}$: 
in this case (\ref{def.compatible.1}) follows from 
(\ref{eq.bord.delta}) and the point $ii)$ of Proposition \ref{prop.2}.
If $\mu +\rho_{c}\notin\Delta^o$, we  choose $\xi$ close to $\mu+\rho_{c}$ and 
not in $\Delta$: in this case (\ref{def.compatible.2}) follows from 
(\ref{eq.bord.delta}) and the point $iii)$ of Proposition 
\ref{prop.2}.

\medskip

Propositions \ref{prop.1} and \ref{prop.2} are proved in the next subsections.

\subsection{Proof of Proposition \ref{prop.1}}\label{preuve.prop.1}

The induction formula of Proposition \ref{prop.1} is essentially
identical to the one we proved in \cite{pep4}. The main difference is 
that the almost complex structure is not assumed to be compatible 
with the symplectic structure.

We identify the coadjoint orbit $\Ocal^{\tilde{\mu}}$ with 
$K/T$. Let $\Hcal^{\tilde{\mu}}$ be the 
Hamiltonian vector field of the function 
$\frac{-1}{2}\parallel\Phi_{\tilde{\mu}}\parallel^{2}:M\times \overline{K/T}\to\R$. 
Here $\Ycal$ denotes a small neighborhood 
of $\Phi^{-1}(\tilde{\mu})$ in the symplectic slice 
$\Phi^{-1}({\rm interior}(\tgot^*_+))$ such that 
the open subset $\Ucal:=(K\times_{T}\Ycal)\times \overline{K/T}$ is a 
neighborhood of $\Phi_{\tilde{\mu}}^{-1}(0)=
K\cdot(\Phi^{-1}(\tilde{\mu})\times\{\bar{e}\})$ 
which verifies $\overline{\Ucal}\cap\{ \Hcal^{\tilde{\mu}}=0\}
=\Phi_{\tilde{\mu}}^{-1}(0)$.

From Definition \ref{def.RR.beta}, the localized Riemann-Roch 
character $RR^{^K}_{0}(M\times \overline{K/T},-)$ is computed by means of 
the Thom class $\Thom_{K}^{\Phi_{\tilde{\mu}}}(\Ucal)\in 
\K_{K}(\T_{K}\Ucal)$. On the other hand, the localized map
$\Qcal^{^T}_{\Phi-\tilde{\mu}}(\Ycal,-)$ is computed by means of 
the class $\sthom_{T}^{\Phi-\tilde{\mu}}(\Ycal)\in 
\K_{T}(\T_{T}\Ycal)$ (see Definition \ref{def.S-thom.loc}). 
Proposition \ref{prop.1} will follow from 
a simple relation between these two transversally elliptic symbols.

First, we consider  the isomorphism  $\phi:\Ucal\to\Ucal'$, 
$\phi([k;y],[h])=[k;[k^{-1}h],y]$,
where $\Ucal':=K\times_{T}(\overline{K/T}\times\Ycal)$. Let
$\phi^{*}:\K_{K}(\T_{K}\Ucal')\to \K_{K}(\T_{K}\Ucal)$ be
the induced isomorphism. Then one considers the inclusion 
$i:T\croc K$ which induces an isomorphism 
$i_{*}:\K_{T}(\T_{T}(\overline{K/T}\times\Ycal))\to
\K_{K}(\T_{K}\Ucal')$ (see \cite{Atiyah74,pep4}). Let 
$j:\Ycal\croc \overline{K/T}\times\Ycal$ be the 
$T$-invariant inclusion map defined by $j(y):=(\bar{e},y)$. We 
have then a pushforward map $j_{!}:
\K_{T}(\T_{T}\Ycal)\to 
\K_{T}(\T_{T}(\overline{K/T}\times\Ycal))$. 
Finally we get  a map 
$$
\Theta:=\phi^{*}\circ i_{*}\circ j_{!}\ :\ 
\K_{T}(\T_{T}\Ycal)\longrightarrow \K_{K}(\T_{K}\Ucal)\ ,
$$
such that $\indice_{\Ucal}^{K}(\Theta(\sigma))=
\indT(\indice_{\Ycal}^{T}(\sigma))$ for every 
$\sigma\in \K_{T}(\T_{T}\Ycal)$ (see Section 3 in \cite{pep4}).

Proposition \ref{prop.1} is an immediate consequence of the following 
\begin{lem}
    We have the equality 
    $$
 \Theta\left(\sthom_{T}^{\Phi-\tilde{\mu}}(\Ycal)\right)=
 \Thom_{K}^{\Phi_{\tilde{\mu}}}(\Ucal)\ .
    $$
 \end{lem} 
 
{\em Proof.}  Let $\Scal$ be the bundle of spinors on 
$K\times_{T}\Ycal$: $\Scal={\rm P}\times_{\spinc_{2k}}\Delta_{2k}$, where 
${\rm P}\to \Pso(\T(K\times_{T}\Ycal))\times\Pu(\Lfibre)$ is the $\spinc$ 
structure induced by the complex structure. From Proposition 
\ref{prop.spinc.induit.1}, we have the reductions 

\begin{equation}\label{eq.2.reduction}
\xymatrix@C=2cm{
{\rm Q}\ar[r]\ar[d] & 
\Pso(\T\Ycal)\times\u(\kgot/\tgot)\times\Pu(\Lfibre\vert_{\Ycal})\ar[d]\\
{\rm P}\vert_{\Ycal}\ar[r] \ar[d]  & 
\Pso(\T\Ycal\oplus[\kgot/\tgot])\times\Pu(\Lfibre\vert_{\Ycal})\ar[d]\\
{\rm P}\ar[r] & \Pso(\T(K\times_{T}\Ycal))\times\Pu(\Lfibre)\ .
  }
\end{equation}
Here $Q/\u_{l}$ is the induced $\spinc$-structure on $\Ycal$.
Let us denote by $p :\T(K\times_{T}\Ycal)\to K\times_{T}\Ycal$,   
$p_{_\Ycal}:\T\Ycal\to \Ycal$ and $p_{_{K/T}}:\T(K/T)\to K/T$ the 
canonical projections. Using (\ref{eq.2.reduction}), we see that
$$
p^{*}\Scal=\left( K\times_{T}p_{_\Ycal}^{*}\Scal(\Ycal) 
\right)\otimes\ p^{*}_{_{K/T}}\wedge_{\C}\T(K/T)\ ,
$$
where $\Scal(\Ycal)$ is the spinor bundle on $\Ycal$.
Hence we get the decomposition
$$
\sthom_{K}(K\times_{T}\Ycal)=\Thom_{K}(K/T)\odot 
K\times_{T}\sthom_{T}(\Ycal)\ .
$$
The transversally elliptic 
symbol $\Thom_{K}^{\Phi_{\tilde{\mu}}}(\Ucal)$ is equal to 
$$
\left[
\Thom_{K}(K/T)\odot \Thom_{K}(\overline{K/T})\odot K\times_{T}\sthom_{T}(\Ycal)
\right]_{{\rm deformed\ by}\ \Hcal^{\tilde{\mu}}} \ ,
$$
hence 
$\sigma_{1}:=(\phi^{-1})^{*}\Thom_{K}^{\Phi_{\tilde{\mu}}}(\Ucal)$ is 
equal to 
$$
\left[
\Thom_{K}(K/T)\odot K\times_{T}
\left(\Thom_{T}(\overline{K/T})\odot \sthom_{T}(\Ycal)\right)
\right]_{{\rm deformed\ by}\ \Hcal'} \ ,
$$
where $\Hcal^{'}=\phi_*(\Hcal^{\tilde{\mu}})$.

Using the decomposition 
$\T\Ucal'\simeq K\times_{T}(\kgot/\tgot\oplus
K\times_{T}(\overline{\kgot/\tgot})\oplus \T\Ycal)$, a small 
computation gives $\Hcal'(m)=pr_{\kgot/\tgot}(h\tilde{\mu})+ 
R(m) + \Hcal_{\tilde{\mu}}(y)+ S(m)$
for $m=[k;[h],y]\in \Ucal'$, where\footnote{A small computation
shows that $R(m)=pr_{\kgot/\tgot}(h^{-1}(pr_{\tgot}(h\tilde{\mu})-\Phi(y)))$, 
and $S(m)=[\tilde{\mu}-pr_{\tgot}(h\tilde{\mu})]_{\Ycal}(y)$.}
$R(m)\in \overline{\kgot/\tgot}$ and  $S(m)\in \T_{y}\Ycal$ vanishes
when $m\in K\times_{T}(\{\bar{e}\}\times\Ycal)$, i.e. when
$[h]=\bar{e}$. Here $\Hcal_{\tilde{\mu}}$ is the Hamiltonian vector field 
of the function $\frac{-1}{2}\parallel\Phi-\tilde{\mu}\parallel^{2}:\Ycal\to\R$. 

The transversally elliptic symbol $\sigma_{1}$ is equal to the exterior product
$$
\sigma_{1}(m, \xi_{1}+\xi_{2}+ v)=
\clif(\xi_{1}-pr_{\kgot/\tgot}(h\tilde{\mu}))\odot
\clif(\xi_{2}-R(m))\odot
\clif(v-\Hcal_{\tilde{\mu}}- S(m))\ ,
$$
with $\xi_{1}\in\kgot/\tgot$,   
$\xi_{2}\in\overline{\kgot/\tgot}$,  and $v\in\T\Ycal$.

Now, we simplify the symbol $\sigma_{1}$ without changing its 
$K$-theoretic class. Since $\Char(\sigma_{1})\cap\T_{K}\Ucal'
= K\times_{T}(\{\bar{e}\}\times\Ycal)$, we transform
$\sigma_{1}$ through the 
$K$-invariant diffeomorphism $h=e^{X}$ from a neighborhood 
of  $0$ in $\overline{\kgot/\tgot}$ to a neighborhood of $\bar{e}$ in 
$\overline{K/T}$. That gives $\sigma_{2}\in 
K_{K}(\T_{K}(K\times_{T}(\kgot/\tgot\times \Ycal)))$ defined by 
\begin{eqnarray*}
\lefteqn{\sigma_{2}([k;X,y], \xi_{1}+\xi_{2}+ v)
=}\\
& & 
\clif(\xi_{1}-pr_{\kgot/\tgot}(e^{X}\tilde{\mu}))\odot
\clif(\xi_{2}-R(m))\odot
\clif(v-\Hcal_{\tilde{\mu}}- S(m))\ .
\end{eqnarray*}    
Now trivial homotopies link $\sigma_{2}$ 
with the symbol 
$$
\sigma_{3}([k,X,y], \xi_{1}+\xi_{2}+ v)= 
\clif(\xi_{1}-[X,\tilde{\mu}])\odot
\clif(\xi_{2})\odot
\clif(v-\Hcal_{\tilde{\mu}})\ ,
$$
where we have removed the terms
$R(m)$ and $S(m)$, and where we have replaced 
$pr_{\kgot/\tgot}(e^{X}\tilde{\mu})=[X,\tilde{\mu}]+{\rm o}([X,\tilde{\mu}])$
by the term $[X,\tilde{\mu}]$. Now we see that $\sigma_{3}=i_{*}
( \sigma_{4})$ where the symbol 
$\sigma_{4}\in K_{T}(\T_{T}(\kgot/\tgot\times 
\Ycal))$ is defined by 
$$
\sigma_{4}(X,y;\xi_{2}+ v)= 
\clif(-[X,\tilde{\mu}])\odot \clif(\xi_{2})\odot \clif(v-\Hcal_{\tilde{\mu}})\ .
$$    
So $\sigma_{4}$ is equal to the exterior product of $(y,v)\to
\clif(v-\Hcal_{\tilde{\mu}})$, which is 
$\sthom_{T}^{\Phi-\tilde{\mu}}(\Ycal)$, with the 
transversally elliptic symbol on $\kgot/\tgot$:
$(X,\xi_{2})\to \clif(-[X,\tilde{\mu}])\odot \clif(\xi_{2})$. But 
the $K$-theoretic class of 
this former symbol is  equal to $k_{!}(\C)$, where $k:\{0\}
\croc\kgot/\tgot$ (see subsection 5.1 in \cite{pep4}). This shows that  
$$
\sigma_{4}=k_{!}(\C)\odot
\sthom_{T}^{\Phi-\tilde{\mu}}(\Ycal)=
j_{!}(\sthom_{T}^{\Phi-\tilde{\mu}}(\Ycal))\ .
$$
$\Box$

\subsection{Proof of Proposition \ref{prop.2}}\label{preuve.prop.2}

In this subsection, $\tilde{\mu}=\mu +\rho_{c}$ is fixed,  
and is assumed to belong to $\Delta$. 
The induced moment map on the symplectic slice $\Phi^{-1}({\rm 
Interior}(\tgot^*_+))$ is still denoted by $\Phi$. Let
$r>0$ be the smallest non zero critical value of 
$\parallel\Phi-\tilde{\mu}\parallel$, and let 
$\Ycal=\Phi^{-1}\{\xi\in {\rm Aff}(\Delta), 
\parallel\xi-\tilde{\mu}\parallel<\frac{r}{2}\}$.

For $\xi\in {\rm Aff}(\Delta)$, we consider $\xi_{t}=t\xi +(1-t)\tilde{\mu}$, 
$0\leq t\leq 1$. If one shows that there exists a compact subset $\Kcal\subset\Ycal$ 
such that $\Cr(\parallel\Phi-\xi_{t}\parallel^2)\cap\Ycal\subset 
\Kcal$, the family $\sthom_{T}^{\Phi-\xi_{t}}(\Ycal)$, $0\leq t\leq 1$, 
defines then  a homotopy of transversally elliptic symbols between 
$\sthom_{T}^{\Phi-\tilde{\mu}}(\Ycal)$ and
$\sthom_{T}^{\Phi-\xi}(\Ycal)$. It shows that \break
$\Qcal^{^T}_{\Phi-\tilde{\mu}}(\Ycal,-)$ and 
$\Qcal^{^T}_{\Phi-\xi}(\Ycal,-)$ are equal.

We describe now $\Cr(\parallel\Phi-\xi_{t}\parallel^2)\cap\Ycal$ using 
a parametrization introduced in \cite{pep1}[Section 6]. 
Let $\Bcal$ be the collection of affine subspaces of 
$\tgot^{*}$ generated by the image under $\Phi$ of submanifolds $\Zcal$ 
of the following type: $\Zcal$ is a connected component of $\Ycal^H$ 
which intersects $\Phi^{-1}(\tilde{\mu})$, $H$ being a subgroup of 
$T$. The set $\Bcal$ is finite since 
$\Phi^{-1}(\tilde{\mu})$ is compact and thus has a finite number of 
stabilizers for the $T$ action. Note that $\Bcal$ is 
reduced to ${\rm Aff}(\Delta)$ if $\tilde{\mu}$ is 
a generic quasi regular value of $\Phi$. For $A\in\Bcal$, 
we denote by $\beta(-,A)$ the orthogonal projection on $A$. Let 
$\Bcal_{\xi}=\{\beta(\xi,A)-\xi\mid A\in\Bcal\}$. Like in 
\cite{pep2}[Section 4.3], we see that 
\begin{equation}\label{eq.decomposition.delta}
  \Cr(\parallel\Phi-\xi\parallel^2)\cap\Ycal
=\bigcup_{\beta\in\Bcal_{\xi}}
\left(\Ycal^{\beta}\cap\Phi^{-1}(\beta +\xi)\right)
\end{equation}
if $\parallel\xi-\tilde{\mu}\parallel<\frac{r}{2}$. If we take 
$\Kcal:=\Phi^{-1}\{\xi\in {\rm Aff}\Delta, 
\parallel\xi-\tilde{\mu}\parallel\leq\frac{r}{3}\}$, we have \break
$\Cr(\parallel\Phi-\xi\parallel^2)\cap\Ycal\subset\Kcal$ for 
$\parallel\xi-\tilde{\mu}\parallel\leq\frac{r}{3}$. Thus point $i)$ 
of Proposition \ref{prop.2} is proved.

\medskip

Now we fix $\xi\in {\rm Aff}(\Delta)$  
close enough to $\tilde{\mu}$. And for each $\beta\in\Bcal_{\xi}$, we 
denote by $\Qcal^{^T}_{\beta}(\Ycal,-)$ the map localized near 
$\Ycal^{\beta}\cap\Phi^{-1}(\beta +\xi)$.  The excision property tells 
us, like in (\ref{localisation-pep}), that 
$$
\Qcal^{^T}_{\Phi-\xi}(\Ycal,-)=\sum_{\beta\in\Bcal_{\xi}}\Qcal^{^T}_{\beta}(\Ycal,-)\ .
$$
Note that $0\in\Bcal_{\xi}$ if and only if  $\Phi^{-1}(\xi)\neq\emptyset$.
Point iii) of Proposition \ref{prop.2} will follow from the following 

\begin{lem}\label{lem.2}  Let $\xi\in {\rm Aff}(\Delta)$  close enough to $\tilde{\mu}$, 
  and let $\beta$ be a non-zero element of $\Bcal_{\xi}$. Then 
  $[\Qcal^{^T}_{\beta}(\Ycal,\tilde{L}\vert_{\Ycal}\otimes \C)]^T=0$. Hence  
  $[\Qcal^{^T}_{\Phi-\xi}(\Ycal,\tilde{L}\vert_{\Ycal}\otimes 
  \C_{-\mu})]^T=0$, if $\xi\notin\Delta$.
\end{lem}

For the point $ii)$ of Proposition \ref{prop.2}, we also need the 

\begin{lem}\label{lem.1}
  If $\xi\in\Delta^o$ is a generic quasi regular value of $\Phi$, we have 
  \break
  $[\Qcal^{^T}_{0}(\Ycal,\tilde{L}\vert_{\Ycal}\otimes 
  \C_{-\mu})]^T=\Qcal(M_{\xi},o(M_{\xi}))$.
\end{lem}  

\medskip

Other versions of Lemmas \ref{lem.1} and \ref{lem.2}
are already known : in 
the $\spin$-case for an $S^1$-action by Vergne \cite{Vergne96}, and by the author 
\cite{pep4} when the $\spinc$-structure comes from an almost complex 
structure. 

We review briefly the arguments, since they work in the same way.  
We consider the $\spinc$ structure on $\Ycal$ defined in Proposition 
\ref{prop.spinc.induit.1}, that we twist by the line bundle 
$\tilde{L}\vert_{\Ycal}\otimes\C_{-\mu}$ : it defines a $\spinc$ 
structure ${\rm Q}^{\mu}$ on $\Ycal$ with canonical line bundle 
$\Lfibre^{\mu}:= L_{2\omega}\otimes\C_{-2\tilde{\mu}}$. We consider 
then the symbol $\sthom_{T,\mu}^{\Phi-\xi}(\Ycal)$ constructed with 
${\rm Q}^{\mu}$ (see Definition \ref{def.S-thom.loc}). For 
$\beta\in\Bcal_{\xi}$, 
the term $\Qcal^{^T}_{\beta}(\Ycal,\tilde{L}\vert_{\Ycal}\otimes 
\C_{-\mu})$ is by definition the $T$-index of 
$\sthom_{T,\mu}^{\Phi-\xi}(\Ycal)\vert_{\Ucal^{\beta}}$, where 
$\Ucal^{\beta}$ is a 
sufficiently small open  neighborhood of 
$\Ycal^{\beta}\cap\Phi^{-1}(\beta +\xi)$ in $\Ycal$.

\bigskip

{\em Proof of Lemma \ref{lem.1}.}  A neighborhood 
$\Ucal^0$ of $\Zcal:=\Phi^{-1}(\xi)$ is diffeomorphic to 
a neighborhood of $\Zcal$ in $\Zcal\times\overrightarrow{\Delta}$, 
where $\Phi-\xi:\Zcal\times\overrightarrow{\Delta}\to 
\overrightarrow{\Delta}$ is the projection to the second factor. 
Let $pr: \Zcal\times\overrightarrow{\Delta}\to \Zcal$ 
be the projection to the first factor. We 
still denote by $\Qcal^{\mu}$ the $\spinc$-structure on 
$\Zcal\times\overrightarrow{\Delta}$ equal to $pr^*(\Qcal^{\mu}\vert_{\Zcal})$.
We easily show that $\Qcal^{^T}_{\beta}(\Ycal,\tilde{L}\vert_{\Ycal}\otimes 
\C_{-\mu})$ is equal to the $T$-index of  
$\sigma_{\Zcal}=\sthom_{T,\mu}^{\Phi-\xi}(\Zcal\times\overrightarrow{\Delta})$. 
Let ${\rm Q}''$ be the reduction of $\Qcal^{\mu}\vert_{\Zcal}$ introduced in 
(\ref{diagram.Q.seconde}). Since 
$\Qcal^{\mu}\vert_{\Zcal}=\spinc_{2k}\times_{(\spinc_{2k'}\times\u_{2l'})}{\rm Q}''$, 
the bundle of spinors $\Scal$ over $\Zcal\times\overrightarrow{\Delta}$ 
decomposes as
$$
\Scal=pr^*\Big(\pi^*\Scal(M_{\xi})\otimes 
\Zcal\times\wedge(\tgot/\tgot_{\Delta}\otimes\C)\Big)\ .
$$
Here $\Scal(M_{\xi})$ is the bundle of spinors on $M_{\xi}$ induces 
by the $\spinc$-structure ${\rm Q}''/\u_{2l'}$, and $\pi:\Zcal\to 
M_{\xi}$ is the quotient map. Inside the trivial bundle 
$\Zcal\times(\tgot/\tgot_{\Delta}\otimes\C)$, we have identified 
$\Zcal\times(\tgot/\tgot_{\Delta}\otimes i\R)$ with the subspace of $\T \Zcal$ formed by 
the vector fields generated by the infinitesimal action of 
$\tgot/\tgot_{\Delta}$, and $\Zcal\times(\tgot/\tgot_{\Delta}\otimes\R)$ with
$\Zcal\times\overrightarrow{\Delta}\subset
\T(\Zcal\times\overrightarrow{\Delta})\vert_{\Zcal}$.
For $(z,f)\in\Zcal\times\overrightarrow{\Delta}$, let us decompose 
$v\in \T_{(z,f)}(\Zcal\times\overrightarrow{\Delta})$ as $v=v_1 +
X+i Y$, where $v_1\in\pi^*(\T M_{\xi})$, and 
$X+i Y\in \tgot/\tgot_{\Delta}\otimes\C$. The map 
$\sigma_{\Zcal}(z,f;v_1 +X+i Y)$ acts on $\Scal(M_{\xi})_z\otimes 
\wedge(\tgot/\tgot_{\Delta}\otimes\C)$ as the product
$$
\clif_z(v_1)\odot \clif(X+i(Y-f))\ ,
$$
which is homotopic\footnote{See \cite{pep4}[Section 6.1].} 
to the  transversally elliptic symbol 
$$
\ \clif_z(v_1)\odot \clif(f+iX)\ .
$$
So we have proved that $\sigma_{\Zcal}=j_!\circ\pi^*(\sthom(M_{\xi}))$, where 
$j_!: K_T(\T_T\Zcal)\to K_T(\T_T(\Zcal\times\overrightarrow{\Delta}))$ 
is induced by the inclusion $j:\Zcal\croc \Zcal\times\overrightarrow{\Delta}$. 
The last equality finishes the proof (see \cite{pep4}[Section 6.1]).

\bigskip

{\em Proof of Lemma \ref{lem.2}.}  The equality 
$[\Qcal^{^T}_{\beta}(\Ycal,\tilde{L}\vert_{\Ycal}\otimes 
\C_{-\mu})]^T=0$ comes from a localization formula on the submanifold 
$\Ycal^{\beta}$ for the map $\Qcal^{^T}_{\beta}(\Ycal,-)$
(see \cite{pep4,Vergne96}). The normal bundle $\Ncal$ 
of $\Ycal^{\beta}$ in $\Ycal$ carries a complex structure $J_{\Ncal}$ on the fibers 
such that each $\tore_{\beta}$-weight $\alpha$ on $(\Ncal,J_{\Ncal})$ 
satisfies $(\alpha,\beta)>0$. The principal bundle ${\rm Q}^{\mu}$, 
when restricted to $\Ycal^{\beta}$ admits the reduction

\begin{equation}\label{diagram.Q.beta}
\xymatrix@C=2cm{
 {\rm Q}'\ar[r]\ar[d] & 
 \Pso(\T \Ycal^{\beta})\times\Pu(\Ncal)
 \times\Pu(\Lfibre\vert_{\Ycal^{\beta}})\ar[d]\\
  {\rm Q}^{\mu}\vert_{\Ycal^{\beta}}\ar[r]   & 
  \Pso(\T \Ycal^{\beta}\oplus\Ncal)
  \times\Pu(\Lfibre^{\mu}\vert_{\Ycal^{\beta}})\ ,
  }
\end{equation}
Hence ${\rm Q}^{\beta}:={\rm Q}'/\u(l)$ is a $\spinc$-structure on 
$\Ycal^{\beta}$ with canonical line bundle equal to $\Lfibre^{\beta}:=
\Lfibre^{\mu}\vert_{\Ycal^{\beta}}\otimes(\det\Ncal)^{-1}= 
L_{2\omega}\vert_{\Ycal^{\beta}}\otimes\C_{-2\tilde{\mu}}\otimes(\det\Ncal)^{-1}$. 
Let $\Qcal^T_{\beta}(\Ycal^{\beta},-)$ be the map defined 
by ${\rm Q}^{\beta}$ and localized near 
$\Phi^{-1}(\beta+\xi)\cap\Ycal^{\beta}$ by $\Phi-\xi$. Following the 
argument of \cite{pep4}[Section 6] one obtains
$$
\Qcal^{^T}_{\beta}(\Ycal,\tilde{L}\vert_{\Ycal}\otimes 
\C_{-\mu})=(-1)^{l}\sum_{k\in\N}\Qcal^{^T}_{\beta}(\Ycal^{\beta},
\det\Ncal\otimes S^k(\Ncal))\ ,
$$
where $S^k(\Ncal)$ is the k-th symmetric 
product of $\Ncal$, and $l={\rm rank}_{\C}\Ncal$. Thus, it is sufficient to prove that 
$[\Qcal^{^T}_{\beta}(\Ycal^{\beta},\det\Ncal\otimes S^k(\Ncal))]^T=0$ for 
every $k\in\N$. For this purpose, we use Lemma \ref{lem.poids.tore}. 
Let $\alpha$ be the $\tore_{\beta}$-weight on $\det\Ncal$. From the 
Kostant formula (\ref{eq.kostant}), the $\tore_{\beta}$-weight on 
$L_{2\omega}\vert_{\Ycal^{\beta}}$ is equal to $2(\beta +\xi)$. Hence 
any $\tore_{\beta}$-weight $\gamma$ on 
$\det\Ncal\otimes S^k(\Ncal)\otimes(\Lfibre^{\beta})^{1/2}$ is of the form 
$$
\gamma=\beta+\xi-\tilde{\mu}+\frac{1}{2}\alpha + \delta
$$
where $\delta$ is a $\tore_{\beta}$-weight on $S^k(\Ncal)$. So 
$(\gamma,\beta)=(\beta+\xi-\tilde{\mu},\beta)+\frac{1}{2}(\alpha,\beta) + 
(\delta,\beta)$. But the $\tore_{\beta}$-weights on $\Ncal$ are 
`positive' for $\beta$, so $(\alpha,\beta)>0$ and $(\delta,\beta)\geq 0$. On the 
other hand, $\beta +\xi=\beta(\xi,A)$ is the orthogonal projection of 
$\xi$ on some affine subspace $A\subset\tgot^*_+$ which contains 
$\tilde{\mu}$: hence $(\beta+\xi-\tilde{\mu},\beta)=0$. This proves 
that $(\gamma,\beta)>0$. $\Box$


\section{Quantization and the discrete series}

In this section we apply Theorem \ref{quantization-non-compact-1} 
to the coadjoint orbits that parametrize the discrete series of a real, 
connected, semi-simple Lie group $G$, with finite center.  Nice 
references on the subject of `the discrete series' are \cite{Schmid97,Duflo77}.

Let $K$ be a maximal compact subgroup of $G$, and $T$ be a 
maximal torus in $K$. For the remainder of this section, we 
assume that $T$ is a Cartan subgroup of $G$. The discrete series of 
$G$ is then non-empty and is parametrized by a subset $\widehat{G}_{d}$ in 
the dual $\tgot^*$ of the Lie algebra of $T$ \cite{Harish-Chandra65-66}.

Let us fix some notation.
Let $\Rgot_c\subset\Rgot\subset\Lambda^{*}$ be respectively the set of 
(real) roots for the action of $T$ on $\kgot\otimes\C$ and 
$\ggot\otimes\C$. We choose a system of positive roots $\Rgot_{c}^{+}$ for 
$\Rgot_{c}$, we denote by $\tgot^{*}_{+}$ the corresponding Weyl chamber, 
and we let $\rho_{c}$ be half the sum of the elements of $\Rgot_{c}^{+}$. 
We denote by $B$ the Killing form on $\ggot$. It induces a scalar product 
(denoted by $(-,-)$) on $\tgot$, and then on $\tgot^{*}$. 
An element $\lambda\in \tgot^{*}$ is called {\em regular} if 
$(\lambda,\alpha)\neq 0$ for every $\alpha\in\Rgot$, or equivalently,  
if the stabilizer subgroup of $\lambda$ in $G$ is $T$. 
Given a system of positive roots $\Rgot^{+}$ for $\Rgot$, consider 
the subset $\Lambda^{*}+\frac{1}{2}\sum_{\alpha\in\Rgot^{+}}\alpha$ of
$\tgot^*$. This does not depend on the choice of $\Rgot^{+}$, and we denote it 
by $\Lambda^{*}_{\rho}$ \cite{Duflo77}. 
 
\medskip 
 
The discrete series of $G$ are parametrized by 
\begin{equation}\label{G-hat}
    \widehat{G}_{d}:=\{\lambda\in\tgot^{*},\lambda\ {\rm regular}\ 
    \}\cap \Lambda^{*}_{\rho}\cap \tgot^{*}_{+}\ .
\end{equation}
When $G=K$ is compact, the set $\widehat{G}_{d}$ equals
$\Lambda^{*}_{+}+\rho_{c}$, and it parametrizes the  set of irreducible 
representations of $K$. Harish-Chandra has associated to any $\lambda\in \widehat{G}_{d}$ 
an invariant eigendistribution on $G$, denoted by $\Theta_{\lambda}$, 
which is shown to be the global trace of an irreducible, square integrable, unitary 
representation of $G$.

On the other hand we associate to $\lambda\in\widehat{G}_{d}$, 
the regular coadjoint orbit $M:=G\cdot \lambda$. It is a symplectic 
manifold with a Hamiltonian action of $K$. Since the vectors 
$X_{M},\, X\in \ggot$, span the tangent space at every $\xi\in M$, the 
symplectic $2$-form is determined by 
$$
\omega(X_{M},Y_{M})_{\xi}=\langle \xi,[X,Y]\rangle\ .
$$
The corresponding moment map 
$\Phi:M\to\kgot^*$ for the $K$-action
is the composition of the inclusion
$i: M\croc\ggot^*$ with the projection 
$\ggot^*\to\kgot^*$. The vector $\lambda$ determines a choice 
$\Rgot^{+,\lambda}$ of positive roots for the $T$-action on 
$\ggot\otimes\C$ : $\alpha\in \Rgot^{+,\lambda}\Longleftrightarrow 
(\alpha,\lambda)>0$. We recall now how the choice of $\Rgot^{+,\lambda}$ determines a complex 
structure on $M$. First take the decomposition
$\ggot\otimes\C=\tgot\otimes\C\,\oplus\,\sum_{\alpha\in\Rgot}\ggot_{\alpha}$ 
where $\ggot_{\alpha}:=\{v\in\ggot\otimes\C \mid \exp(X).v=
e^{i\langle\alpha,X\rangle}v\ {\rm for\ any}\ X\in\tgot$\}. 
It gives the following $T$-equivariant 
decomposition of the complexified tangent space of $M$ at $\lambda$ :
$$
\T_{\lambda}M\otimes\C=\sum_{\alpha\in\Rgot}\ggot_{\alpha}=\ngot\oplus
\overline{\ngot}\ ,
$$
with $\ngot=\sum_{\alpha\in\Rgot^{+,\lambda}}\ggot_{\alpha}$. We have then a 
$T$-equivariant isomorphism $\Ical:\T_{\lambda}M\to \ngot$ equal to the 
composition of the inclusion 
$\T_{\lambda}M\croc\T_{\lambda}M\otimes\C$ with the projection 
$\ngot\oplus\overline{\ngot}\to \ngot$. The $T$-equivariant complex structure 
$J_{\lambda}$ on $\T_{\lambda}M$ is determined by the relation 
$\Ical(J_{\lambda}v)=i\,\Ical(v)$. Hence, the set of real 
infinitesimal weights for the $T$-action on 
$(\T_{\lambda}M,J_{\lambda})$ is $\Rgot^{+,\lambda}$. Since $M$ is 
a homogeneous space, $J_{\lambda}$ defines  an invariant almost 
complex structure $J$ on $M$, which is in fact integrable.
Through the isomorphism $M\cong G/T$, the canonical line bundle 
$\kappa=\det_{\C}(\T M)^{-1}$ is equal to
$\kappa=G\times_{T}\C_{-2\rho}$ with 
$\rho=\frac{1}{2}\sum_{\alpha\in\Rgot^{+,\lambda}}\alpha$. 

If $\lambda\in \widehat{G}_{d}$, then $\lambda-\rho$ is a weight, and 
\begin{equation}\label{L.tilde.G.lambda}
\tilde{L}:=G\times_{T}\C_{\lambda-\rho}\to G/T
\end{equation}
is a $\kappa$-prequantum line bundle over $(M,\omega,J)$. We have shown in 
\cite{pep3}, that \break $\Cr(\parallel\Phi\parallel^2)$ is compact,  
equal to the $K$-orbit $K\cdot\lambda$. Then the generalized 
Riemann-Roch character $RR_{\Phi}^{^K}(M,-)$ is well defined 
(see Definition \ref{def.RR.phi}). The main result of this section is 
the following 

\medskip

\begin{theo}\label{th.theta=RR}
  We have the following equality of tempered distributions on $K$
  $$
  \Theta_{\lambda}\vert_{K}=(-1)^{\frac{\dim(G/K)}{2}}\, 
  RR_{\Phi}^{^K}(G\cdot \lambda,\tilde{L})\ ,
  $$
  where $\Theta_{\lambda}\vert_{K}$ is the restriction of the 
  eigendistribution $\Theta_{\lambda}$ to the subgroup $K$.
\end{theo}

\medskip

The proof of Theorem \ref{th.theta=RR} is given in Subsection 
\ref{subsec.preuve.theoreme}. It uses the Blattner formulas 
in an essential way (see Subsection \ref{subsec.K.multiplicite}).

With Theorem \ref{th.theta=RR} at our disposal we can exploit 
the result of Theorem \ref{quantization-non-compact-1} 
to compute the $K$-multiplicities, $\mm_{\mu}(\lambda)\in\N$, of 
$\Theta_{\lambda}\vert_{K}$ in term of the reduced spaces. 
By definition we have
\begin{equation}\label{pi-lambda-K}
\Theta_{\lambda}\vert_{K}=\sum_{\mu\in\Lambda^{*}_{+}}
\mm_{\mu}(\lambda)\, \chi_{_{\mu}}^{_K} \quad {\rm in}  \quad
R^{-\infty}(K)\ .
\end{equation}

The moment map $\Phi:G\cdot \lambda\to\kgot^{*}$ is {\em proper} since 
the $G\cdot \lambda$ is closed in $\ggot^*$ \cite{pep3}. We show (Lemma \ref{forme-normale}) 
that the moment polyhedron $\Delta=\Phi(G\cdot \lambda)\cap\tgot^*_+$ 
is of dimension $\dim T$. Thus on the relative interior $\Delta^o$ of the 
moment polyhedron, the notions of {\em generic quasi-regular values} 
and {\em regular values} coincide : they concern the elements 
$\xi\in\Delta^o$ such that 
$\Phi^{-1}(\xi)$ is a smooth submanifold with a locally free action 
of $T$. We have shown (Subsection \ref{def.Q.M.mu}) how to define the quantity 
$Q((G\cdot \lambda)_{\mu+\rho_c})\in\Z$ as the index of 
a suitable $\spinc$ Dirac operator on $\Phi^{-1}(\xi)/T$, where 
$\xi\in\Delta^o$ is a regular value of $\Phi$ close enough to $\mu+\rho_c$.  

\begin{prop}\label{prop.K.multiplicite}
 For every $\mu\in\Lambda^*_+$, we have 
 $$
 \mm_{\mu}(\lambda)= \Qcal\left((G\cdot \lambda)_{\mu+\rho_c}\right)\ .
 $$
In particular $\mm_{\mu}(\lambda)=0$ if $\mu+\rho_c$ does not
belong to the relative interior of the moment polyhedron $\Delta$. 
\end{prop} 

\medskip

{\em Proof.}  A small check of orientations shows that 
$\esp=(-1)^{\frac{\dim(G/K)}{2}}$, thus this proposition follows from 
Theorems \ref{quantization-non-compact-1} and \ref{th.theta=RR} if one checks that 
the following holds:  $(G\cdot \lambda,\Phi)$ satisfies Assumption 
\ref{hypothese.phi.t.carre}, and the infinitesimal $K$-stabilizers are Abelian.
The first point will be handled in Subsection \ref{subsec.hyp.phi.2}. 
The second point is obvious since $M\cong G/T$: all the $G$-stabilizers are 
conjugate to $T$, so all the $K$-stabilizers are Abelian. $\Box$  

\bigskip


\subsection{Blattner formulas}\label{subsec.K.multiplicite}

In this section, we fix  $\lambda\in \widehat{G}_{d}$. Let 
$\Rgot^{+,\lambda}$ be the system of positive roots defined $\lambda$: 
$\alpha\in \Rgot^{+,\lambda}\Longleftrightarrow 
(\alpha,\lambda)>0$. Then $\Rgot^{+}_{c}
\subset\Rgot^{+,\lambda}$, and $\rho=\frac{1}{2}\sum_{\alpha\in\Rgot^{+,\lambda}}
\alpha$ decomposes in $\rho=\rho_{c}+\rho_{n}$ where
$\rho_{n}=\frac{1}{2}\sum_{\alpha\in\Rgot^{+}_{n}}\alpha$ and
$\Rgot^{+}_{n}=\Rgot^{+,\lambda}-\Rgot^{+}_{c}$.

Let $\Pcal:\Lambda^{*}\mapsto \N$ be the partition function associated 
to the set $\Rgot^{+}_{n}$ : for $\mu\in \Lambda^{*}$, $\Pcal(\mu)$ 
is the number of distinct ways we can write $\mu=
\sum_{\alpha\in\Rgot^{+}_{n}}k_{\alpha}\alpha$ with 
$k_{\alpha}\in\N$ for all $\alpha$. The following Theorem is known as the 
Blattner formulas and was first proved by Hecht and Schmid \cite{Hecht-Schmid}.

\medskip

\begin{theo}\label{Blattner}
    For $\mu\in\Lambda^{*}_{+}$, we have 
$$
\mm_{\mu}(\lambda)=\sum_{w\in 
W}(-1)^w\Pcal\Big(w(\mu+\rho_{c})-(\mu_{\lambda}+\rho_{c})\Big)\ ,
$$
where\footnote{We shall note that $\mu_{\lambda}\in\Lambda^{*}_{+}$ 
(see \cite{Duflo77}, section 5).} 
$\mu_{\lambda}:=\lambda-\rho_{c}+\rho_{n}$. Here $W$ is the Weyl group 
of $(K,T)$.
\end{theo}

\medskip

Using Theorem \ref{Blattner}, we can describe 
$\Theta_{\lambda}\vert_{K}$ through the holomorphic induction map 
$\HolT: R^{-\infty}(T)\to R^{-\infty}(K)$. Recall that $\HolT$
is characterized by the following properties:
i) $\HolT(t^{\mu})=\chi_{_{\mu}}^{_K}$ for every dominant weight 
$\mu\in\Lambda^{*}_{+}$; 
ii) $\HolT(t^{w\circ\mu})=(-1)^w\HolT(t^{\mu})$ for every $w\in W$ and 
$\mu\in\Lambda^{*}$; 
iii) $\HolT(t^{\mu})=0$ if $W\circ\mu\cap \Lambda^{*}_{+}=\emptyset$. 
Using these properties we have 
\begin{equation}\label{eq.induit.i}
\sum_{\mu\in\Lambda^{*}}\Rcal(\mu)\, 
\HolT(t^{\mu})=\sum_{\mu\in\Lambda^{*}_{+}}
\Big[\sum_{w\in W}(-1)^w\Rcal(w\circ\mu)\Big]\, 
\chi_{_{\mu}}^{_K}\ ,
\end{equation}
for every map $\Rcal:\Lambda^{*}\to \Z$.

For a weight $\alpha\in\Lambda^{*}$, with $(\lambda,\alpha)\neq 0$, 
let us denote the oriented inverse of $(1-t^{\alpha})$
in the following way 
$$
[1-t^{\alpha}]^{-1}_{\lambda}=\left\{
\begin{array}{ll}
\sum_{k\in\N}t^{k\alpha}\ ,              & {\rm if} \quad (\lambda,\alpha)> 0 \\  
-t^{-\alpha}\sum_{k\in\N}t^{-k\alpha}\ , & {\rm if} \quad (\lambda,\alpha)<0\ .
\end{array}\right.
$$
Let $A=\{\alpha_{1},\cdots,\alpha_{l}\}$ 
be a set of weights with $(\lambda,\alpha_{i})\neq 0,\,\forall i$.
We denote by $A^{+}=\{\esp_{1}\alpha_{1},\cdots,\esp_{l}\alpha_{l}\}$ 
the corresponding set of polarized weights: $\esp_{i}=\pm 1$ and
$(\lambda,\esp_{i}\alpha_{i})> 0$ for all $i$. 
The product $\Pi_{\alpha\in A} [1-t^{\alpha}]^{-1}_{\lambda}$
is well defined in $R^{-\infty}(T)$, and is denoted by 
$[\Pi_{\alpha\in A}(1-t^{\alpha})]^{-1}_{\lambda}$. A small
computation shows that
\begin{eqnarray}\label{produit-infini} 
\Big[\Pi_{\alpha\in A}(1-t^{\alpha})\Big]^{-1}_{\lambda}
&=& (-1)^{r}\,t^{-\gamma}\Big[\Pi_{\alpha\in A^{+}}
(1-t^{\alpha})\Big]^{-1}_{\lambda} \nonumber\\
&=&
(-1)^{r}\,t^{-\gamma}\sum_{\mu\in\Lambda^{*}}\Pcal_{A^+}(\mu)\, 
t^{\mu}\ .
\end{eqnarray}
Here $\Pcal_{A^+}:\Lambda^{*}\mapsto \N$ is the partition function 
associated to $A^+$, $\gamma=\sum_{(\lambda,\alpha)<0}\alpha$, and 
$r=\sharp\{\alpha\in A,\, (\lambda,\alpha)<0\}$. This notation is 
compatible with the one we used in \cite{pep4}[Section 5]. If $V$ is 
a complex $T$-vector space where the subspace fixed by $\lambda$ is
reduced to $\{0\}$, then $\wedge^{\bullet}_{\C}V\,\in R(T)$ admits a 
polarized inverse $[\wedge^{\bullet}_{\C}V]^{-1}_{\lambda}=
[\Pi_{\alpha\in \Rgot(V)}(1-t^{\alpha})]^{-1}_{\lambda}$, where 
$\Rgot(V)$ is the set of real infinitesimal $T$-weights on $V$.

\medskip

\begin{lem}\label{lem.theta.RR} 
  We have the following equality in $R^{-\infty}(K)$
 $$
 \Theta_{\lambda}\vert_{K}=\HolT\Big(t^{\mu_{\lambda}}
 \Big[\Pi_{\alpha\in \Rgot_{n}^{+}}(1-t^{\alpha})\Big]^{-1}_{\lambda}\Big)\ .
 $$
\end{lem}

\medskip

{\em Proof.}  Let $\Theta\in R^{-\infty}(K)$ be the RHS in the 
equality of the Lemma. From (\ref{produit-infini}), we have 
$\Theta=\sum_{\mu\in\Lambda^{*}}\Pcal(\mu)\, 
\HolT(t^{\mu+\mu_{\lambda}})=\sum_{\mu\in\Lambda^{*}}\Pcal(\mu-\mu_{\lambda})\, 
\HolT(t^{\mu})$.  If we use now (\ref{eq.induit.i}), we see that multiplicity of
$\Theta$ relative to the highest weight $\mu\in\Lambda^{*}_{+}$ 
is $\sum_{w\in W}(-1)^w\Pcal(w(\mu +\rho_{c})-(\mu_{\lambda}+\rho_{c}))$. 
From Theorem \ref{Blattner} we conclude that 
$\Theta_{\lambda}\vert_{K}=\Theta$. $\Box$

\medskip


\subsection{ Proof of Theorem 
\ref{th.theta=RR}}\label{subsec.preuve.theoreme}

In Lemma \ref{lem.theta.RR} we have used the Blattner formulas to 
write $\Theta_{\lambda}\vert_{K}$ in term of the holomorphic induction 
map $\HolT$. Theorem \ref{th.theta=RR} is then proved if one shows 
that
$RR_{\Phi}^{^K}(G\cdot \lambda,\tilde{L})=(-1)^r \HolT(t^{\mu_{\lambda}}
[\Pi_{\alpha\in \Rgot_{n}^{+}}(1-t^{\alpha})]^{-1}_{\lambda})$, with
$\mu_{\lambda}=\lambda-\rho_c+\rho_n$, and 
$r=\frac{1}{2}\dim(G/K)$. More generally, we show in this section 
that for any $K$-equivariant vector bundle $V\to G\cdot \lambda$ 
\begin{equation}\label{eq.theta=RR}
RR_{\Phi}^{^K}(G\cdot \lambda,V)= (-1)^r\HolT\Big(V_{\lambda}.t^{2\rho_n}. 
\Big[\Pi_{\alpha\in \Rgot_{n}^{+}}(1-t^{\alpha})\Big]^{-1}_{\lambda}\Big)\ ,
\end{equation}
where $V_{\lambda}\in R(T)$ is the fiber of $V$ at $\lambda$.

\medskip

First we recall why $\Cr(\parallel\Phi\parallel^{2})=K\cdot\lambda$ in
$M:=G\cdot\lambda$ (see \cite{pep3} for the general case of closed 
coadjoint orbits). One can work with an adjoint orbit 
$M:=G\cdot\tilde{\lambda}$  
through the $G$-identification $\ggot^{*}\simeq\ggot$ 
given by the Killing form; then $\Phi:M\to\kgot$ is just the 
restriction on $M$ of the (orthogonal) projection $\ggot\to\kgot$. Let $\pgot$ be 
the orthogonal complement of $\kgot$ in $\ggot$. Every $m\in M$ 
decomposes as $m=x_{m}+y_{m}$, with $x_{m}=\Phi(m)$ and $y_{m}\in\pgot$. 
The Hamiltonian vector field of $\frac{-1}{2}\parallel\Phi\parallel^{2}$ is
$\Hcal_{m}=[x_{m},m]=[x_{m},y_{m}]$ (see \ref{def.H}). Thus
$$
\Cr(\parallel\Phi\parallel^{2})=\{\Hcal=0\}=
\{m\in M,\ [x_{m},y_{m}]=0\}\ .
$$
Now, since $\tilde{\lambda}$ is elliptic, every $m\in M$ is also 
elliptic. If $m\in\Cr(\parallel\Phi\parallel^{2})$, $[m,x_{m}]=0$ and 
$m,x_{m}$ are elliptic, hence $y_{m}=m-x_{m}$ is elliptic and so is 
equal to $0$. Finally $\Cr(\parallel\Phi\parallel^{2})=G\cdot\tilde{\lambda}
\cap\kgot =K\cdot\tilde{\lambda}$. 

\medskip

According to Definition \ref{def.RR.phi}, the computation of $RR_{\Phi}^{^K}(G\cdot 
\lambda,\tilde{L})$ holds on a small $K$-invariant 
neighborhood of $K\cdot\lambda$ of $G\cdot\lambda$.  
Our model for the computation will be 
$$
\tilde{M}:=K\times_T\pgot
$$ 
endowed with the canonical $K$-action. The tangent bundle $\T\tilde{M}$ is 
isomorphic to  $K\times_T(\rgot\oplus\T\pgot)$ where $\rgot$ is the 
$T$-invariant complement of $\tgot$ in $\kgot$. One has a  symplectic form 
$\tilde{\Omega}$ on $\tilde{M}$ defined by $\tilde{\Omega}_{m}(V,V')=
\langle\lambda,[X,X']+[v,v']\rangle$. Here $m=[k,x]\in K\times_T\pgot$, and 
$V=[k,x;X+v]$, $X'=[k,x;X'+v']$ are two tangent vectors, with $X,X'\in\rgot$ and 
$v,v'\in\pgot$. A small computation shows that the $K$-action on 
$(K\times_T\pgot,\tilde{\Omega})$ is Hamiltonian with 
moment map $\tilde{\Phi}:\tilde{M}\to\kgot^*$ defined by 
$$
\tilde{\Phi}([k,x])=k\cdot\left(\lambda-\frac{1}{2}pr_{\tgot^*}(\lambda
\circ{\rm ad}(x)\circ{\rm ad}(x))\right)\ .
$$
Here ${\rm ad}(x)$ is the adjoint action of $x$, and 
$pr_{\tgot^*}:\ggot^*\to\tgot^*$ is the projection. Note first that 
the tangent space $\T_{\lambda}M$ and $\T_{[1,0]}\tilde{M}$ are 
canonically isomorphic to $\rgot\oplus\pgot$.

\begin{lem}\label{forme-normale} There exists a $K$-Hamiltonian 
isomorphism $\Upsilon:\Ucal\simeq\tilde{\Ucal}$, where $\Ucal$ is a 
$K$-invariant neighborhood of $K\cdot\lambda$ in $M$, $\tilde{\Ucal}$ is a 
$K$-invariant neighborhood of $K/T$ in $\tilde{M}$, and 
$\Upsilon(\lambda)=[1,0]$. We impose furthermore that the differential 
of $\Upsilon$ at $\lambda$ is the identity. 
\end{lem}

\begin{coro}\label{coro.forme-normale}
   The cone $\lambda+\sum_{\alpha\in\Rgot_n^+}\R^+\alpha$ coincides 
   with $\Delta=(G\cdot\lambda)\cap\tgot^*_+$ in a neighborhood of 
   $\lambda$. The polyhedral set 
   $\Delta$ is of dimension $\dim T$. 
\end{coro}

{\em Proof of the Corollary.}  The first assertion is an immediate consequence of Lemma 
\ref{forme-normale} and of the convexity Theorem \cite{L-M-T-W}. 
Let $X_o\in \tgot$ such that $\xi(X_o)=0$ for all 
$\xi\in\overrightarrow{\Delta}$, that is $\alpha(X_o)=0$ for all 
$\alpha\in\Rgot_n^+$ : $X_o$ commutes with all elements in $\pgot$. 
Let $\agot$ be a maximal Abelian subalgebra of $\pgot$, and let 
$\Sigma$ be the set of weights for the adjoint action of $\agot$ on 
$\ggot$ : $\ggot=\sum_{\alpha\in\Sigma}\ggot_{\alpha}$, where 
$\ggot_{\alpha}=\{Z\in\ggot,\,[X,Z]=\alpha(X)Z\ {\rm for \ all\ 
}X\in\agot\}$. Since $[X_o,\agot]=0$, we have 
$[X_o,\ggot_{\alpha}]\subset\ggot_{\alpha}$ for all 
$\alpha\in\Sigma$. But since $[X_o,\pgot]=0$ and 
$\ggot_{\alpha}\cap\kgot=0$ for all $\alpha\neq 0$, we see that 
$[X_o,\ggot_{\alpha}]=0$ for all $\alpha\neq 0$. But $X_o$ belongs to 
the Abelian subalgebra $\ggot_{0}$, so $[X_o,\ggot_{\alpha}]=0$. We 
see finally that $X_o$ belongs to the center of $\ggot$, and 
that implies $X_o=0$ since $G$ has a {\em finite center}. We have 
proved that $\overrightarrow{\Delta}^{\perp}=0$, or equivalently 
$\overrightarrow{\Delta}=\tgot$.  $\Box$

\medskip

{\em Proof of Lemma \ref{forme-normale}.}  The {\em symplectic cross-section Theorem} 
\cite{Guillemin-Sternberg84} asserts that the pre-image 
$\Ycal:=\Phi^{-1}({\rm interior}(\tgot^*_+))$ is a symplectic 
submanifold provided with a Hamiltonian action of $T$. The 
restriction $\Phi\vert_{\Ycal}$ is the moment map for the $T$-action 
on $\Ycal$. Moreover, the set $K.\Ycal$ is a $K$-invariant 
neighborhood of $K\cdot\lambda$ in $M$ diffeomorphic to 
$K\times_T\Ycal$. Since $\lambda$ is a fixed $T$-point of $\Ycal$, a 
Hamiltonian model for $(\Ycal,\omega\vert_{\Ycal},\Phi\vert_{\Ycal})$ 
in a neighborhood of $\lambda$ is 
$(\T_{\lambda}\Ycal,\omega_{\lambda},\Phi_{\lambda})$ 
where $\omega_{\lambda}$ is the linear symplectic form of the tangent 
space $\T_{\lambda}M$ restricted to $\T_{\lambda}\Ycal$, and 
$\Phi_{\lambda}:\T_{\lambda}\Ycal\to\tgot^*$ is the unique 
moment map with $\Phi_{\lambda}(0)=\lambda$. A small computation shows 
that $x\to\lambda\circ{\rm ad}(x)$ is an isomorphism from $\pgot$ to
$\T_{\lambda}\Ycal$, and $\Phi_{\lambda}(x)=\lambda -\frac{1}{2}pr_{\tgot^*}(\lambda
\circ{\rm ad}(x)\circ{\rm ad}(x))$. $\Box$

\medskip

We still denote the almost complex structure transported on 
$\tilde{\Ucal}\subset K\times_T\pgot$ through $\Upsilon$ by $J$. Since 
$d\Upsilon(\lambda)$ is the identity, 
$J_{[1,0]}:\rgot\oplus\pgot\to\rgot\oplus\pgot$ is equal to $J_{\lambda}$. Let 
$\pi:K\times_T\pgot\to K/T$, and $\pi_{\tilde{\Ucal}}:\tilde{\Ucal}\to K/T$ 
be the fibering maps. Remark that for any 
equivariant vector bundle $V$ over $M$ the vector bundle 
$(\Upsilon^{-1})^*(V\vert_{\Ucal})\to \tilde{\Ucal}$ is isomorphic to 
$\pi_{\tilde{\Ucal}}^*(K\times_{T}V_{\lambda})$. At this 
stage, we have according to Definition \ref{def.RR.phi}
\begin{equation}\label{eq.RR.tilde.1}
RR_{\Phi}^{^K}(G\cdot \lambda,V)=
    \indice_{\tilde{\Ucal}}^K\left(\Thom_{_{K}}^{\tilde{\Phi}}(\tilde{\Ucal},J)\otimes 
    \pi_{\tilde{\Ucal}}^*(K\times_{T}V_{\lambda})\right)\ .
\end{equation}    

With the help of Lemma \ref{lem.deformation}, we  define now a simpler 
representative of the class defined by  $\Thom_{_{K}}^{\tilde{\Phi}}(\tilde{\Ucal},J)$ 
in $\K_{K}(\T_{K}\tilde{\Ucal})$. Consider the map
\begin{eqnarray*}
\underline{\lambda}:K\times_T\pgot&\longrightarrow&\kgot^*\\ 
 (k,x) &\longmapsto& k\cdot\lambda\ ,
\end{eqnarray*}
and let $\underline{\lambda}_{\tilde{M}}$ be 
the vector field on $\tilde{M}$ generated  $\underline{\lambda}$ (see 
(\ref{def.H})). Note that $\underline{\lambda}_{\tilde{M}}$ never 
vanishes outside the zero section of $K\times_T\pgot$. Let $(-,-)_{_{\tilde{M}}}$ 
be the Riemannian metric on $\tilde{M}$ defined by 
$(V,V')_{_{\tilde{M}}}=(X,X')+(v,v')$ for 
$V=[k,x;X+v]$, $V=[k,x;X'+v']$. A small computation shows that 
 $$
 (\tilde{\Hcal},\underline{\lambda}_{\tilde{M}})_{_{\tilde{M}}}
 =\parallel\underline{\lambda}_{\tilde{M}}\parallel^2
 +\,o(\parallel\underline{\lambda}_{\tilde{M}}\parallel^2)
 $$ 
in the neighborhood of the zero section in $K\times_T\pgot$. 
Hence, if we take $\tilde{\Ucal}$ small enough, 
$(\tilde{\Hcal},\underline{\lambda}_{\tilde{M}})_{_{\tilde{M}}}>0$ on
$\tilde{\Ucal}-\{{\rm zero\ section}\}$, hence 
$\Thom_{_{K}}^{\tilde{\Phi}}(\tilde{\Ucal},J)=
\Thom_{_{K}}^{\underline{\lambda}}(\tilde{\Ucal},J)$ in 
$\K_{K}(\T_{K}\tilde{\Ucal})$ (see Lemma
\ref{lem.deformation}).

\medskip

Let $\tilde{J}$ be the $K$-invariant 
almost complex structure on $\tilde{M}$, constant on the fibers of 
$\tilde{M}\to K/T$, and equal to $J_{\lambda}$ at $[1,0]$ so that for 
$[k,x]\in K\times_T\pgot$, $\tilde{J}_{[k,x]}(V)=$ $[k,x,J_{\lambda}(X+v)]$ for 
$V=[k,x,X+v]$. Since the set $\{\underline{\lambda}_{\tilde{M}}=0\}=K/T$ 
is compact, using $\tilde{J}$ and the map $\underline{\lambda}$, one defines  
the localized Thom symbol
$$
\Thom_{_{K}}^{\underline{\lambda}}(\tilde{M},\tilde{J})
\in \K_{K}(\T_{K}\tilde{M})\ .
$$
Through the canonical identification of the tangent 
spaces at $[k,x]$ and $[k,0]$, one can write 
$\tilde{J}_{[k,x]}=\tilde{J}_{[k,0]}=J_{[k,0]}$ for any 
$[k,x]\in\tilde{\Ucal}$. We note that $J$ and $\tilde{J}$ are related on 
$\tilde{\Ucal}$ by the  homotopy $J^t$ of almost complex structures: 
$J^t_{[k,x]}:=J_{[k,tx]}$ for $[k,x]\in\tilde{\Ucal}$. From   
Lemma \ref{lem.deformation}, we conclude that the localized Thom 
symbols $\Thom_{_{K}}^{\underline{\lambda}}(\tilde{\Ucal},J)$ and
$\Thom_{_{K}}^{\underline{\lambda}}(\tilde{M},\tilde{J})\vert_{\tilde{\Ucal}}$
define the same class in $\K_{K}(\T_{K}\tilde{\Ucal})$, thus 
 (\ref{eq.RR.tilde.1}) becomes
\begin{equation}\label{eq.RR.tilde.2}
RR_{\Phi}^{^K}(G\cdot \lambda,V)=
    \indice_{\tilde{M}}^K\left(\Thom_{_{K}}^{\underline{\lambda}}
    (\tilde{M},\tilde{J})\otimes 
    \pi^*(K\times_{T}V_{\lambda})\right)\ .
\end{equation} 

\medskip

In order to compute (\ref{eq.RR.tilde.2}), we use the induction
morphism 
$$
i_*:\K_{T}(\T_{T}\pgot)\longrightarrow\K_{K}(\T_{K}(K\times_T\pgot))
$$ 
defined by Atiyah in \cite{Atiyah74} (see \cite{pep4}[Section 3]). 
The map $i_*$ enjoys two properties: first, $i_*$ is an isomorphism 
and the $K$-index of $\sigma\in\K_{K}(\T_{K}(K\times_T\pgot))$ 
can be computed with the $T$-index of $(i_*)^{-1}(\sigma)$.

Let $\sigma : p^*(E^+)\to p^*(E^-)$ be a $K$-transversally elliptic 
symbol on $K\times_T\pgot$, where $p:\T(K\times_T\pgot)\to 
K\times_T\pgot$ is the projection, and $E^{+},E^{-}$ are equivariant 
vector bundles over $K\times_T\pgot$. So for any $[k,x]\in 
K\times_T\pgot$, we have a collection of linear maps 
$\sigma([k,x,X+v]):E^{+}_{[k,x]}\to E^{-}_{[k,x]}$ depending on the 
tangent vectors $X+v$. The symbol $(i_*)^{-1}(\sigma)$ is defined 
by
\begin{equation}\label{eq.def.i.star}
(i_*)^{-1}(\sigma)(x,v)=\sigma([1,x,0+v]):E^{+}_{[1,x]}
\longrightarrow  E^{-}_{[1,x]}\quad {\rm for\ any}\quad (x,v)\in\T\pgot.
\end{equation}
For $\sigma=\Thom_{_{K}}^{\tilde{\Phi}}(\tilde{M},\tilde{J})$, the vector bundle
$E^{+}$ (resp. $E^{-}$) is $\wedge^{odd}_{\C}\T\tilde{M}$ (resp. 
$\wedge^{even}_{\C}\T\tilde{M}$). Since the complex structure leaves 
$\rgot\cong \kgot/\tgot$ and $\pgot$ invariant one gets 
$$
(i_*)^{-1}(\Thom_{_{K}}^{\underline{\lambda}}(\tilde{M},\tilde{J}))=
\Thom_{_{T}}^{\lambda}(\pgot,J_{\lambda})\, 
\wedge^{\bullet}_{\C}\kgot/\tgot\ ,
$$
and
\begin{equation}\label{eq.i.star.thom}
(i_*)^{-1}\left(
\Thom_{_{K}}^{\underline{\lambda}}(\tilde{M},\tilde{J})
\otimes \pi^*(K\times_{T}V_{\lambda})\right)=
\Thom_{_{T}}^{\lambda}(\pgot,J_{\lambda})\, V_{\lambda} 
\wedge^{\bullet}_{\C}\kgot/\tgot \ ,
\end{equation}
where $\Thom_{_{T}}^{\lambda}(\pgot,J_{\lambda})$ is the 
$T$-equivariant Thom symbol on the complex vector space
$(\pgot,J_{\lambda})$ deformed by the constant map 
$\pgot\to\tgot,\,x\mapsto\lambda$. In (\ref{eq.i.star.thom}), 
our notation uses the structure of $R(T)$-module for
$\K_{T}(\T_{T}\pgot)$, hence we can multiply 
$\Thom_{_{T}}^{\lambda}(\pgot,J_{\lambda})$ by 
$V_{\lambda} \wedge^{\bullet}_{\C}\kgot/\tgot$.

Theorem 4.1 of Atiyah in \cite{Atiyah74} tells us that 
\begin{equation}\label{diagram}
\xymatrix{
K_{T}(\T_{T}\pgot)\ar[r]^{i_{*}}\ar[d]_{\indice_{\pgot}^T} & 
K_{K}(\T_{K}\tilde{M})\ar[d]^{\indice_{\tilde{M}}^K}\\
\fgene(T)\ar[r]_{\indT} &  \fgene(K)^{K}\ .
   }
\end{equation}
is a commutative diagram, with $\tilde{M}=K\times_T\pgot$, and  
where $\indT$ is the induction map (see (\ref{eq.induction.T.K})).
In other words, $\indice_{\tilde{M}}^K(\sigma)=
\indT(\indice_{\pgot}^T((i_*)^{-1}(\sigma)))$.
With (\ref{eq.RR.tilde.2}), (\ref{eq.i.star.thom}), and  
(\ref{diagram}), we find
\begin{eqnarray*}
RR_{\Phi}^{^K}(G\cdot \lambda,V)&=&
\indT\left(\indice_{\pgot}^T(\Thom_{_T}^{\lambda}(\pgot,J_{\lambda}))
\, V_{\lambda}\wedge^{\bullet}_{\C}\kgot/\tgot\right)\\
&=&
\HolT\left(\indice_{\pgot}^T(\Thom_{_T}^{\lambda}(\pgot,J_{\lambda}))
\, V_{\lambda}\right)\ .
\end{eqnarray*}
(See the Appendix in \cite{pep4} for the relation 
$\HolT(-)=\indT(\,-\,\wedge^{\bullet}_{\C}\kgot/\tgot)$.) 
But the index $\indice_{\pgot}^T(\Thom_{_T}^{\lambda}(\pgot,J_{\lambda}))$ is 
computed in Section 5 of \cite{pep4}:
\begin{eqnarray*}
\indice_{\pgot}^T(\Thom_{_T}^{\lambda}(\pgot,J_{\lambda}))
&=&
\Big[\Pi_{\alpha\in 
\Rgot_{n}^{+}}(1-t^{-\alpha})\Big]^{-1}_{\lambda}\\
&=&
(-1)^{r}\,t^{2\rho_{n}}\Big[\Pi_{\alpha\in 
\Rgot_{n}^{+}}(1-t^{\alpha})\Big]^{-1}_{\lambda}\ ,
\end{eqnarray*}
with $r=\frac{1}{2}\dim(G/K)$. Equality (\ref{eq.theta=RR}) in then 
proved.


\subsection{\protect $(G\cdot \lambda,\Phi)$ satisfies Assumption 
\ref{hypothese.phi.t.carre}} \label{subsec.hyp.phi.2}

Let $M$ be a regular elliptic coadjoint orbit for $G$, with the 
canonical Hamiltonian $K$-action. The goal of this section is to show 
that $M$ satisfies Assumption \ref{hypothese.phi.t.carre} at every 
$\mu$. 

Let $\ggot=\kgot\oplus\pgot$ be the Cartan decomposition of $\ggot$. 
The Killing form $B$ provides a $G$-equivariant identification 
$\ggot\simeq\ggot^{*}$ and $K$-equivariant identifications  
$\kgot\simeq\kgot^{*},\,\pgot\simeq\pgot^{*}$. The Killing form $B$ 
provides also a $K$-invariant Euclidean structure on $\ggot$ such that 
$B(X,X)=-\parallel X_{1}\parallel^{2}+\parallel X_{2}\parallel^{2}$ 
and $\parallel X\parallel^{2}= \parallel X_{1}\parallel^{2}+
\parallel X_{2}\parallel^{2}$, for $X=X_{1}+X_{2}$, with $X_{1}\in\kgot,\,X_{2}\in\pgot$.
 
Hence we can and we 
shall consider $M$ as a adjoint orbit of $G$: $M=G\cdot\lambda$ 
where $\lambda\in\kgot$ is a regular element, i.e. $G_{\lambda}=K_{\lambda}$ 
is a maximal torus in $K$ (in this section $\cdot$ means the adjoint action). 
The moment map $\Phi:M\to\kgot$ is then the restriction to $M$ of the orthogonal 
projection $\kgot\oplus\pgot\to\kgot$. For $\mu\in\kgot$, we consider the map 
$\Phi_{\mu} :M\times K\cdot\mu\to\kgot,\,(m,n)\mapsto\Phi(m)-n$.

This section is devoted to the proof of the following

\medskip

\begin{prop}\label{prop.Phi.mu.compact.1}
The set $\Cr(\parallel\Phi_{\mu}\parallel^{2})$ of critical points of 
$\parallel\Phi_{\mu}\parallel^{2}$ is a compact subset of $M\times K\cdot\mu$. 
More precisely, for any $r\geq 0$, there exists $c(r)>0$ such that
 $$
 \Cr(\parallel\Phi_{\mu}\parallel^{2})\subset
 \Big(M\cap\Big\{\xi\in\ggot,\, 
 \parallel\xi\parallel\leq c(r)\Big\}\Big)\times K\cdot\mu\ ,
 $$
 whenever $\parallel\mu\parallel\leq r$.
\end{prop} 
Note that $M\cap\{\xi\in\ggot,\,\parallel\xi\parallel\leq 
c(r)\}$ is compact, thus Proposition \ref{prop.Phi.mu.compact.1}
shows that $M$ satisfies Assumption \ref{hypothese.phi.t.carre} at every 
$\mu$. Let $\agot$ be a maximal Abelian subalgebra of $\pgot$, and 
consider the map
$$
F^{\mu}:(K\cdot\lambda)\times(K\cdot\mu)\times \agot\to \R
$$
defined by $F^{\mu}(m,n,X)=\frac{1}{2}\parallel e^{X}\cdot m\parallel^{2}-
2<e^{X}\cdot m,n>$.  

\medskip

\begin{prop}\label{prop.Phi.mu.compact.2}
    For any $r\geq 0$, there exists $c(r)>0$ such that
 $$
 (m,n,X)\in\Cr(F^{\mu})\ \Longrightarrow\  
 \parallel e^{X}\cdot m\parallel\leq c(r) \ ,
 $$
 whenever $\parallel\mu\parallel\leq r$.
\end{prop} 

We first show that Proposition \ref{prop.Phi.mu.compact.2} implies 
Proposition \ref{prop.Phi.mu.compact.1}, and then we concentrate on the 
proof of Proposition \ref{prop.Phi.mu.compact.2}.

\medskip

{\em Proposition \ref{prop.Phi.mu.compact.2}\; $\Longrightarrow$ 
\ Proposition \ref{prop.Phi.mu.compact.1}.}  
Consider the map $\Phi-\mu:M\to\kgot$. \break  
One easily sees that  $\Cr(\parallel\Phi_{\mu}\parallel^{2})= K\cdot
(\Cr(\parallel\Phi_{\mu}\parallel^{2})\cap(M\times\{\mu\}))$, 
and \break 
$\Cr(\parallel\Phi_{\mu}\parallel^{2})\cap(M\times\{\mu\})
\subset \Cr(\parallel\Phi-\mu\parallel^{2})\times\{\mu\}$. Thus 
Proposition \ref{prop.Phi.mu.compact.1} is proved if one shows that 
for any $r\geq 0$, there exists $c(r)>0$ such that
$$
 \Cr(\parallel\Phi-\mu\parallel^{2})\subset
 M\cap\Big\{\xi\in\ggot,\,\parallel\xi\parallel\leq c(r)\Big\}\ ,
$$
whenever $\parallel\mu\parallel\leq r$.
Since the bilinear form $B$ is $G$-invariant, the map $m\to B(m,m)$ is 
constant on $M$, equal to $-\parallel\lambda\parallel^{2}$, and thus 
$\parallel\Phi(m)\parallel^{2}=\frac{1}{2}\parallel m\parallel^{2}
+\frac{1}{2}\parallel\lambda\parallel^{2}$ for any $m\in M$. Finally 
we have on $M$ the equality
$$
\parallel\Phi(m)-\mu\parallel^{2}=\frac{1}{2}\parallel m\parallel^{2}-
2<m,\mu> + {\rm cst}
$$
where ${\rm cst}=\frac{1}{2}\parallel\lambda\parallel^{2} +
\parallel\mu\parallel^{2}$. If we use now the Cartan decomposition 
$G=K\cdot\exp(\pgot)$, and the fact that $\pgot=\cup_{k\in K}k\cdot\agot$, 
we see that every element $M$ is of the form $m=(k_1^{-1} e^X k_2)\cdot\lambda$ 
with $k_1,k_2\in K$ and $X\in \agot$. It follows that \break 
$\parallel\Phi((k_1^{-1} e^X k_2)\cdot\lambda)-\mu\parallel^{2}=
F^{\mu}(m',n,X)+ {\rm cst}$ with $m'=k_2\cdot\lambda$, $n=k_1\cdot\mu$. It is 
now obvious that if $m=(k_1^{-1} e^X k_2)\cdot\lambda\in 
\Cr(\parallel\Phi-\mu\parallel^{2})$ then 
$(k_2\cdot\lambda,k_1\cdot\mu,X)\in \Cr(F^{\mu})$. Finally, if
Proposition \ref{prop.Phi.mu.compact.2} holds we get 
$\parallel m\parallel=\parallel e^{X}\cdot 
m'\parallel\leq c(r)$. $\Box$

\medskip

{\em Proof of Proposition \ref{prop.Phi.mu.compact.2}.}   
Let $(m,n,X)\in \Cr(F^{\mu})$. Then, the identity \break
$\frac{d}{dt}F^{\mu}(m,n,X+tX)\vert_{t=0}=0$ gives 
\begin{equation}\label{eq.A.B}
<e^{X}\cdot m,e^{X}\cdot[X,m]>=2<e^{X}\cdot [X,m],n>\ .
\end{equation}

The proof of Proposition \ref{prop.Phi.mu.compact.2} is then reduced 
to the 

\begin{lem}
  \begin{itemize}
  \item[i)] For any $r\geq 0$, there exists $d(r)>0$ such that 
$\parallel e^{X}\cdot[X,m]\parallel\leq d(r)\parallel X\parallel\,$ 
holds for every $(m,n,X)\in K\cdot\lambda\times\kgot\times\agot$ 
satisfying (\ref{eq.A.B}) and $\parallel n\parallel\leq r$.

  \item[ii)] For any $d>0$ there exists $c >0$, such that
  for every $(m,X)\in K\cdot\lambda\times\agot$, we have
$\parallel e^{X}\cdot[X,m]\parallel\leq d\parallel X\parallel\ 
\Longrightarrow\ \parallel e^{X}\cdot m\parallel\leq c$.
   \end{itemize}
\end{lem}

{\em Proof of i).}  Let $\Sigma$ be the set of weights for the adjoint 
action of $\agot$ on $\ggot$: $\ggot=\sum_{\alpha\in\Sigma}\ggot_{\alpha}$, where 
$\ggot_{\alpha}=\{Z\in\ggot,\,[X,Z]=\alpha(X)Z\ {\rm for \ all\ 
}X\in\agot\}$. Each $m\in\ggot$ amits a decomposition 
$m=\sum_{\alpha}m_{\alpha}$, with $m_{\alpha}\in\ggot_{\alpha}$, which 
is stable relatively to the Cartan involution:
\begin{equation}\label{theta-stable}
\theta(m_{\alpha})=m_{-\alpha}, \quad {\rm for\ every}\; m\in\ggot\ .
\end{equation}

Suppose now that $v:=(m,n,X)\in K\cdot\lambda\times\kgot\times\agot$ 
satisfies (\ref{eq.A.B}). We decompose $m\in K\cdot\lambda$ into 
$m=\sum_{\alpha}m_{\alpha}$ with $m_{\alpha}\in\ggot_{\alpha}$. Let 
$\Sigma^{\pm}_{v}:=\{\alpha\in\Sigma,\, m_{\alpha}\neq 0\ {\rm and}\ \pm \alpha(X)>0\}$.
The LHS of (\ref{eq.A.B}) decomposes in ${\rm LHS}=\sum_{\alpha}e^{2\alpha(X)}\alpha(X)
\parallel m_{\alpha}\parallel^{2}$, and  
\begin{eqnarray*}
{\rm LHS} &=& \sum_{\alpha\in\Sigma^{+}_{v}}e^{2\alpha(X)}\alpha(X)
\parallel m_{\alpha}\parallel^{2}  +  
\sum_{\alpha\in\Sigma^{-}_{v}}e^{2\alpha(X)}\alpha(X)
\parallel m_{\alpha}\parallel^{2}  \\
&\geq& 
\sum_{\alpha\in\Sigma^{+}_{v}}
e^{2\alpha(X)}\frac{\alpha(X)^{2}}{R\parallel X\parallel}
\parallel m_{\alpha}\parallel^{2}- R \parallel X\parallel
\sum_{\alpha\in\Sigma^{-}_{v}}
\parallel m_{\alpha}\parallel^{2}\quad   [1]
\end{eqnarray*}
with $R:=\sup_{\alpha,\parallel X\parallel\leq 
1}\vert\alpha(X)\vert$. But 
\begin{eqnarray*}
    \sum_{\alpha\in\Sigma^{+}_{v}}
e^{2\alpha(X)}\alpha(X)^{2}\parallel m_{\alpha}\parallel^{2}
&=&\parallel e^{X}\cdot[X,m]\parallel^{2}
-\sum_{\alpha\in\Sigma^{-}_{v}}
e^{2\alpha(X)}\alpha(X)^{2}
\parallel m_{\alpha}\parallel^{2}\\
&\geq&
\parallel e^{X}\cdot[X,m]\parallel^{2}- R^{2}\parallel X\parallel^{2}
\sum_{\alpha\in\Sigma^{-}_{v}}\parallel 
m_{\alpha}\parallel^{2}\quad [2]
\end{eqnarray*}
Since $\alpha\in\Sigma^{+}_{v}\Leftrightarrow -\alpha\in\Sigma^{-}_{v}$,
we have $2\sum_{\alpha\in\Sigma^{-}_{v}}\parallel 
m_{\alpha}\parallel^{2}\leq\sum_{\alpha\in\Sigma}\parallel 
m_{\alpha}\parallel^{2}=\parallel m\parallel^{2}$\break
$=\parallel \lambda\parallel^{2}$. So, the inequalities $[1]$ and $[2]$ give
\begin{equation}\label{eq.A.plusgrand}
{\rm LHS}\geq \frac{\parallel e^{X}\cdot[X,m]\parallel^{2}}{R\parallel 
X\parallel}- R\parallel X\parallel. 
\parallel \lambda\parallel^2\ .
\end{equation}
Since the RHS of (\ref{eq.A.B}) satisfies obviously 
${\rm RHS}\leq 2 \parallel e^{X}\cdot[X,m]\parallel.
\parallel n\parallel$, (\ref{eq.A.B}) and (\ref{eq.A.plusgrand}) yield
$$
2 \parallel e^{X}\cdot[X,m]\parallel.\parallel n\parallel\ \geq\ 
\frac{\parallel e^{X}\cdot[X,m]\parallel^{2}}{R\parallel 
X\parallel}- R\parallel X\parallel. 
\parallel \lambda\parallel^2\ .
$$
In other words $E:=\parallel e^{X}\cdot[X,m]\parallel$ satisfies the
polynomial inequality $E^{2}-2aE-b^2\leq 0$, with $b= R
\parallel X\parallel. \parallel \lambda\parallel$ and 
$a=R\parallel X\parallel. \parallel n\parallel$. A direct
computation gives
$$
\parallel e^{X}\cdot[X,m]\parallel\leq d\parallel X\parallel \ ,
$$
with $d=R(\parallel n\parallel+
\sqrt{\parallel n\parallel^{2}+\parallel \lambda\parallel^{2}})$. 
$\Box$

\medskip

{\em Proof of ii).}  Suppose that $ii)$ does not hold. So there is a sequence 
$(m_{i},X_{i})_{i\in\N}$ in $K\cdot\lambda\times\agot$ such that
$\parallel e^{X_{i}}\cdot[X_{i},m_{i}]\parallel\leq d
\parallel X_{i}\parallel$ but 
$\lim_{i\to\infty}\parallel e^{X_{i}}\cdot m_{i}\parallel=\infty$. We write 
$X_{i}=t_{i}v_{i}$ with $t_{i}\geq 0$ and $\parallel v_{i}\parallel=1$. We 
can assume moreover that $v_{i}\to v_{\infty}\in\agot$ with 
$\parallel v_{\infty}\parallel=1$, and 
$m_{i}\to m_{\infty}\in K\cdot\lambda$ when $i\to\infty$. 

But $m_{\infty}\in K\cdot\lambda$ is a regular element of $G$, and 
rank$(G)=$rank$(K)$, thus $[v_{\infty},m_{\infty}]=\sum_{\alpha\in \Sigma}
\alpha(v_{\infty}) m_{\infty,\alpha}\neq 0$: 
there exists $\alpha_o\in \Sigma$ such that 
$\alpha_o(v_{\infty}) m_{\infty,\alpha_o}\neq 0$, and then also  
$\alpha_o(v_{\infty}) m_{\infty,-\alpha_o}\neq 0$ (see 
(\ref{theta-stable})).

On one hand the sequence $e^{t_{i}v_{i}}.m_{i}=
\sum_{\alpha}e^{t_{i}\alpha(v_{i})}m_{i,\alpha}$ diverges. Hence 
$(t_{i})_{i\in\N}$ is not bounded and so can be assumed to be divergent. 
On the other hand \break $e^{t_{i}v_{i}}\cdot[v_{i},m_{i}]=\sum_{\alpha} 
e^{t_{i}\alpha(v_{i})}\alpha(v_{i})m_{i,\alpha}$ is bounded, so the 
sequences \break 
$e^{t_{i}\pm\alpha_o(v_{i})}\alpha(v_{i})m_{i,\pm\alpha_o}$ are 
also bounded. But $\lim_{i\to\infty}\alpha(v_{i})m_{i,\pm\alpha_o}=
\alpha_o(v_{\infty}) m_{\infty,\pm\alpha_o}$  
$\neq 0$, hence the sequences 
$e^{t_{i}\pm\alpha_o(v_{i})}$ are bounded. This contradicts 
the fact that $\lim_{i\to\infty}t_i=+\infty$ and 
$\lim_{i\to\infty}\alpha_o(v_{i})\neq 0$. $\Box$



{\small

\begin{thebibliography}{99}

\bibitem{Atiyah74} {\sc M.F. Atiyah}, Elliptic operators and 
compact groups, Springer, 1974. Lecture notes in Mathematics,  
{\bf 401}.

\bibitem{Atiyah.82} {\sc M.F. Atiyah}, Convexity and commuting
hamiltonians, {\em Bull. London Math. Soc.} {\bf 14}, 1982, p. 1-15. 

\bibitem{Atiyah-Bott-Shapiro} {\sc M.F. Atiyah},  {\sc R. Bott} and 
{\sc A. Shapiro}, Clifford modules, {\em Topology} {\bf  3}, Suppl. 1, 
p. 3-38. 

\bibitem{Atiyah-Segal68} {\sc M.F. Atiyah}, {\sc G.B. Segal},
The index of elliptic operators II, {\em Ann. Math.} {\bf 87}, 1968,
p. 531-545.

\bibitem{Atiyah-Singer-1} {\sc M.F. Atiyah}, {\sc I.M. Singer},
The index of elliptic operators I, {\em Ann. Math.} {\bf 87}, 1968,
p. 484-530.

\bibitem{Atiyah-Singer-2} {\sc M.F. Atiyah}, {\sc I.M. Singer},
The index of elliptic operators III, {\em Ann. Math.} {\bf 87}, 1968,
p. 546-604.

\bibitem{Atiyah-Singer-3} {\sc M.F. Atiyah}, {\sc I.M. Singer},
The index of elliptic operators IV, {\em Ann. Math.} {\bf 93}, 1971,
p. 139-141.




\bibitem{B-G-V} {\sc N. Berline}, {\sc E. Getzler} and {\sc M. Vergne}, 
{\em Heat kernels and Dirac operators}, Grundlehren, vol. 298, Springer, 
Berlin, 1991.

\bibitem{B-V.inventiones.96.1} {\sc N. Berline} and {\sc M. Vergne}, 
The Chern character of a transversally elliptic symbol and the 
equivariant index, {\em Invent. Math.} {\bf 124}, 1996, p. 11-49.

\bibitem{B-V.inventiones.96.2} {\sc N. Berline} and {\sc M. Vergne}, 
L'indice \'equivariant des op\'erateurs transversalement elliptiques, 
{\em Invent. Math.} {\bf 124}, 1996, p. 51-101.



\bibitem{CdS-K-T}{\sc A. Cannas da Silva}, {\sc Y. Karshon} and {\sc S. Tolman}, 
Quantization of presymplectic manifolds and circle actions, 
{\em Trans. Amer. Math. Soc.} {\bf 352}, 2000, p. 525-552.


\bibitem{Duflo77} {\sc M. Duflo}, Repr\'esentations de carr\'e int\'egrable des 
groupes semi-simples r\'eels, {\em S\'eminaire Bourbaki}, 
Vol. 1977/78, Expos\'e No.508, Lect. Notes Math. 710, 22-40 (1979). 

\bibitem{Duflo-Heckman-Vergne} {\sc M. Duflo}, {\sc G. Heckman} and {\sc 
M. Vergne}, Projection d'orbites, formule de Kirillov et 
formule de Blattner, {\em M\'emoires de la S.M.F.} {\bf 15}, 1984, p. 
65-128.

\bibitem{Duistermaat96} {\sc J. J. Duistermaat}, 
{\em The heat lefschetz fixed point formula for the 
Spin\textsuperscript{c}-Dirac operator}, 
Progress in Nonlinear Differential Equation and Their Applications, 
vol. 18, Birkhauser, Boston, 1996.




\bibitem{Guillemin-Sternberg82.bis} {\sc V. Guillemin} and {\sc S. 
Sternberg}, Convexity properties of the moment mapping, 
{\em Invent. Math.} {\bf 67}, 1982, p. 491-513.

\bibitem{Guillemin-Sternberg82} {\sc V. Guillemin} and {\sc S. 
Sternberg}, Geometric quantization and multiplicities of group 
representations, {\em Invent. Math.} {\bf 67}, 1982, p. 515-538.

\bibitem{Guillemin-Sternberg84} {\sc V. Guillemin} and {\sc S. 
Sternberg}, A normal form for the moment map, 
in {\em Differential Geometric Methods in Mathematical Physics}(S. 
Sternberg, ed.), Reidel Publishing Company, Dordrecht, 1984.


\bibitem{Jeffrey-Kirwan97} {\sc L. Jeffrey} and {\sc F. Kirwan}, 
Localization and quantization conjecture, {\em Topology} {\bf 36}, 
1997, p. 647-693.



\bibitem{Harish-Chandra65-66} {\sc Harish-Chandra} Discrete series 
for semi-simple Lie group, I and II, {\em Acta Mathematica}, 
{\bf 113} (1965) p. 242-318, and {\bf 116} (1966) p. 1-111.

\bibitem{Hecht-Schmid} {\sc H. Hecht} and {\sc W. Schmid}, A 
proof of Blattner's conjecture, {\em Invent. Math.}, {\bf 31}, 1975, 
p. 129-154.

\bibitem{Kawasaki81} {\sc T.Kawasaki}, The index of elliptic operators 
over V-manifolds, {\em Nagoya Math. J.}, {\bf 84}, 1981, p. 135-157.

\bibitem{Kirwan84} {\sc F. Kirwan}, {\em Cohomology of quotients in 
symplectic and algebraic geometry}, Princeton Univ. Press, Princeton, 1984.

\bibitem{Kirwan.84.bis} {\sc F. Kirwan}, Convexity properties of the 
moment mapping III, {\em Invent. Math.} {\bf 77}, 1984, p. 547-552.


\bibitem{Kostant70} {\sc B. Kostant}, Quantization and unitary 
representations, in {\em Modern Analysis and 
Applications}, Lecture Notes in Math., Vol. 170, Spinger-Verlag, 1970,
p. 87-207.

\bibitem{Lawson-Michel} {\sc H. Lawson} and {\sc M.-L. Michelsohn}, 
{\em Spin geometry}, Princeton Math. Series, 38. Princeton Univ. Press,
Princeton, 1989.


\bibitem{L-M-T-W} {\sc E. Lerman}, {\sc E. Meinrenken}, {\sc S. Tolman} 
and {\sc C. Woodward}, Non-Abelian convexity by symplectic cuts, 
{\em Topology} {\bf 37}, 1998, p. 
245-259.

\bibitem{Meinrenken96} {\sc E. Meinrenken}, 
On Riemaan-Roch formulas for multiplicities,
{\em J. Amer. Math. Soc.}, {\bf 9}, 1996, p. 373-389.

\bibitem{Meinrenken98} {\sc E. Meinrenken}, 
Symplectic surgery and the Spin\textsuperscript{c}-Dirac operator, 
{\em Adv. in Math.} {\bf 134}, 1998, p. 240-277.

\bibitem{Meinrenken-Sjamaar} {\sc E. Meinrenken}, {\sc S. Sjamaar},
Singular reduction and quantization, {\em Topology} {\bf 38}, 1999, 
p. 699-762. 

\bibitem{pep1} {\sc P-E. Paradan}, Formules de localisation en 
cohomologie \'equivariante, {\em Compositio Mathematica}  
{\bf 117}, 1999, p. 243-293.


\bibitem{pep2} {\sc P-E. Paradan}, The moment map and equivariant 
cohomology with generalized coefficient, {\em Topology}, 
{\bf 39}, 2000, p. 401-444. 

\bibitem{pep3} {\sc P-E. Paradan}, The Fourier transform of 
semi-simple coadjoint orbits, {\em J. Funct. Anal.}  {\bf 163}, 1999, 
p. 152-179.

\bibitem{pep4} {\sc P-E. Paradan}, Localization of the Riemann-Roch 
character, {\em J. Funct. Anal.} {\bf 187}, 2001, p. 442-509..


\bibitem{Schmid71} {\sc W. Schmid}, On a conjecture of Langlands, 
{\em Ann. of Math.} {\bf 93}, 1971, p. 1-42.

\bibitem{Schmid76}{\sc W. Schmid}, $L\sp{2}$-cohomology and the
  discrete series, {\em Ann. of Math.} {\bf 103}, 1976, p. 375-394.

\bibitem{Schmid97} {\sc W. Schmid}, Discrete Series, {\em Proc. of Symp. 
in Pure Math.} {\bf 61}, 1997, p. 83-113.


\bibitem{Sjamaar96} {\sc R. Sjamaar}, Symplectic reduction and Riemann-Roch formulas for 
multiplicities, {\em Bull. Amer. Math. Soc.} {\bf 33}, 1996,
p. 327-338.

\bibitem{Sjamaar98} {\sc R. Sjamaar}, Convexity properties of the moment mapping 
re-examined, {\em Adv. Math.}  {\bf 138}, 1998, p. 46-91. 

\bibitem{Segal68} {\sc G. Segal}, Equivariant K-Theory, {\em Publ.
Math. IHES} {\bf 34}, 1968, p. 129-151.

\bibitem{Tian-Zhang98} {\sc Y. Tian}, {\sc W. Zhang}, An analytic
proof of the geometric quantization conjecture of
Guillemin-Sternberg, {\em Invent. Math.} {\bf 132}, 1998, p. 229-259.


\bibitem{Vergne94} {\sc M. Vergne}, Geometric quantization and
equivariant cohomology, First European Congress in Mathematics, vol. 
1, {\em Progress in Mathematics} {\bf 119}, Birkhauser, Boston, 1994, 
p. 249-295.

\bibitem{Vergne96} {\sc M. Vergne}, Multiplicity formula for geometric
quantization, Part I, Part II, and Part III, 
{\em Duke Math. Journal}  {\bf 82}, 1996, 
p. 143-179, p 181-194, p 637-652.

\bibitem{Vergne01} {\sc M. Vergne}, Quantification g\'eom\'etrique et 
r\'eduction symplectique,  {\em S\'eminaire Bourbaki}  {\bf 888}, 2001. 

\bibitem{Witten} {\sc E. Witten}, Two dimensional gauge theories 
revisited, {\em J. Geom. Phys.} {\bf 9}, 1992, p. 303-368.

\bibitem{Woodhouse} {\sc Woodhouse}, {\em Geometric quantization}, 2nd ed. 
Oxford Mathematical Monographs. Oxford: Clarendon Press, 1997.

\end{thebibliography}
}
\end{document}